\documentclass[11pt]{article}

\usepackage{tabularx}
\usepackage{amssymb,amsmath,natbib}
\usepackage{amsthm}
\usepackage{enumerate}
\usepackage{multirow}
\usepackage{graphicx}
\usepackage{subfig}
\usepackage{caption}
\usepackage{float}
\usepackage{graphics,color}
\usepackage{bm}
\usepackage{secdot}
\usepackage{multibib}
\usepackage{url}
\usepackage{hyperref}

\newcites{ECtex}{Electronic Companion References}

\newcommand{\eProof}{\hspace*{.1in}\hfill  $\blacksquare$}
\numberwithin{equation}{section}

\numberwithin{figure}{section}
\numberwithin{table}{section}

\newtheorem{theorem}{Theorem}
\newtheorem{corollary}{Corollary}
\newtheorem{lemma}{Lemma}

\newtheorem{proposition}{Proposition}
\newtheorem{definition}{Definition}
\newtheorem{remark}{Remark}
\newtheorem{claim}{Claim}
\DeclareMathOperator*{\argmin}{arg\,min}

\begin{document}
%%%%%%%%%%%%%%%%

%\documentclass[11pt]{article}
%\usepackage{amssymb,amsmath,latexsym,fullpage,tabularx}
%\usepackage{epsfig}

%\input{amsvect.tex}
%\numberwithin{equation}
%\renewcommand{\theequation}{\thechapter.\arabic{equation}}
%\renewcommand{\thefigure}{\arabic{figure}}
%\renewcommand{\thetable}{\arabic{table}}
%\renewcommand{\baselinestretch}{2.5}
%\parskip 0.01in
%\setcounter{equation}{0}
%\setcounter{chapter}{0}
%\newcommand{\qed}{\hfill $\fbox{}$}  \def\M{\hspace*{1.5em}}

\newcommand{\ls}[1]
{\dimen0=\fontdimen6\the\font \lineskip=#1\dimen0
\advance\lineskip.5\fontdimen5\the\font \advance\lineskip-\dimen0
\lineskiplimit=.9\lineskip \baselineskip=\lineskip
\advance\baselineskip\dimen0 \normallineskip\lineskip
\normallineskiplimit\lineskiplimit
\normalbaselineskip\baselineskip \ignorespaces }

\begin{center}
  {\LARGE Staffing Call Centers with Uncertain Arrival Rates and Co-sourcing }\\[12pt]
\end{center}

\begin{center}
{\textbf {Ya\c{s}ar Levent Ko\c{c}a\u{g}a} \\
{Syms School of Business, Yeshiva
University, Belfer Hall 403/A, New York, NY 10033, kocaga@yu.edu}\\
\bigskip
\textbf{Mor Armony}\\
{Stern School of Business, New York University, KMC 8--62, New York, NY, 10012, marmony@stern.nyu.edu}\\
\bigskip
\textbf{Amy R. Ward}\\
Marshall School of Business, University of Southern California, Bridge Hall 401H Los Angeles, CA, 90089, amyward@marshall.usc.edu }
\end{center}

% Here is the abstract:
\noindent\hrulefill

\noindent
In a call center, staffing decisions must be made before the call arrival rate is known with certainty.  Then, once the arrival rate becomes known, the call center may be over-staffed, in which case staff are being paid to be idle, or under-staffed, in which case many callers hang-up in the face of long wait times.  Firms that have chosen to keep their call center operations in-house can mitigate this problem by co-sourcing; that is, by sometimes outsourcing calls.  Then, the required staffing $N$ depends on how the firm chooses which calls to outsource in real-time, after the arrival rate realizes and the call center operates as a $M/M/N+M$ queue with an outsourcing option. Our objective is to find a {\em joint} policy for staffing and call outsourcing that minimizes the long run average cost of this two-stage stochastic program when there is a linear staffing cost per unit time and linear costs associated with abandonments and outsourcing.

We propose a policy that uses a square-root safety staffing rule, and outsources calls in accordance with a threshold rule that characterizes when the system is ``too crowded''.  Analytically, we establish that our proposed policy is asymptotically optimal, as the mean arrival rate becomes large, when the level of uncertainty in the arrival rate is of the same order as the inherent system fluctuations in the number of waiting customers for a known arrival rate. Through an extensive numerical study, we establish that our policy is extremely robust.  In particular, our policy performs remarkably well over a wide range of parameters, and far beyond where it is proved to be asymptotically optimal.

\smallskip


\baselineskip 20pt plus .3pt minus .1pt

\noindent\hrulefill

%\HISTORY{}

%    \thispagestyle{empty}
%\end{titlepage}}

%\maketitle

%\ls{2}

\section{Introduction}
Call centers have become ubiquitous in business. Today, every
Fortune 500 company has at least one call center, and the average
Fortune 500 company employs 4500 call center agents (who may be
distributed across more than one site) \citep{GiKh2005}. For many companies, the call center is a primary point-of-contact
with their customers. Hence a well-run call center
promotes good customer relations, and a poorly-managed one hurts
them. But call center management is difficult.

A call
center manager faces the classical operational challenge of
determining appropriate staffing levels throughout the day and week
in order to meet a random and time-varying call volume. This is extremely difficult especially when call arrival {\em rate} is itself {\em random}, as was empirically shown in \cite{Brown05} and \cite{Maman09}, among others.
When
staffing levels are too low, customers are put on hold, and many
hang up in frustration while waiting for an agent to take their
call. But when staffing levels are too high, the call center manager
ends up paying staff to be idle.

One option in managing this uncertainty is for a company to outsource its call center operations. Then, the
challenges of call center management can be handled by a vendor firm whose
primary focus is call center operations. That vendor can pool demand
amongst various companies, thereby lowering variability, which
should allow for more accurate demand forecasts, and so better
staffing decisions. However, it is also true that many companies are reluctant to relinquish control of their call center operations.  This is evidenced by a recent
survey from the Incoming Call Management Institute~\citep{ICMI2006}: only 7.9\% out of 279 call center professionals
 used an outside vendor to handle most or
all of their calls. One reason for that are the ``hidden costs of outsourcing'' \citep{Kharif2003}, which include service quality costs that are hard to explicitly quantify.  As a result, many of these companies prefer to co-source; that is, to outsource some, but not all, of their calls.

We study a co-sourcing structure in which the vendor charges the company a fee per call outsourced, which is consistent with the pay-per-call (PPC) co-sourcing structure analyzed in~\cite{Aksin2008}.   Then, the company can decide on a call-by-call basis which calls to answer in-house and which calls to route to the vendor.  This is helpful because call centers typically make their daily staffing decisions at least a week in advance, before the actual arrival rate to the call center for a given day is known.  If the planned staffing is sufficient to handle the mean arrival rate, then the company needs to outsource only a small fraction of calls in order to handle the inherent variability that results in congestion every so often.  On the other hand, if the planned staffing is insufficient to handle the mean arrival rate, then the company can outsource a large fraction of its calls, thereby preventing high congestion levels.

The relevant question for this paper is:  how do we decide on staffing levels when the arrival rate is uncertain and the aforementioned co-sourcing option is present?  To answer this question, we begin with one simple and widely used queueing model of a call center, the Erlang A or $M/M/N+M$ (see, for example, Section 4.2.2 in~\cite{Gans03}), and add uncertainty in the arrival rate and an outsourcing option.
Then, our model for the call center is a multi-server queue with a doubly stochastic time-homogeneous Poisson arrival process, exponential service times, and exponential times to abandonment.  Although this model ignores the time-varying nature of the arrival rate over the course of each day, there is call center literature that discusses how to use the Erlang A model to make staffing decisions for time-varying arrival rates, using the ``stationary independent period by period approach'' (SIPP); see~\cite{Green_etal_2001},~\cite{Gans03},~\cite{Aksin07}, and~\cite{LiuWhitt2012} for more discussion.  We suppose that a similar approach can be adopted here to accommodate the added feature of arrival rate uncertainty.

 The control decisions in our model are (a) an upfront staffing decision and (b) real-time call outsourcing (routing) decisions.
Recall that staffing decisions are made on a much longer time horizon  and well before the timing of the control decisions.  In particular, these decisions are made on two different time scales.  This means that we have a two stage stochastic program:  The staffing decisions are made in the first stage, before the arrival rate is known, and the outsourcing decisions are made in the second stage, after the arrival rate is known.  Then, the outsourcing decisions can depend on the actual arrival rate even though the staffing decisions cannot.

 Our objective is to propose a policy for staffing and outsourcing under the assumption of linear staffing cost and linear abandonment and outsourcing costs.  Then, for each arriving customer that cannot be immediately served, there is a tension between choosing to outsource that customer (and paying the outsourcing fee) or having that customer wait for an in-house agent (and risking incurring an abandonment cost). In summary, we are solving a joint staffing and routing control problem for a (modified) Erlang A model with an uncertain arrival rate and an outsourcing option.

The three main contributions in this paper are:

\begin{itemize}
 \item  The {\em modeling} contribution is the formulation of a joint staffing and outsourcing problem for a call center that has access to co-sourcing, and must make staffing decisions when there is arrival rate uncertainty. This modeling framework can be used to study more general joint staffing and control problems in call centers that have been previously studied in the literature under the assumption of a known arrival rate.
 \item The {\em application} contribution is the development of a square root safety staffing and threshold outsourcing policy that we numerically show to be extremely {\em robust} over the entire parameter space. This robustness may come as no surprise for readers who are familiar with related literature such as~\cite{Borst04} and~\cite{GHM2012}. However, the existing literature has not addressed the issue of robustness in the context of random arrival rates and dynamic control, nor can this robustness be readily explained using existing results.
 \item The {\em technical} contribution is the proof that our proposed square root safety staffing and threshold outsourcing policy is asymptotically optimal, as the mean arrival rate becomes large, when the level of uncertainty in the arrival rate is of the same order as the inherent system stochasticity (which is of the order of the square root of the mean of the arrival rate).
\end{itemize}

The remainder of this paper is organized as follows.  First, we review the most relevant literature.  Next, in Section \ref{Section: Model}, we describe our model in detail.  In Section \ref{Section:Exact}, we present the exact (non-asymptotic) analysis which leads to an algorithm to compute the optimal policy numerically.  However, that algorithm does not provide insight into the structure of an optimal policy, and so, in  Section \ref{Section:aoHW}, we perform an asymptotic analysis under the assumption that the level of uncertainty in the arrival rate is of the same order as the inherent system stochasticity.  That asymptotic analysis motivates us to propose, in Section~\ref{Section:ProposedPolicy}, a square root safety staffing and threshold outsourcing policy that is universal in the sense that there is no assumption on the level of uncertainty in the arrival rate.
We evaluate the performance of our universal policy numerically in Section \ref{Section:Numerics}. We make concluding remarks in Section \ref{Section: Conclusion}. All proofs and additional numerical results can be found in the electronic companion (EC).

\noindent{\bf{Literature Review}}

Previous work on joint staffing and routing problems in call centers includes \cite{GurvichArmonyMandelbaum08} who study staffing and dynamic routing in call centers with multiple customer classes and a single server pool, \cite{ArmonyMandelbaum11} who consider the symmetric case of a single customer type and a heterogeneous server pool, and \cite{GurvichWhitt10} who consider multiple customer classes and a heterogeneous server pool.
These papers study the staffing and dynamic routing problems within the Halfin-Whitt heavy traffic regime, pioneered by \cite{HalfinWhitt81}, and extended to include abandonments by \cite{GarnettMandelbaumReiman02}. This is also known as the Quality and Efficiency Driven (QED) heavy traffic regime. The key idea is to approximate the behavior of call centers that are modeled as multi-server queues with that of their limiting diffusions. The limiting diffusion arises from a specific relationship between the arrival rate and the staffing level as both grow large without bound. Our work is different than the aforementioned papers in that our model is pertinent to situations where the arrival rate is not known when staffing decisions  are made, and thus has to be inferred or forecasted using available historical data.

Given a staffing level and a realized arrival rate, our dynamic outsourcing decision is equivalent to the admission control problem studied in \cite{KocagaWard10}. The equivalence follows because we do not explicitly model the vendor firm and assume it has ample capacity to handle the outsourced calls. \cite{KocagaWard10} show that a threshold admission control policy is optimal, and characterize a simple form for the threshold level that is asymptotically optimal when the staffing level is {\em assumed} to be such that the system operates in the QED regime.  In contrast, this paper explicitly models the staffing decisions and has a random arrival rate.

There is a growing body of literature that studies staffing for call centers with uncertain arrival rates including (in chronological order) \cite{Chen01}, \cite{Jongbloed01}, \cite{Ross01}, \cite{BassambooHarrisonZeevi05}, \cite{BassambooHarrisonZeevi05}, \cite{Whitt06}, \cite{BassambooHarrisonZeevi06}, \cite{Steckley2009}, \cite{Maman09},  \cite{BassambooZeevi09}, \cite{GurvichLuedtkeTezcan09}, \cite{RobbinsHarrison2010}, \cite{BassambooRandhawaZeevi10}, \cite{MehrotraOzlukSaltzman2010},  \cite{Gans12}, and \cite{ZanHasenbeinMorton2013}.  The two works most closely related to ours are \cite{Maman09} and \cite{BassambooRandhawaZeevi10}, and we discuss each in turn.

The focus of \cite{Maman09} is to extend the QED staffing formula under a general form of arrival uncertainty.  Our asymptotic optimality result assumes a special case of the form of the arrival rate uncertainty presented in that paper. However,  that paper does not explicitly study the cost minimizing staffing and does not model routing decisions, as we do.

\cite{BassambooRandhawaZeevi10} propose a staffing policy for an $M/M/N+M$ queue in which the arrival rate is random, and there is no outsourcing option. They establish that a simple newsvendor based staffing policy performs extremely well when the order of uncertainty in the arrival rate exceeds the order of the inherent system stochasticity.  In contrast, we establish the asymptotic optimality of our proposed policy when the aforementioned two magnitudes are the same.  Then, in our numeric study, we adapt their policy to our setting with outsourcing, in order to evaluate policy robustness. 

In relation to the literature on call center outsourcing (see, for example, \cite{ZhouRen2011}), our paper is most similar to \cite{Aksin2008}. In contrast to most papers in this literature, which assumes all calls will be outsourced, \cite{Aksin2008} considers the contract design problem of a company which faces an uncertain call volume, and can outsource part of its calls by choosing between a capacity-based and volume-based contract that is pay-per-call. Although both \cite{Aksin2008} and our model study co-sourcing decisions which are driven by call volume uncertainty, \cite{Aksin2008} focuses on the optimal contract choice, whereas we focus on the in-house staffing and dynamic routing decisions and assume the contract.

\section{Model Description}
\label{Section: Model}
We  model the in-house call center (which we henceforth refer to as ``the call center" or ``the system")  as an $M/M/N+M$ queueing system in which the arrival rate is uncertain. We let $\Lambda$ denote the random arrival rate with cdf $F_{\Lambda}$, and we let $l$ denote a particular realization of $\Lambda$. We assume that $\Lambda$ is a non-negative random variable with mean $E\left[\Lambda\right]=\lambda$.  For ease of exposition, we assume the mean service time is 1, so that we can think of measuring time in terms of the mean time to serve an arrival.  The mean patience time is $1/\gamma$.

The call center manager must make two decisions:  the upfront staffing level $N$, and the dynamic outsourcing decision. The staffing level $N:=N(F_{\Lambda})$ must be set before the arrival rate $\Lambda$ is realized,  based on the knowledge of its distribution.  After the arrival rate $\Lambda$ is realized as $l$, every arriving call can be either accepted into the system, or routed to the outsourcing vendor.  Then, the routing control policy $\pi :=\pi(N,l)$ is in general a function of the staffing level $N$ and $l$. The notation $\pi(N,\Lambda)$ refers to a routing policy that may depend on the actual realization $l$ of $\Lambda$.  Any stationary routing control policy $\pi = (\pi_n : n \in \{0,1,\ldots\})$ is a vector, where $\pi_n \in [0,1]$ denotes the probability that a customer is accepted into the system when there are $n$ customers currently present there.  We let $\Pi$ be the set of all such vectors.  An admissible policy
\[
u:=(N, \pi(N,\Lambda)) = (N, (\pi_n(N,\Lambda): n \in \{0,1,\ldots\}))
 \]
 sets the staffing level as a non-negative integer $N$, and, after the arrival rate $\Lambda$ realizes as $l$, controls outsourcing decisions dynamically by routing calls according to the policy $\pi(N,l)\in \Pi$.

After the arrival rate $l$ realizes, the system operates as a birth and death process with birth rate $l \pi_n$ and death rates
\[
\mu_n = \min(N,n) + \gamma [n-N]^+.
\]
It is straightforward to solve the balance equations for the steady-state distribution for the number-in-system process, from which it follows that the following performance measures  are well-defined:
\begin{eqnarray*}
P_{\pi}(\mbox{ab};l) & = & \mbox{ the probability an entering customer abandons;} \\
P_{\pi}(\mbox{out};l) & = & \mbox{ the probability an arriving customer is routed to the outsourcer}.
\end{eqnarray*}

The objective of the system manager is to minimize the expected long-run average cost, when there are costs due to customer abandonment, routing customers to the outsourcing vendor, and staffing costs.  Every customer that abandons the system before receiving service costs $a$ and the per call cost of routing to the outside vendor is $p$. Then, the long-run average operating cost associated with $\pi \in \Pi$ when the arrival rate realizes as $l$ and the staffing level is $N$ is
\begin{equation} \label{eq:stage2cost}
z_\pi := z_\pi(N,l) = p l P_\pi(\mbox{out};l) + a l P_\pi(\mbox{ab};l).
\end{equation}
This is expressed as a random variable by writing each of the performance measures of interest as random variables, as follows
\[
z_\pi(N,\Lambda) = p \Lambda  P_\pi(\mbox{out};\Lambda) + a \Lambda P_\pi(\mbox{ab};\Lambda).
\]
For a given realization $l$ of $\Lambda$,
 $z^{\mbox{\tiny{opt}}}(N,l)$ denotes the minimum cost, and $\pi^{\mbox{\tiny{opt}}}(N,l) \in \Pi$ is a policy that achieves that minimum cost.
The associated random varaiables are
 $z^{\mbox{\tiny{opt}}}(N,\Lambda)$ and $\pi^{\mbox{\tiny{opt}}}(N,\Lambda)$.
The expected long-run average cost under the policy $u=(N,\pi)= (N,\pi(N,\Lambda))$, with respect to the random arrival rate $\Lambda$, is
\begin{equation} \label{eq:objective_function}
\mathcal{C}(u) = cN + E \left[ z_\pi\left(N,\Lambda\right) \right].
\end{equation}
 We would like to find a staffing level $N^{\mbox{\tiny{opt}}}$ and a routing control policy $\pi^{\mbox{\tiny{opt}}}(N,\Lambda)$ that achieves the minimum long-run average cost
\begin{equation} \label{eq:objective}
\mathcal{C}^{\mbox{\tiny{opt}}} := \inf_{u} \mathcal{C}(u) = \min_{N \in \{0,1,2,\ldots\}} cN + E\left[ z^{\mbox{\tiny{opt}}}(N,\Lambda) \right].
\end{equation}

\begin{remark} \label{remark:waitingCost} (Including waiting costs.)
The objective function in (\ref{eq:objective_function}) can be modified to include a customer waiting cost  by modifying (\ref{eq:stage2cost}) as follows.  Suppose the cost for one customer to wait one time unit is $w \geq 0$.  Then, (\ref{eq:stage2cost}) becomes
\[
z_\pi := z_\pi(N,l) = p l P_\pi(\mbox{out};l) + a l P_\pi(\mbox{ab};l) + wl(1-P_{\pi}(\mbox{out};l))\overline{W}_\pi(l),
\]
where $\overline{W}_\pi(l)$ is the steady-state average waiting time, including both abandoning and served customers.  Letting $\overline{Q}_\pi(l)$ denote the steady-state average number of customers waiting in queue to be served, it follows from Little's law that
\[
l(1-P_{\pi}(\mbox{out};l))\overline{W}_\pi(l) = \overline{Q}_\pi(l),
\]
and so
\[
z_\pi = p l P_\pi(\mbox{out};l) + a l P_\pi(\mbox{ab};l) + w \overline{Q}_\pi(l).
\]
Also, since the steady-state rate at which abandoning customers arrive must equal the steady-state abandonment rate
\[
l P_\pi(\mbox{ab};l) = \gamma \overline{Q}_\pi(l),
\]
and so
\[
z_\pi = p l P_\pi(\mbox{out};l) + \left( a +\frac{w}{\gamma} \right) l P_\pi(\mbox{ab};l).
\]
The analysis in this paper is valid with $a$ replaced by $a' := a + w/\gamma$.
Therefore, to include a customer waiting cost, the only change is to replace $a$ in (\ref{eq:stage2cost}) by $a'$.
\end{remark}

\begin{table}
\centering
\begin{tabular}{|c|c|c|c|}
\cline{3-4}
\multicolumn{2}{c}{} & \multicolumn{2}{|c|}{} \\
\multicolumn{2}{c}{} & \multicolumn{2}{|c|}{In-house} \\
\multicolumn{2}{c}{} & \multicolumn{2}{|c|}{capacity investment} \\
\cline{3-4}
\multicolumn{2}{c|}{} &  &  \\
\multicolumn{2}{c|}{} & $c<\min(a,p)$ & $c\geq \min(a,p)$ \\ \hline
&  &  &  \\
& $a>p$ & co-sourcing & complete outsourcing \\ \cline{2-4}
outsourcing &  &  &  \\
leverage & $a\leq p$ & no outsourcing & no operation \\ \hline
\end{tabular}
\caption{Optimal In-house Staffing and Outsourcing Decisions}
\label{tab:optdecisions}
\end{table}

It makes intuitive sense that if $c\geq p$ the system manager will not invest in any in-house capacity because serving an arrival is more costly than routing that arrival to the outsourcing vendor (recall that $\mu=1$, so $c$, $p$ and $a$ are comparable). Continuing with such intuitive comparisons (see Table 1), if $c \geq \min(a,p)$, then the system manager will either route every call to the outsourcer ($a>p$) or will let every call abandon ($p\geq a$).  This suggests that it is only when $c < \min(a,p)$ that the system manager will invest in in-house capacity.  Then, he will not route calls to the outsourcer if $a \leq p$.  In summary, we expect the system manager to invest in capacity and route some calls to the outsourcing vendor only if  $c< p < a$. The following proposition confirms this observation.
\begin{proposition}\label{prop: prelimit optimal} (Characterizing the parameter regimes.)
\begin{enumerate}[(i)]
\item Suppose that $c \geq \min(a,p)$.  Then, $N^{\mbox{\tiny{opt}}}=0$ solves (\ref{eq:objective}).
\begin{enumerate}
\item In addition, if $a>p$, then $\pi^{\mbox{\tiny{opt}}}(N,l) = (0,0,0,\ldots)$, so that all calls are routed to the outsourcing vendor.
\item Otherwise, if $a\leq p$, then $\pi^{\mbox{\tiny{opt}}}(N,l) = (1,1,1,\ldots)$, so that all calls are left to abandon.
\end{enumerate}
\item If $c<\min(a,p)$ and $a \leq p$, then the optimal control policy is $\pi^{\mbox{\tiny{opt}}}(N,l) = (1,1,1,\ldots)$ for any given $N\geq0$.
\end{enumerate}
\end{proposition}

\noindent From Proposition~\ref{prop: prelimit optimal}, is follows that the only cases with a non-trivial optimal staffing are when $c<\min(a,p)$.
Hence, for the remainder of the paper we {\bf assume that } $\mathbf{c<\min(a,p)}$.

%[LEVENT AND MOR:  In the remainder of the paper, we assume $c<\min(a,p)=p$ not $c<\min(a,p)$, right?  Also, I don't think we need both a paragraph before Proposition~\ref{prop: prelimit optimal} explaining it and one after.  I have changed so that all that intuition appears before.  Also, can we separate out the two cases in part (i) of Proposition~\ref{prop: prelimit optimal}?]
%{\color{blue}[AMY AND MOR:  I belive the AO result is also valid when $c<\min(a,p)$ while the numerics assume  $c<\min(a,p)=p$, right? I agree with the rest of the changes]}[LEVENT:  I agree with you.  However, it seems wierd to use the Proposition to rule out trivial cases, and then to include the trivial case.  Why do we need to prove asymptotic optimality when $c<\min(a,p)=a$ when we have just proved the EXACT optimal policy?]
%{\color{blue}[AMY: We have proved only the routing control will keep all calls in-house,  the staffing level still remains to be determined. Also, please see Remark 2.]}
\section{Exact Analysis}
\label{Section:Exact}

For a fixed staffing level $N$ and realized arrival rate $l$, the problem of minimizing $z_\pi$ is a Markov decision problem (MDP), having solution $z^{\mbox{\tiny{opt}}}(N,l)$. This MDP has been solved in \cite{KocagaWard10} in the context of an admission control problem. It follows from Theorems 3.1, Theorem 3.2, and Theorem 3.3 in \cite{KocagaWard10} that the optimal policy is a deterministic threshold policy (with a potentially infinite threshold level).  Hence we can restrict ourselves to the class of threshold control policies
\[
\tau(T) = ( \tau_n(T) : n \in \{0,1,2,\ldots\} ),
\]
defined for threshold level $T := T(N,l) \in [0,\infty]$ as
\[
\tau_{n}(T) := \left\{ \begin{array}{ll} 1 & \mbox{ if } n  < T, \\ 0 & \mbox{ if } n \geq T. \end{array} \right.
\]

Under the threshold policy $\tau_{n}(T)$, after the arrival rate realizes, the system operates as an $M$/$M$/$N$/$T+M$ queue.  Then, the process tracking the number of customers in the system is a birth-and-death process on $\{0,1,\ldots,N-1,N,N+1,\ldots,T\}$ with birth rate $l$ and death rate in state $n$
\[
\mu_n = \mbox{min}(n,N)+\gamma[n-N]^+.
\]
Then, we can solve the balance equations to find the steady-state probabilities
\[
\theta_k(l) = \left( \prod_{i=1}^k \frac{l}{\mu_i} \right) \theta_0(l)
\]
for
\[
\theta_0(l) = \frac{1}{\sum_{k=0}^T \left( \prod_{i=1}^k \frac{l}{\mu_i} \right) },
\]
and develop the expressions for the performance measures\footnote{See also Section 7 in~\cite{Whitt05} for exact expressions for other performance measures of interest, such as the expected wait time conditioned on an arrival being served, in a more general model that allows for state-dependent abandonment rates.}
\begin{eqnarray*}
P_{\tau(T)}(\mbox{out};l) & = & \theta_T(l) \\
\overline{Q}_{\tau(T)}(l) & = & \sum_{k=0}^T [k-N]^+ \theta_k(l).
\end{eqnarray*}
Then, recalling from Remark~\ref{remark:waitingCost} that $lP_\pi(\mbox{ab};l)=\gamma \overline{Q}_\pi(l)$,
 we can express the long-run average cost in terms of the steady-state probabilities as
\[
z_{\tau(T)} = p l \theta_T(l) + a\gamma \overline{Q}_{\tau(T)}(l).
\]
Hence we can optimize over $T$ to find
\[
T^{\mbox{\tiny{opt}}} := \argmin_{T \in \{0,1,2,\ldots\}} z_{\tau(T)},
\]
for which
\[
 z_{\tau(T^{\mbox{\tiny{opt}}})}(N,l)= z^{\mbox{\tiny{opt}}}(N,l).
\]

Unfortunately, the resulting expression for $z_{\tau(T)}$ is not simple, and so the above minimization over $T$ to find $T^{\mbox{\tiny{opt}}}$ must be performed numerically.  Furthermore, it still remains to take expectations, and minimize over the staffing level to numerically solve for
\[
N^{\mbox{\tiny{opt}}} = \argmin_{N \in \{0,1,2,\ldots\}} cN + E \left[ z_{\tau(T^{\mbox{\tiny{opt}}})} (N,\Lambda) \right],
\]
and the associated minimum cost $\mathcal{C}^{\mbox{\tiny{opt}}}$.  To do this, we must perform an exhaustive search over $N$. The reason an exhaustive search should be performed is that it is very difficult to establish, in general, that the cost in (\ref{eq:objective_function}) is convex in the staffing level $N$.  In fact, even for a system where the arrival rate is known and there is no outsourcing option, convexity results are yet to be established when $\mu < \gamma$ (\cite{Armony09} and \cite{KoolePot11}).

 The exhaustive search to find $N^{\mbox{\tiny{opt}}}$ involves a numeric integration to calculate the expectation with respect to the arrival rate $\Lambda$, and,  for each value $l$ used in the numeric integration, there is another search that must be performed to find $z^{\mbox{\tiny{opt}}}(N,l)$. In other words, the exhaustive search algorithm is not a simple line search as it includes three nested layers of enumeration that correspond to the staffing level, the arrival rate used in numeric integration and the outsourcing threshold.
The exhaustive search algorithm is formally described as follows:

\noindent\textbf{Initialization:} Set $N^0=0$, $\mathcal{C}^0=\mathcal{C}((0, \pi^{\mbox{\tiny{opt}}}(0,\Lambda)))=\min(a,p)\lambda$, and $N=1$.  Decide on the maximum possible staffing level to allow, $N_{\max}$.

\noindent\textbf{Step 1:} Compute $\mathcal{C}((N, \pi^{\mbox{\tiny{opt}}}(N,\Lambda)))=cN + E\left[ z^{\mbox{\tiny{opt}}}\left(N,\Lambda\right) \right]$ via numeric integration. For each possible arrival rate realization $l$ in the numeric integration, initialize\footnote{Lemma 3.2 in \cite{KocagaWard10} implies $T^{\mbox{\tiny{opt}}}\geq N$.} $T=N$, decide on the stopping criterion\footnote{For example, Theorem 3.4 in \cite{KocagaWard10} provides a bound on the difference between the current cost and the minimum cost.}, and then compute $T^{\mbox{\tiny{opt}}} = T^{\mbox{\tiny{opt}}}(N,l)$ as follows.
\begin{enumerate}
\item[(A)] Solve for $z_{\tau(T)}$ from the steady-state probabilities $\{ \theta_n(T): n=0,1,2,\ldots,T \}$.
\item[(B)] If $z_{\tau(T+1)} \geq z_{\tau(T)}$\footnote{If $z_{\tau(T+1)} \geq z_{\tau(T)}$, Theorem 3.2 in  \cite{KocagaWard10} implies that $T^\star=T$.}, or if the stopping criterion holds, set  $T^{\mbox{\tiny{opt}}}=T$ and stop.  Otherwise, increase $T$ by 1 and go to step (A).
\end{enumerate}

\noindent\textbf{Step 2:} If $\mathcal{C}((N, \pi^{\mbox{\tiny{opt}}}(N,\Lambda)))<\mathcal{C}^{0}$, then $N^0\leftarrow N$ and $\mathcal{C}^0\leftarrow\mathcal{C}((N, \pi^{\mbox{\tiny{opt}}}(N,\Lambda)))$

\noindent\textbf{Step 3:}  If $N=N_{\max}$, then set $N^{\mbox{\tiny{opt}}}=N^0$ and $\mathcal{C}^{\mbox{\tiny{opt}}}=\mathcal{C}^0$ and stop. Otherwise, increase $N$ by 1 and go to Step 1.

Although the algorithm above can compute the optimal policy numerically, it does not provide any insight with regards to the structure of the optimal policy. Furthermore, the computation time to obtain an optimal policy can be several hours for large system sizes. (We implemented the algorithm in Matlab.) Therefore we take the following approach to develop our proposed policy:  we evaluate the performance of the family of square root safety staffing policies combined with threshold routing.   To do this, we first assume that the form of the arrival rate uncertainty is on the order of the square root of the mean arrival rate, and then show that square root safety staffing combined with threshold routing is asymptotically optimal (Section \ref{Section:aoHW}).  Second, we propose a universal policy $U$ that is based on that asymptotic optimality result (Section \ref{Section:ProposedPolicy}), and show numerically that  not only its computation time is in the order of seconds, it also has a very good performance  even outside of the regime in which we proved its asymptotic optimality (Section \ref{Section:Numerics}).

\section{Asymptotic Analysis}
\label{Section:aoHW}

In this section of the paper only, we assume that the order of uncertainty in the arrival rate is the same as the square-root of the mean of the arrival rate.  To do this, we consider a sequence of systems indexed by the mean arrival rate $\lambda$, and let $\lambda \rightarrow \infty$.  We assume that the random arrival rate $\Lambda$ throughout this sequence of systems can be expressed as
\begin{equation}
{\Lambda} = \lambda +  {X}\sqrt{\lambda},
\label{eq:QED}
\end{equation}
where $X$ is a random variable with mean zero and has $E|X|<\infty$. For this section only, the expectation operator is with respect to $X$ (instead of $\Lambda$).   We note that this form for the arrival rate is a special case of the model assumed in~\cite{Maman09}.   Our convention is to use the superscript $\lambda$ to denote a process or quantity associated with the system having random arrival rate $\Lambda=\Lambda^\lambda(X)$ given in (\ref{eq:QED}).  The notation
\[
l^\lambda(x) = \lambda + x \sqrt{\lambda}
\]
denotes the realized arrival rate in the system having mean arrival rate $\lambda$; the $\lambda$ superscript should remind the reader that $l^\lambda(x) \rightarrow \infty$ as $\lambda \rightarrow \infty$ for any $x \in (-\infty,\infty)$.

An admissible policy  $\mathbf{u} = (\mathbf{N},\boldsymbol{\pi}):= \{(N^\lambda,\pi^\lambda): \lambda \geq 0 \}$ refers to an entire sequence that specifies an admissible policy for each $\lambda$. In particular, $N^\lambda$ is a non-negative integer and $\pi^\lambda = \pi^\lambda(N^\lambda,l^\lambda(x)) \in \Pi$ for each $\lambda$ and any realization $x$ of $X$ (so that the system arrival rate is $l^\lambda(x)$).  The notation $z_{\boldsymbol{\pi}}^\lambda(\mathbf{N},l^\lambda(x))$ is the long-run average operating cost, as defined in (\ref{eq:stage2cost}), for the system with realized arrival rate $l^\lambda(x)$,  and  $z_{\boldsymbol{\pi}}^\lambda(\mathbf{N},\Lambda^\lambda(X))$ is the associated random variable.  Similarly, $P_{\boldsymbol{\pi}}^\lambda(\mbox{out};l^\lambda(x))$  ($P_{\boldsymbol{\pi}}^\lambda(\mbox{ab};l^\lambda(x))$) is the steady-state probability an arriving customer is routed to the outsourcer (probability of abandonment) when the realized arrival rate is $l^\lambda(x)$, and  $P_{\boldsymbol{\pi}}^\lambda(\mbox{out};\Lambda^\lambda(X))$  ($P_{\boldsymbol{\pi}}^\lambda(\mbox{ab};\Lambda^\lambda(X))$) is the associated random variable.

In this section, we first define what we mean by asymptotic optimality (Section~\ref{subsection:aoDef}).   Then, we perform an asymptotic analysis in order to understand the behavior of the family of square-root staffing policies combined with threshold routing (Section~\ref{subsection:AsymptoticSqS}).  Finally, we optimize over the aforementioned policy class to obtain our proposed policy (Section~\ref{subsection:ProposedPolicy}), and we establish its asymptotic optimality.

\subsection{The Asymptotic Optimality Definition} \label{subsection:aoDef}
Our definition of asymptotic optimality is motivated by first observing that the lowest achievable cost on fluid scale is $c\lambda + o(\lambda)$, where the notation $f^\lambda = o(g^\lambda)$ means that $\lim_{\lambda \rightarrow \infty} f^\lambda/ g^\lambda = 0$.

\begin{proposition}\label{prop:fluid_bound} (Fluid-scaled cost.)
Under the assumption (\ref{eq:QED}) we have that:
\begin{enumerate}[(i)]
\item Any admissible policy $\mathbf{u} = (\mathbf{N},\boldsymbol{\pi})$ has
\[
\liminf_{\lambda\rightarrow\infty} \frac{cN^\lambda+E\left[z_{\boldsymbol{\pi}}^\lambda(\mathbf{N},\Lambda^\lambda(X))\right]}{\lambda}\geq c.
\]
\item If $N^\lambda = \lambda + \beta\sqrt{\lambda} + o\textstyle(\lambda)$, then, under the routing policy $\boldsymbol{\tau(\infty)}$ that outsources no customers,
\[
\lim_{\lambda\rightarrow\infty} \frac{cN^\lambda+E\left[z_{\boldsymbol{\tau(\infty)}}^\lambda(\mathbf{N},\Lambda^\lambda(X))\right]}{\lambda}= c.
\]
\end{enumerate}
\end{proposition}

The following refined and diffusion scaled cost function (defined for the any admissible policy $\mathbf{u} = (\mathbf{N},\boldsymbol{\pi})$)
\begin{equation}
\label{Eq: Centered and Scaled Cost}
    \hat{\mathcal{C}}^{\lambda}(\mathbf{u}) := \sqrt{\lambda} \left( \frac{cN^\lambda+E\left[z_{\boldsymbol{\pi}}^\lambda(\mathbf{N},\Lambda^\lambda(X))\right]}{\lambda} - c \right) \geq 0
\end{equation}
captures both the cost of additional staffing (above the offered load level $\lambda$) and the cost of the routing control.

\begin{definition} \label{asymptotic optimality} (Asymptotic optimality.)
An admissible policy $\mathbf{u}^\star = (\mathbf{N^\star},\boldsymbol{\pi^\star}) = \{ (N^{\lambda,\star}, \pi^{\lambda,\star}(N^{\lambda,\star},\Lambda^\lambda(X))): \lambda \geq 0 \}$ is asymptotically optimal if

\[\limsup_{\lambda \rightarrow \infty} \hat{\mathcal{C}}^{\lambda}(\mathbf{u})< \infty
\]
 and
\[
    \limsup_{\lambda \rightarrow \infty} \hat{\mathcal{C}}^{\lambda}(\mathbf{u^\star})\leq \liminf_{\lambda \rightarrow \infty} \hat{\mathcal{C}}^{\lambda}(\mathbf{u}),
\]
for any admissible policy $\mathbf{u}$.
\end{definition}

\subsection{The Asymptotic Behavior of Square Root Safety Staffing Combined with Threshold Routing} \label{subsection:AsymptoticSqS}

It has been shown in the extensive literature on staffing in large-scale service systems (e.g. \cite{HalfinWhitt81,Borst04,MandelbaumZeltyn2009}) that when the arrival rate is deterministic, square root safety staffing performs extremely well in minimizing both the staffing plus delay costs as well as staffing costs subject to performance constraints. When the arrival rate $\lambda$ is large, under square root safety staffing, the waiting times are small (at the order of $1/\sqrt{\lambda}$), so that the percentage of customers that should be routed to the outsourcer (\cite{KocagaWard10}), as well as the percentage of customers that abandon (\cite{GarnettMandelbaumReiman02}) are both small. This suggests that square root safety staffing should be also relevant when the arrival rate is random.  Similarly to \cite{KocagaWard10}, to route calls, we use a threshold routing policy, $\boldsymbol{\tau} = \{\tau(T^\lambda): \lambda \geq 0\}$, as defined in Section~\ref{Section:Exact}.  The threshold level $T^\lambda = T^\lambda(N^\lambda, \Lambda^\lambda(X))$ is determined after the arrival rate realizes as $l^\lambda(x)$.

The following lemma establishes the asymptotic behavior of square root safety staffing combined with threshold routing for a fixed realization $x$ of $X$.  Let $\phi$ and $\Phi$ be the standard normal pdf and cdf, respectively.
\begin{lemma} \label{lemma:KW2010} (Asymptotic behavior with deterministic arrival rate.)
Suppose the random variable $X$ realizes as the value $x \in (-\infty,\infty)$.  Assume the policy $\mathbf{u} = (\mathbf{N},\boldsymbol{\tau})$ is such that
\begin{eqnarray}
N^\lambda & = & \lambda + \beta \sqrt{\lambda} + o(\sqrt{\lambda}) \label{eq:KW2010lemmaN} \\
T^\lambda & = & N^\lambda + \hat{T} \sqrt{l^\lambda(x)}, \mbox{ where } \hat{T} := \hat{T}(\beta,x) \in [0,\infty). \label{eq:KW2010lemmaT}
\end{eqnarray}
Suppose the initial number of customers in the system $Y_0^\lambda$ is such that
 $\frac{Y_0^\lambda-N^\lambda}{\sqrt{\lambda}} \Rightarrow \hat{Y}(0)$ as $\lambda \rightarrow \infty$, for some random variable $\hat{Y}(0)$ that is finite with probability 1.  Then,
 \[
 \frac{1}{\sqrt{\lambda}}z_{\boldsymbol{\tau}}^\lambda(\mathbf{N},l^\lambda(x)) \rightarrow \hat{z}(\beta-x,\hat{T}), \mbox{ as } \lambda \rightarrow \infty,
 \]
 where
 \begin{equation} \label{eq:zdef}
 \hat{z}(m,\hat{T}) := \frac{A(m,\hat{T})}{B(m,\hat{T})}
 \end{equation}
 for
\begin{eqnarray*}
A(m,\hat{T}) & := & p \phi\left( \sqrt{\gamma} \left( \hat{T}+\frac{m}{\gamma} \right) \right)  \\
& & + \left(a+\frac{w}{\gamma} \right) \left[ \phi\left(\frac{m}{\gamma} \right) -  \phi\left( \sqrt{\gamma} \left( \hat{T}+\frac{m}{\gamma} \right) \right)
 + \frac{m}{\sqrt{\gamma}} \left(\Phi\left( \frac{m}{\sqrt{\gamma}}\right) -  \Phi\left( \sqrt{\gamma} \left( \hat{T}+\frac{m}{\gamma} \right) \right) \right)\right] \\
B(m,\hat{T}) & := & \frac{\phi\left( \frac{m}{\sqrt{\gamma}} \right)}{\phi(m)}\Phi(m) + \frac{1}{\sqrt{\gamma}} \left( \Phi\left( \sqrt{\gamma}\left( \hat{T}+\frac{m}{\gamma} \right) \right) - \Phi\left( \frac{m}{\sqrt{\gamma}} \right) \right)
\end{eqnarray*}
\end{lemma}
\noindent The appearance of $\beta-x$ as an argument in $\hat{z}$ in Lemma~\ref{lemma:KW2010} occurs because under (\ref{eq:QED}) the staffing $N^\lambda$ in (\ref{eq:KW2010lemmaN}) is such that the system operates in the QED regime regardless of the realization $x$ of $X$; in particular,
\[
\frac{N^\lambda - l^\lambda(x)}{\sqrt{\lambda}} = \frac{N^\lambda - \lambda}{\sqrt{\lambda}} - x \rightarrow \beta - x \mbox{ as } \lambda \rightarrow \infty.
\]
Furthermore, the following Corollary to Lemma~\ref{lemma:KW2010} highlights that the dependence of the threshold level on the realized arrival rate is through the definition of $\hat{T}$, and not through its multiplier (which is always of order $\sqrt{\lambda}$ under the assumption (\ref{eq:QED})).

\begin{corollary} \label{corollary:KW2010}
Lemma~\ref{lemma:KW2010} continues to hold when $T^\lambda$ in (\ref{eq:KW2010lemmaT}) is re-defined as
\[
T^\lambda  =  N^\lambda + \hat{T} \sqrt{\lambda}.
\]
\end{corollary}

% \begin{equation*}
%\begin{split}
%A(m,\hat{T}) & :=
% p \exp\left( -\frac{\gamma}{2}\left(\hat{T}^2+2\frac{m}{\gamma}\hat{T}\right)\right)\\
%&+ (a+\frac{w}{\gamma}) \left( 1-\exp\left( -\frac{\gamma}{2}\left(\hat{T}^2+2\frac{m}{\gamma}\hat{T}\right)\right)
%+ \frac{\sqrt{2\pi}}{\sqrt{\gamma}}m\left(\Phi\left(\frac{m}{\sqrt{\gamma}}\right)-\Phi\left(  \sqrt{\gamma}(\hat{T}+\frac{m}{\gamma})\right) \right)\right) \\
%B(m,\hat{T}) & := \sqrt{2\pi}\left(\exp\left(\frac{m^2}{2}\right)\Phi\left(m\right) + \frac{1}{\sqrt{\gamma}} \exp\left(\frac{m^2}{2\gamma}\right)\left( \Phi\left( \sqrt{\gamma}
%\left( \hat{T} + \frac{m}{\gamma} \right) \right) - \Phi\left( \frac{m}{\sqrt{\gamma}} \right) \right)\right).
%\end{split}
%\end{equation*}

The issue is that in order to analyze the performance of square root safety staffing combined with threshold routing, we require that Lemma~\ref{lemma:KW2010} and Corollary \ref{corollary:KW2010} hold when the fixed value $x$ is replaced by the random variable $X$.
\begin{theorem}\label{Thm: Cost Convergence} (Asymptotic cost convergence.)
 Assume the policy $\mathbf{u} = (\mathbf{N},\boldsymbol{\tau})$ is as defined by the equations (\ref{eq:KW2010lemmaN}) and (\ref{eq:KW2010lemmaT}).  Then, under the conditions of Lemma~\ref{lemma:KW2010},
\[
    \hat{\mathcal{C}}^{\lambda}(\mathbf{u}) \rightarrow \hat{\mathcal{C}}(\mathbf{u}):=c \beta + E\left[ \hat{z}(\beta - X,\hat{T}) \right], \mbox{ as } \lambda\rightarrow\infty.
\]
\end{theorem}

\subsection{The Proposed Policy} \label{subsection:ProposedPolicy}

 It is sensible to set the parameters $\beta$ and $\hat{T}$ of Lemma \ref{lemma:KW2010} in order to minimize the limiting cost $\hat{\mathcal{C}}(\mathbf{u})$ of Theorem \ref{Thm: Cost Convergence}.  The first step is to observe that, for $p<a$ and any given $\beta$, Proposition 4.1 in~\cite{KocagaWard10} shows that for the realized arrival rate $l^{\lambda}(x)$, the unique $\hat{T}^\star = \hat{T}^\star(\beta-x) < \infty$ that solves
 \begin{equation}\label{Eq: T def}
    (a-p)\gamma \hat{T} - \hat{z}(\beta-x,\hat{T}) = p(\beta-x)
\end{equation}
has the property that
\begin{equation}\label{Eq: Dif Cost LB}
\hat{z}(\beta-x,\hat{T}^\star) \leq \hat{z}(\beta-x, \hat{T})
\end{equation}
for any other $\hat{T} \geq 0$. Otherwise, for $a\leq p$, $\hat{T}^\star=\infty$ and
\begin{equation}\label{Eq: Dif Cost Inf LB}
\hat{z}(\beta-x,\infty):=\lim_{\hat{T}\rightarrow\infty}\hat{z}(\beta-x,\hat{T}) \leq \hat{z}(\beta-x, \hat{T}_0)
\end{equation}
for any finite $\hat{T}_0 \geq 0$. The second step is to plug $\hat{T}^\star$ into the limiting expression in Theorem~\ref{Thm: Cost Convergence}, and to optimize over $\beta$ to find
\begin{equation}\label{Eq: QED beta}
\beta^{\star}:=\argmin\limits_{\beta }\left\{ c\beta +E\left[ \hat{z }(\beta -{X},\hat{T}^\star(\beta -{X}))\right] \right\}.
\end{equation}

It is important to observe that $\beta^{\star}$ is well-defined in the sense that $\beta^{\star}$ is finite and
\[
 c\beta^{\star} +E\left[ \hat{z }(\beta^{\star} -{X},\hat{T}^\star(\beta^{\star} -{X}))\right]=\inf\limits_{\beta\in(-\infty,\infty)} c\beta +E\left[ \hat{z }(\beta -{X},\hat{T}^\star(\beta -{X}))\right]<\infty.
\]
This follows from the next two propositions.
\begin{proposition}\label{prop: beta_cost_finite}
For any $\beta\in(-\infty,\infty)$, $c\beta +E\left[ \hat{z }(\beta -{X},\hat{T}^\star(\beta -{X}))\right]<\infty$.
\end{proposition}
\begin{proposition}\label{prop: beta_star_finite}
The infimum in $\inf\limits_{\beta\in(-\infty,\infty)} c\beta +E\left[ \hat{z }(\beta -X,\hat{T}^\star(\beta -X))\right]$ is attained by a finite $\beta~\in~(-\infty,\infty)$.
\end{proposition}

%\subsection{The Proposed Policy} \label{subsection:ProposedPolicy}

We are now in a position to define our proposed policy: We let
\begin{equation} \label{eq:proposed_square_root_staffing}
\mathbf{u^\star}=(\mathbf{N^\star},\boldsymbol{\tau^\star}) := \{ (N^{\lambda,\star},\tau(T^{\lambda,\star})) : \lambda \geq 0 \}
\end{equation}
that has the staffing level
\begin{equation} \label{eq:NstarDef}
N^{\lambda,\star} = \lambda + \beta^\star{\sqrt{\lambda}}
\end{equation}
and sets the threshold level $T^{\lambda,\star} = T^{\lambda,\star}(N^{\lambda,\star},\Lambda^\lambda(X))$ when the arrival rate realizes as $l^\lambda(x)$ as
\begin{equation} \label{eq:l_value}
T^{\lambda,\star} = N^{\lambda,\star} + \hat{T}^\star(\beta^\star-x)\times \sqrt{l^{\lambda}(x)},
\end{equation}
for $\hat{T}^\star(\beta^\star-x)$ defined by (\ref{Eq: T def}) with $\beta^\star$ replacing $\beta$.

Theorem~\ref{Thm: Cost Convergence} is valid for $\mathbf{u^\star}$,
and so
\[
\hat{\mathcal{C}}^{\lambda}(\mathbf{u^\star}) \rightarrow \hat{\mathcal{C}}^\star :=  c \beta^{\star} + E \left[ \hat{z}(\beta^{\star} - X, \hat{T}^\star(\beta^\star - X)) \right]
\]
Our next result confirms that $\mathcal{C}^\star$ is the minimum achievable cost, meaning that the policy $\mathbf{u}^\star$ is asymptotically optimal.
\begin{theorem}\label{Thm: AO} (Asymptotic optimality of our proposed policy.)
 The policy $\mathbf{u^\star}$, defined through  (\ref{Eq: QED beta}), (\ref{eq:proposed_square_root_staffing}), (\ref{eq:NstarDef}), and (\ref{eq:l_value}) is asymptotically optimal under (\ref{eq:QED}); i.e., under any other admissible policy $\mathbf{u}$
\[\liminf_{\lambda \rightarrow \infty}  \hat{\mathcal{C}}^{\lambda}(\mathbf{u})\geq \hat{\mathcal{C}}^{\star}.
\]
 Furthermore, it follows that our proposed policy has associated cost that is $o(\sqrt{\lambda})$ higher than the minimum achievable cost for a given $\lambda$; i.e., that
\[
\frac{cN^{\lambda,\star} + E \left[ z_{\boldsymbol{\tau}^\star}^\lambda\left( \mathbf{N^\star}, \Lambda^\lambda(X) \right) \right]-\mathcal{C}^{\lambda,\mbox{\tiny{opt}}}}{\sqrt{\lambda}} \rightarrow 0, \mbox{ as } \lambda \rightarrow \infty,
\]
where $\mathcal{C}^{\lambda,\mbox{\tiny{opt}}} := \mathcal{C}^{\mbox{\tiny{opt}}}$ for $\mathcal{C}^{\mbox{\tiny{opt}}}$ defined in (\ref{eq:objective}) for the system with mean arrival rate $\lambda$.
\end{theorem}

\begin{remark} (Performance under the optimal threshold.)
Another asymptotically optimal policy
\[
(\mathbf{N^\star},\boldsymbol{\pi^{\mbox{\tiny{opt}}}}) := \{ (N^{\lambda,\star},\pi^{\lambda,\mbox{\tiny{opt}}}) : \lambda \geq 0 \}
\]
has staffing levels $N^{\lambda,\star}$ defined in (\ref{eq:NstarDef}), and, after the arrival rate $\Lambda^\lambda(X)$ realizes as $l^\lambda(x)$, solves the relevant Markov decision problem for the routing control policy $\pi^{\lambda,\mbox{\tiny{opt}}} = \pi^{\lambda,\mbox{\tiny{opt}}}(N^{\lambda,\star},l^\lambda(x))$ that achieves the minimum long-run average operating cost $z^{\lambda,{\mbox{\tiny{opt}}}}(\mathbf{N^\star},l^\lambda(x))$.  To see this, it is enough to observe that
\[
z_{\boldsymbol{\tau^\star}}^\lambda(\mathbf{N^\star},l^\lambda(x)) \geq z^{\lambda,\mbox{\tiny{opt}}}(\mathbf{N^\star},l^\lambda(x)),
\]
for every $\lambda$ and any realization $x$ of $X$.
\end{remark}

In the following,  $f^\lambda=\mathcal{O}(g^\lambda)$ means that $\limsup_{\lambda\rightarrow\infty}|f^\lambda / g^\lambda| < \infty$.
\begin{remark} \label{remark:BRZ} (Comparison to~\cite{BassambooRandhawaZeevi10}.)
When $a<p$ (in addition to our assumption that $c < \min(a,p)$), it follows from Proposition~\ref{prop: prelimit optimal} part (ii) that the optimal control policy does not outsource any calls.  Then, the cost minimization problem (\ref{eq:objective_function}) is a pure staffing problem (instead of a joint staffing and routing problem), which is equivalent to the problem solved in~\cite{BassambooRandhawaZeevi10}.  Theorem 1 part (c) of that paper, adapted to our setting, shows that a policy based on a newsvendor prescription can have associated cost that is $\mathcal{O}(\sqrt{\lambda})$ higher than the minimum achievable cost for a given $\lambda$.  In comparison, our proposed policy has associated cost that is $o(\sqrt{\lambda})$ higher than the minimum achievable cost for a given $\lambda$ by Theorem~\ref{Thm: AO}.
 Hence we expect our policy to provide significant improvements over that of \cite{BassambooRandhawaZeevi10}, as the arrival rate uncertainty decreases.  
\end{remark}

\section{The Proposed Universal Policy}
\label{Section:ProposedPolicy}

For models that do not assume uncertain arrival rates, square root safety staffing is known in the literature to be very robust.  For an $M$/$M$/$N$ queue with no abandonments, no dynamic routing decisions, and known arrival rate, ~\cite{Borst04} show that square root safety staffing performs extremely well, {\em both} inside and outside of the parameter regime (linear staffing and waiting costs) in which they prove it to be asymptotically optimal (see their numerical experiments in Section 10).  In a more recent paper, \cite{GHM2012}  prove that performance approximations that are based on the premise that the staffing is of a square-root safety form are asymptotically universally accurate, as the arrival rate becomes large.  This latter paper is also limited to the case of deterministic arrival rates and no dynamic control.

This leads us to propose the universal  policy $U$ when there are {\em no restrictions} on the form of the arrival rate uncertainty, as in (\ref{eq:QED}).  We define
$U$ for the model as specified in Section~\ref{Section: Model},  and analyzed exactly in Section~\ref{Section:Exact}, without considering a sequence of systems as in Section~\ref{Section:aoHW}.  To do this, we begin with the non-negative random variable $\Lambda$ that represents the system arrival rate and has mean $E[\Lambda]=\lambda$.  Then, we make the transformation
\begin{equation}\label{Eq: UniversalX}
X := \frac{\Lambda - \lambda}{\sqrt{\lambda}},
\end{equation}
and use the random variable $X$ to define $U$
\[
U = (N_U,\pi_{U}(N_U,\Lambda)).
\]
The proposed staffing level is
\[
N_{U} = \left[ \lambda + \beta^\star{\sqrt{\lambda}} \right],
\]
for $\beta^\star$ that satisfies (\ref{Eq: QED beta}), with $X$ in that expression defined by (\ref{Eq: UniversalX}), and the function $[\cdot]$ rounds the expression inside the brackets to the nearest integer.  The proposed routing policy when the arrival rate $\Lambda$ realizes as $l$ is the threshold routing policy
\[
\pi_{U}(N,l) = \tau(T_{U})
\]
for
\[
T_{U} = N_{U} + \hat{T}^\star \sqrt{l},
\]
and $\hat{T}^\star$ defined by (\ref{Eq: T def}), with $x$ in that expression replaced by $(l-\lambda)/\sqrt{\lambda}$.

The universal $U$ policy ``pretends'' that the magnitude of the uncertainty in the arrival rate $\Lambda$ is on the order of $\sqrt{\lambda}$, as in (\ref{eq:QED}), and sets $\beta^\star$ and $\hat{T}^\star$ accordingly. In contrast to the policy defined in Section~\ref{subsection:ProposedPolicy} under assumption (\ref{eq:QED}), the magnitude of the second order term appearing in the definitions of $N_{U}$ and $T_{U}$ may not be of order $\sqrt{\lambda}$.  In particular, depending on the distribution of $\Lambda$, the value of $\beta^\star$ may end up being of the same order of $\sqrt{\lambda}$, so that the second term in $N_U$ is of order $\lambda$ (see discussion in Section \ref{subsec:discussion}).  This flexibility suggests that $U$ may perform well, even outside of the regime in which it is proved to be asymptotically optimal, as we indeed observe in the next section.

\section{Numerical Evaluation of the Proposed Policy}
\label{Section:Numerics}
Theorems \ref{Thm: Cost Convergence} and \ref{Thm: AO} establish when the order of uncertainty in the arrival rate is the same as the square-root of the mean arrival rate, so that (\ref{eq:QED}) holds, $U$ staffs and routes in a way that achieves minimum cost for large enough $\lambda$. However, Theorems \ref{Thm: Cost Convergence} and \ref{Thm: AO} do not provide guidance on: how large $\lambda$ must be, what happens when (\ref{eq:QED}) does not hold, or how $U$ performs in comparison to alternative benchmark policies. In this section we show that  $U$ generally achieves within 0.1\% of the minimum cost even when these assumptions are relaxed, and its robustness in comparison to two benchmark policies is the highest. To do this, we first vary the system size expressed by the mean arrival rate $\lambda$ (Section \ref{Subsection: finite_system_size}) and the level of uncertainty in $\Lambda$ (Section \ref{Subsection: varying_CV}) to gain an initial conclusion that $U$  performs {\em remarkably} well. Then we show that this conclusion is, to a large degree, insensitive to changes in the cost parameters (Section \ref{Subsection: varying_cost}) and the asymmetry of the arrival rate distribution (Section \ref{Subsection: varying_skewness}). In summary, $U$ performs extremely well, even when the system is far away from the regime in which it is proved to be asymptotically optimal.

Throughout our numerical examples, we set the mean service time and the mean patience time equal to 1 and fix the cost parameters at $c=0.1$, $p=1$ and $a=5$ unless specified otherwise. It follows from Proposition \ref{prop: prelimit optimal} that our choice of cost parameters is such that it is optimal for the system manager to set a nonzero staffing level and routes some calls to the outsourcer.

\subsection{Finite System Size}\label{Subsection: finite_system_size}
Having established that  $U$  is asymptotically optimal as the system size grows without bound, under the form of uncertainty in $\Lambda$ as in (\ref{eq:QED}), we proceed to evaluate its performance for finite size systems. This evaluation is done by comparing $U$ to the numerically computed optimal staffing policy. We compute the optimal staffing level $N^{\mbox{\tiny{opt}}}$ via an exhaustive search, as described in Section~\ref{Section:Exact}.

Table \ref{tab:system size} illustrates the performance of our proposed staffing policy with respect to the optimal staffing level by varying the system size and letting the distribution of the arrival rate $\Lambda$ be in accordance with (\ref{eq:QED}). Specifically, we assume that $X$ follows a Uniform distribution on $[-1,1]$, and increase the mean arrival rate $\lambda$ from 1 to 1600. Then $\Lambda$ follows a Uniform distribution with its support interval increasing from $[0,2]$ to $[1560,1640]$. This is consistent with our assumption in Theorems 1 and 2 that prove asymptotic optimality of $U$ as $\lambda$ becomes large under (\ref{eq:QED}). The first and second columns in Table \ref{tab:system size} show the resulting distribution for $\Lambda$.
%We assume that $\Lambda$ follows a uniform distribution as reported in the first columns of Table \ref{tab:system size} and is chosen such that $\Lambda \buildrel d \over = \lambda+\sqrt{\lambda}X$ where $X\sim U[-1,1]$. As a result, the deviations in the random arrival rate are of order $\sqrt{\lambda}$ so that (\ref{eq:QED}) is satisfied and the system stays in the QED regime.
The third and fourth columns in Table \ref{tab:system size} show the optimal staffing level and the associated optimal average cost, while columns five and six show our proposed approximate staffing policy, along with its average cost. Column seven displays the staffing error which is the difference between the optimal staffing level and our approximate staffing level. Finally, column eight displays the percentage cost error with respect to the optimal policy.

We see from Table \ref{tab:system size} that $U$ performs extremely well for {\em all} system sizes, that are consistent with the assumption that the uncertainty in the arrival rate is of the same order as the square-root of the mean arrival rate. Notice that the percentage cost error may be nonzero even when the staffing error equals zero. This is because $U$ sets the threshold level according to (\ref{eq:l_value}) which may not equal to the optimal threshold. In light of Theorems \ref{Thm: Cost Convergence} and \ref{Thm: AO}, that establish asymptotic optimality, it is not surprising that our policy performs extremely well for large $\lambda$. The less expected numerical insight is that $U$ also performs extremely well for small $\lambda$. (We note that there is a chance that the rounding can go the wrong way, and subsequently may cause a large cost error in extremely small system size. However such small systems sizes are not realistic for most call center applications). In summary, $U$ is very robust to system size, provided the order of uncertainty in the arrival rate is as assumed in (\ref{eq:QED}).

\begin{table}[htbp]
  \centering
    \begin{tabular}{c|c||cc|cc|cc}
\hline

    $\lambda$& Distribution of  & \multicolumn{2}{c|}{Optimal Policy}  & \multicolumn{2}{c|}{$U$}  & \multicolumn{2}{c}{Difference}\\\cline{3-8}
    & $\Lambda$  & $N^{\mbox{\tiny{opt}}}$ & $\mathcal{C}^{\mbox{\tiny{opt}}}$ &  $N_{U}$  & $\mathcal{C}(N_{U})$ & $ N^{\mbox{\tiny{opt}}}-N_{U}$ &$\frac{\mathcal{C}(N_{U})-\mathcal{C}^{\mbox{\tiny{opt}}}}{\mathcal{C}^{\mbox{\tiny{opt}}}}$ \\ \hline\hline
    1     & U[0,2] & 3     & 0.4149 & 3     & 0.4188 & 0     & 0.9400\% \\
    9     & U[6,12] & 16    & 1.7702 & 15    & 1.7786 & 1     & 0.4745\% \\
    25    & U[20,30] & 36    & 3.8979 & 36    & 3.8998 & 0     & 0.0487\% \\
    100   & U[90,110] & 121   & 12.7131 & 121   & 12.7149 & 0     & 0.0142\% \\
    226   & U[210,240] & 257   & 26.5227 & 257   & 26.5236 & 0     & 0.0034\% \\
    400   & U[380,420] & 443   & 45.3338 & 442   & 45.3355 & 1     & 0.0037\% \\
    625   & U[600,650] & 678   & 69.1435 & 678   & 69.1441 & 0     & 0.0009\% \\
    900   & U[870,930] & 964   & 97.9536 & 963   & 97.9553 & 1     & 0.0017\% \\
    1600  & U[1560,1640] & 1685  & 170.5732 & 1684  & 170.5750 & 1     & 0.0011\%
 \\\hline\hline

 \end{tabular}
 \caption{Performance of $U$: increasing system size.}
  \label{tab:system size}
\end{table}

\subsection{Varying Arrival Rate Uncertainty}\label{Subsection: varying_CV}
Next, we evaluate the robustness of $U$  with respect to changes in the level of uncertainty in the arrival rate. This is important because the proof of asymptotic optimality of $U$  requires the assumption that the level of uncertainty in the arrival rate is of the same order as the square-root of the mean arrival rate (i.e., that (\ref{eq:QED}) holds). We measure the level of uncertainty in the arrival rate through its coefficient of variation $CV:=CV_{\Lambda}=\frac{\sqrt{\mbox{Var}[\Lambda]}}{\mbox{E}[\Lambda]}$. We are interested in both cases where the level of uncertainty in the arrival rate is lower than that assumed in (\ref{eq:QED}) and where it is higher.

In this subsection, we keep the mean arrival rate fixed at $\lambda=100$, and we assume that $\Lambda$ follows a Uniform distribution with support $[a,b]$. Then
\[
CV= \frac{1}{\sqrt{3}}\frac{b-a}{a+b}\leq \frac{1}{\sqrt{3}}=0.5774.
\]
In comparison, under assumption (\ref{eq:QED}), when $X$ follows a uniform distribution with support $[-1,1]$ as in Section \ref{Subsection: finite_system_size}, the coefficient of variation of the arrival rate $\Lambda=\Lambda(X)$ is
\begin{equation}\label{eq:CV Lambda}
CV_{\Lambda(X)}= \frac{\sqrt{Var\left[X\right]}}{\sqrt{\lambda}}= \frac{1}{10\sqrt{3}}=0.0577.
\end{equation}
Then by varying $CV$ from 0 to approximately 1/2, we cover both the cases where the level of uncertainty in the arrival rate is lower that (\ref{eq:QED}) and higher.

It is sensible to compare the performance of $U$ to two other possible staffing policies: one that is expected to perform well when the level of uncertainty in the arrival rate is low and the other that is expected to perform well when the level of uncertainty in the arrival rate is high. The first alternative policy we consider is $D$, a square root safety staffing policy that has the same form as (2), but chooses the coefficient of $\sqrt{\lambda}$ differently by assuming that the arrival rate is {\em deterministic} and fixed at the mean arrival rate $\lambda=100$. Specifically when the mean arrival rate is $\lambda$, $D$ staffs
\[
N_{D}:=\left[\lambda+\beta_1^\star \sqrt{\lambda}\right],
\]
where
\[
\beta_1^\star:=\argmin_{\beta}{c\beta+\hat{z}(\beta,\hat{T}^\star(\beta))}
\]
for $\hat{z}$ as defined in (\ref{eq:zdef}) and $\hat{T}^\star(\beta)$ that satisfies (\ref{Eq: T def}).  Note that  $D$ is exactly the proposed policy $U$ in the case $P(X=0)=1$. It is intuitive to expect that the performance of $D$ deteriorates significantly as $CV$ increases.

The second alternative policy we consider is $NV$, a newsvendor based prescription that is a modification of the policy proposed in \cite{BassambooRandhawaZeevi10} to include co-sourcing. The $NV$ policy follows a fluid approximation which ignores stochastic queueing effects and, as a result, when $a>p$, the abandonment cost is irrelevant. That is, in the fluid scale, all customers who cannot be served in-house immediately upon arrival will be outsourced. Similarly, when $a\leq p$ no calls will be outsourced. In particular, in newsvendor terminology the overage cost is $c$ (because of extra staffing) and the underage cost is $\min\{a,p\}-c$ (because we incur the cost of routing or cost of abandonment but do not incur the cost of an additional person for staffing). Then the critical ratio is $\frac{\min\{a,p\}-c}{\min\{a,p\}}$, and the newsvendor based staffing prescription is
\[
N_{NV}:=\left[ F_{\Lambda}^{-1}\left(\frac{\min\{a,p\}-c}{\min\{a,p\}}\right) \right].
\]
 We observe that when (\ref{eq:QED}) holds, $N_{NV}$ can also be written as
$N_{NV}:=\left[\lambda+\beta_2^\star\sqrt{\lambda}\right]$, where $\beta_2^\star:=F_{X}^{-1}(\frac{p-c}{p})$ and $F_{X}$ is the cumulative distribution function of $X$. Notice that, in sharp contrast to $D$, which disregards the uncertainty in the arrival rate, $NV$ disregards the inherent stochasticity of the system that produces queueing. Hence, we expect the performance of the $NV$ policy to deteriorate when $CV$ decreases.

We have specified the staffing rules, $N_{D}$ and $N_{NV}$, of two alternative policies $D$ and $NV$. There is still the question of what should be the routing policy. For this, we recall that after the arrival rate realizes, the optimal routing policy can be found by solving the relevant MDP (see \cite{KocagaWard10}). Hence, in our numerical experiments, after the arrival rate realizes, we operate both comparison policies under the {\em optimal} routing policy. The $U$ policy follows the threshold routing policy $\tau(T^{\lambda,\star})$ where $T^{\lambda,\star} = T^{\lambda,\star}(N^{\lambda,\star},\Lambda(X))$ is as defined in Section \ref{Section:ProposedPolicy} (although we observe that the performance of the diffusion based threshold routing policy and the exact solution to the relevant MDP are almost indistinguishable).

\begin{figure}[htb*]
\centering
\subfloat[\% cost error]{
\includegraphics[scale=0.5]{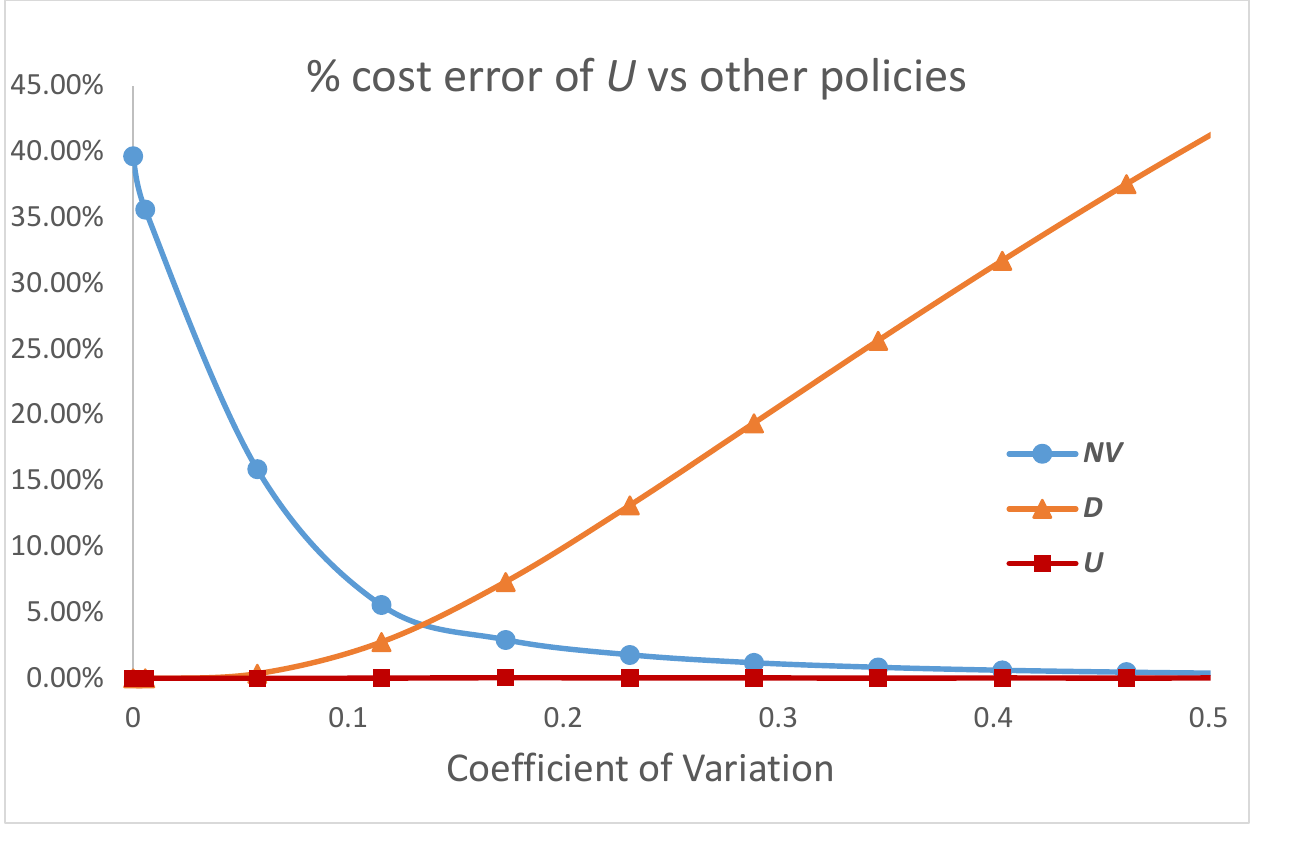}
}
\subfloat[staffing error (i.e., $N^{\mbox{\tiny{opt}}} - \triangle$ where $\triangle \in \{ N_U, N_{D}, N_{NV} \}$)]{
\includegraphics[scale=0.5]{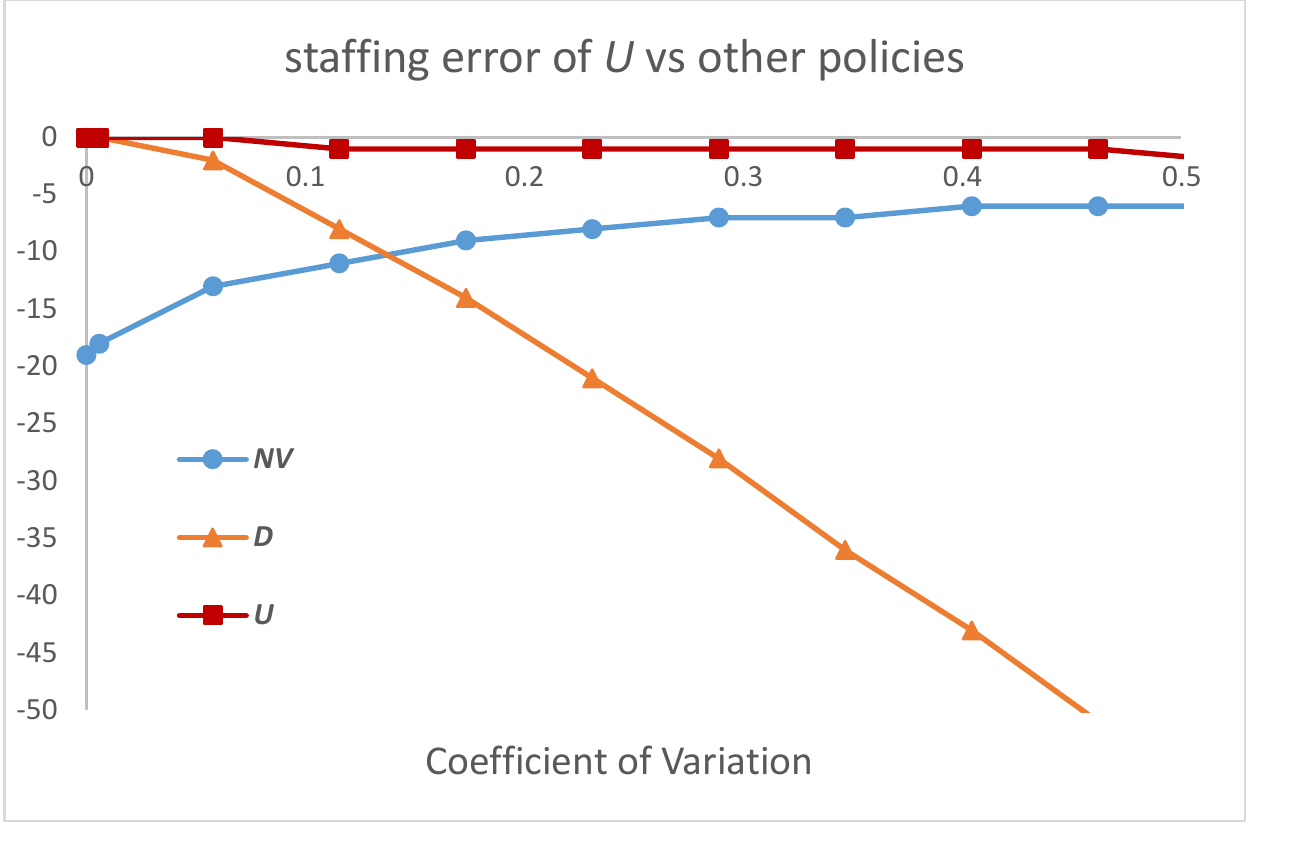}
}
\caption{Performance of $U$ and other policies for increasing CV}
\label{fig: T2}
\end{figure}

Figure \ref{fig: T2} plots the relative percentage cost error and staffing error of $U$, $NV$ and $D$. The staffing error and percentage cost error for $U$ is defined as in columns seven and eight in Table \ref{tab:system size}, and is defined similarly for $D$ and $NV$. Table \ref{tab: CV1} in EC contains further details regarding this study, such as the exact costs and staffing levels. We see from Figure \ref{fig: T2} that $U$ staffs very close to the optimal staffing level and therefore performs well even for very high CV values. We also see that $U$ outperforms $NV$ for lower CV values and outperforms $D$ for higher CV values. Furthermore, in both cases, the staffing and percentage cost error can be arbitrarily large. Hence we conclude that $U$ is robust and performs extremely well even in parameter settings beyond which it has been proven to be asymptotically optimal.

Figure \ref{fig: T2} is a first step in concluding that $U$ is very robust, and performs extremely well over a large range of parameter settings much beyond where Theorems \ref{Thm: Cost Convergence} and \ref{Thm: AO} establish its asymptotic optimality. The next step in establishing the aforementioned conclusion is to explore the effect of varying other parameters, for example the staffing cost.

\subsection{Varying Staffing Costs}\label{Subsection: varying_cost}

\begin{figure}[htbp]
\subfloat[low CV ($CV_\Lambda=0.0144$)]{
\includegraphics[scale=0.475]{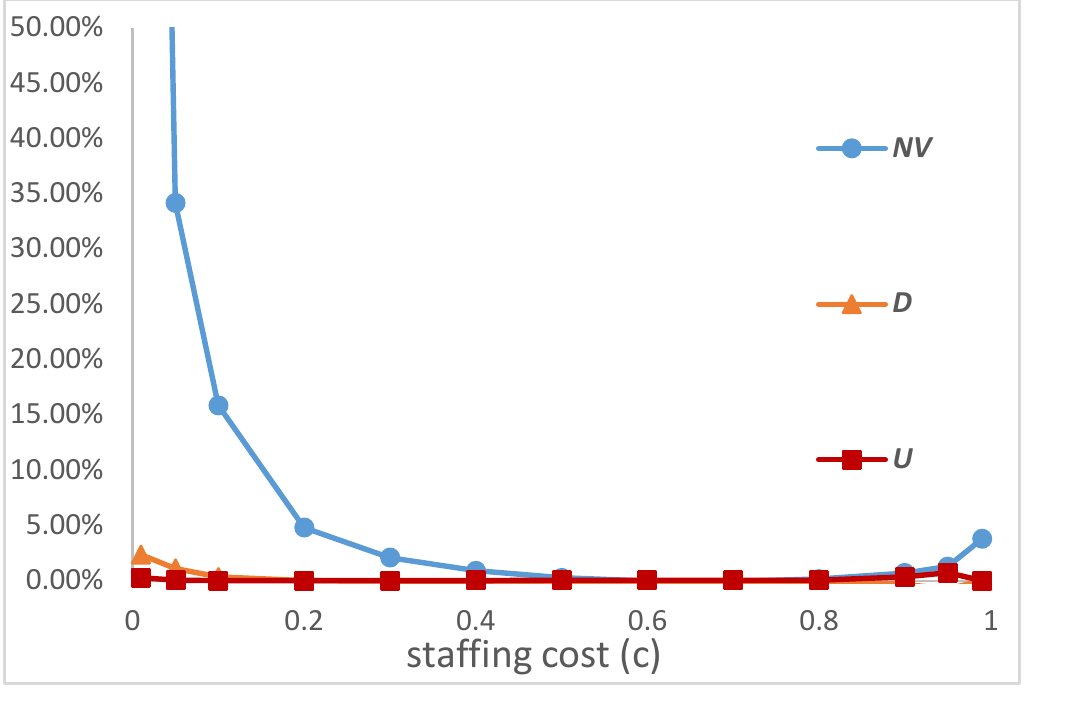}
\label{fig:crlowcv_cost}
}
\subfloat[moderate CV ($CV_\Lambda=0.0722$)]{
\includegraphics[scale=0.475]{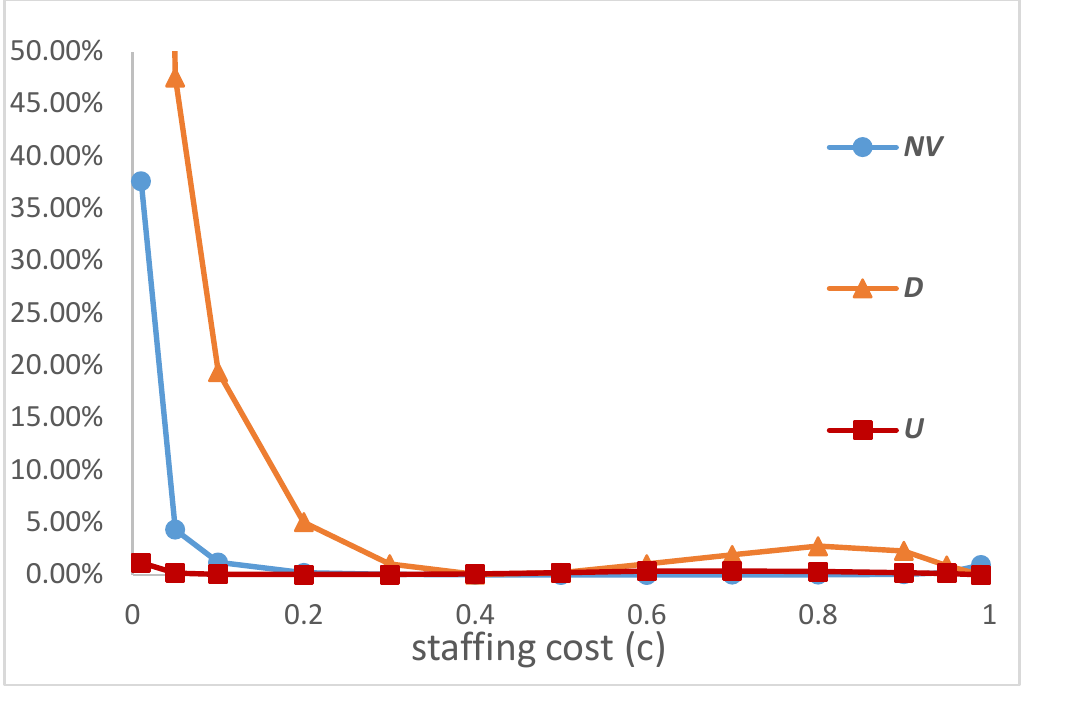}
\label{fig:crmedcv_cost}
}
\subfloat[high CV ($CV_\Lambda=0.1299$)]{
\includegraphics[scale=0.475]{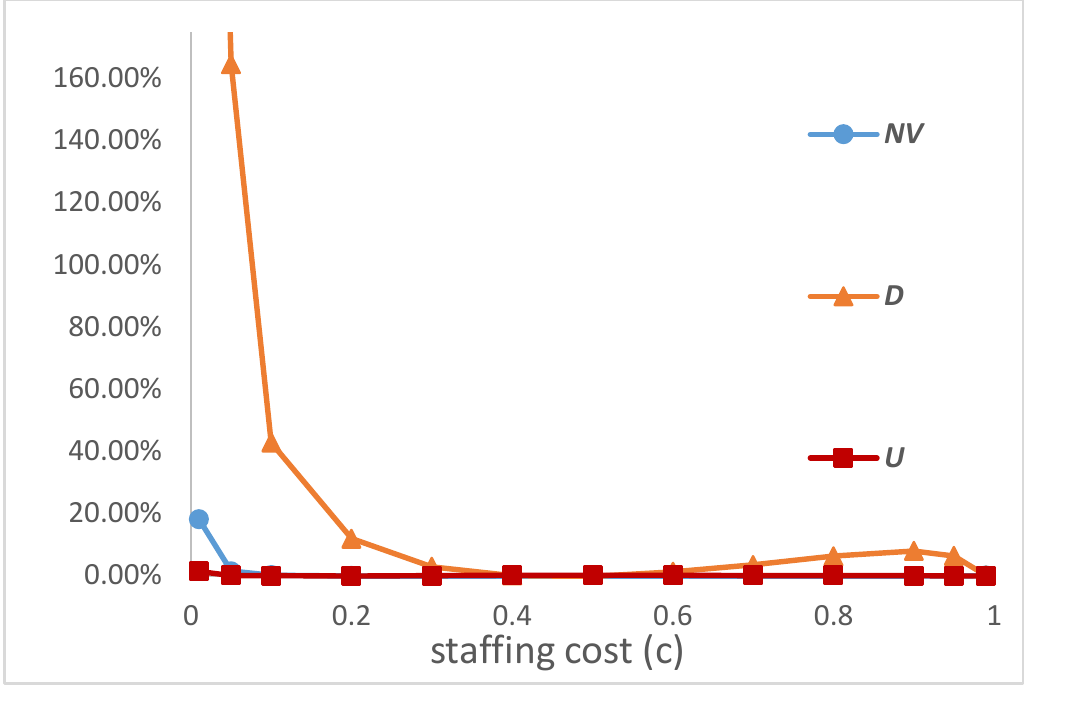}
\label{fig:crhighcv_cost}
}
\caption{\% cost error for changing staffing costs at three levels of arrival rate uncertainty}
\label{fig: cr_cost}
\end{figure}

\begin{figure}[htb*]
\subfloat[low CV ($CV_\Lambda=0.0144$)]{
\includegraphics[scale=0.475]{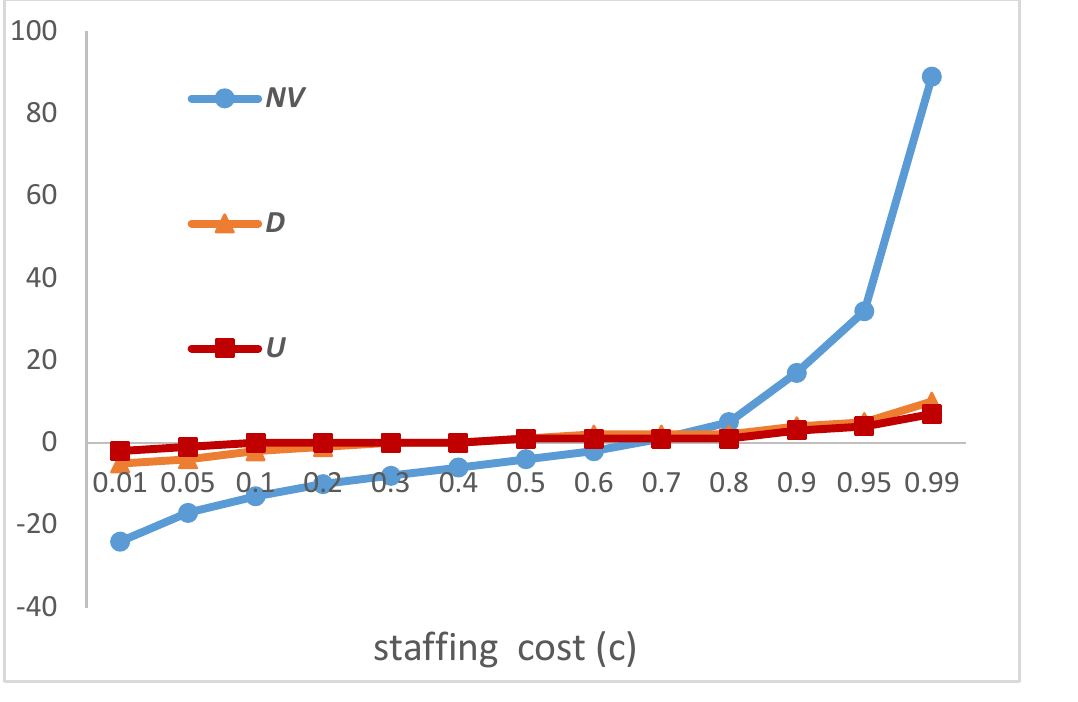}
\label{fig:crlowcv_staff}
}
\subfloat[moderate CV ($CV_\Lambda=0.0722$)]{
\includegraphics[scale=0.475]{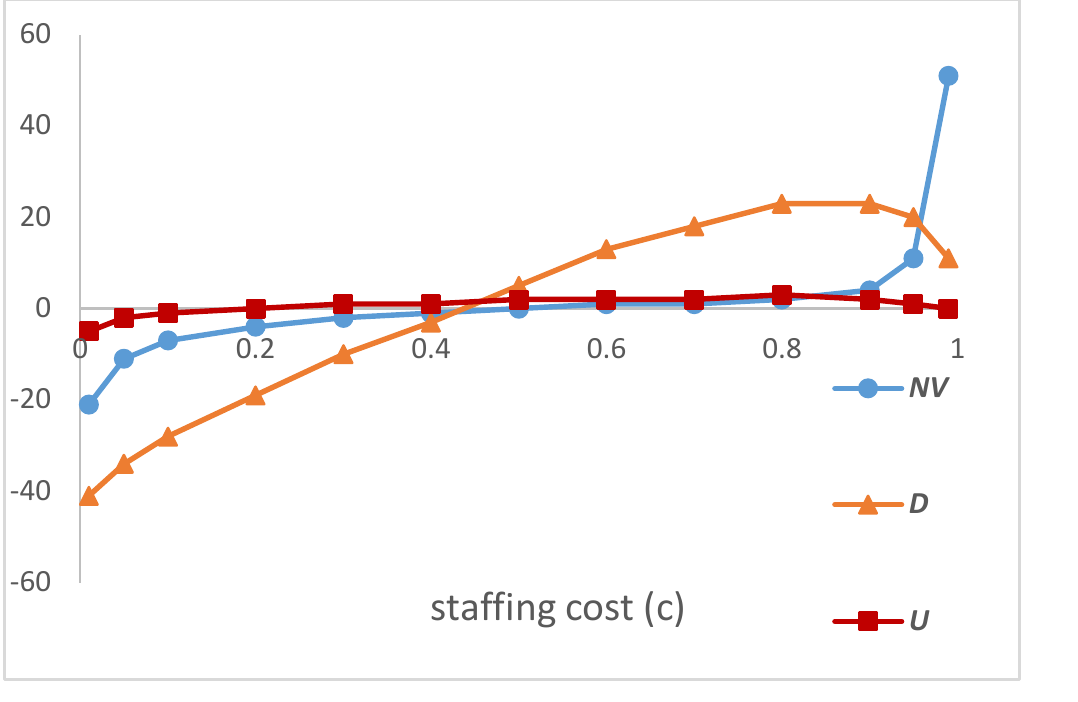}
\label{fig:crmedcv_staff}
}
\subfloat[high CV ($CV_\Lambda=0.1299$)]{
\includegraphics[scale=0.475]{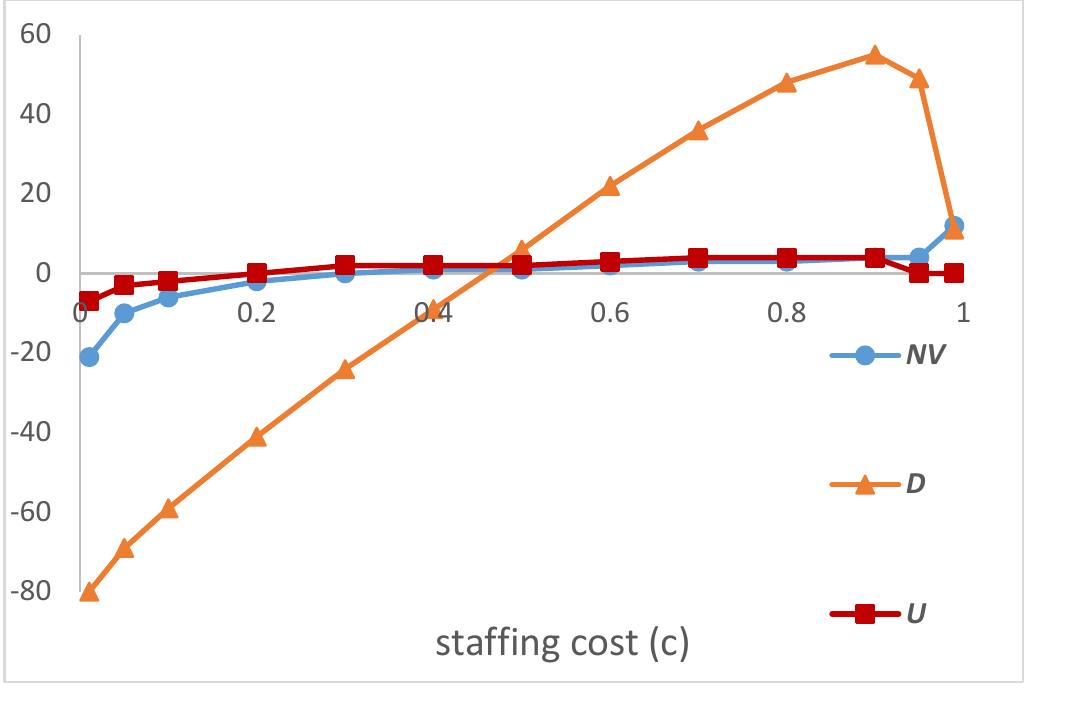}
\label{fig:crhighcv_staff}
}
\caption{staffing error (i.e., $N^{\mbox{\tiny{opt}}} - \triangle$ where $\triangle \in \{ N_U, N_{D}, N_{NV} \}$) for changing staffing costs at three levels of arrival rate uncertainty}
\label{fig: cr_staff}
\end{figure}

Next, we explore the effect of the staffing cost, which determines  the associated critical ratio of the newsvendor policy. To do this, we  change the staffing cost $c$ while holding the other parameters constant. We perform three separate studies by fixing the arrival rate distribution at three separate levels of uncertainty; low CV, moderate CV and high CV. In particular, we assume $\Lambda\sim \mathcal{U}[90,110]$ to produce low CV (Figures \ref{fig:crlowcv_cost} and \ref{fig:crlowcv_staff}), $\Lambda\sim \mathcal{U}[50,150]$ to produce moderate CV (Figures \ref{fig:crmedcv_cost} and \ref{fig:crmedcv_staff}), and $\Lambda\sim \mathcal{U}[10,190]$ to produce high CV (Figures \ref{fig:crhighcv_cost} and \ref{fig:crhighcv_staff}). We plot the percentage cost errors in Figure \ref{fig: cr_cost} and the staffing errors in Figure \ref{fig: cr_staff} for $U$, $D$ and $NV$. We refer the reader to Tables \ref{tab: T2R1CR}-\ref{tab: T2R9CR} in EC for further details (exact costs and staffing levels).

%\begin{figure}[htbp]
%\centering
%\includegraphics[trim=10cm 10cm 10cm 0cm]{T2R1CR.pdf}
%\vspace{-7.5cm}
%\caption{performance increasing staffing cost low CV}
%\label{fig: T2R1CR}
%\end{figure}

We first observe that for all staffing costs and across all levels of CV, $U$ staffs very close to the optimal policy and thus performs extremely well. On the other hand, $NV$ performs poorly when the staffing cost is low although the effect gets less pronounced for higher CV values. This is because $NV$ tends to understaff for low staffing costs when the CV is low. As a result, when $c\ll p<a$, $NV$ incurs higher routing control and abandonment costs. As the in-house staffing cost $c$ increases to $p=1$, $NV$ tends to overstaff, although the adverse effects of overstaffing are not as detrimental. We see that  $D$ performs very poorly as the CV level increases because it fails to capture the effect of randomness. Overall, we see that $U$ is robust and performs well across different critical ratio and CV combinations while the alternative policies can perform arbitrarily poorly.

%\begin{figure}[h*]
%\centering
%\includegraphics[trim=10cm 10cm 10cm 0cm]{T2R5CR.pdf}
%\vspace{-7.5cm}
%\caption{performance}
%\label{fig: T2R5CR}
%\end{figure}

%\begin{figure}[h*]
%\centering
%\includegraphics[trim=10cm 10cm 10cm 0cm]{T2R9CR.pdf}
%\vspace{-7.5cm}
%\caption{performance}
%\label{fig: T2R9CR}
%\end{figure}

%Tables \ref{tab: CV2} and \ref{tab: CV3} present the same numerical study as in Table \ref{tab: CV1} except for different $c$ values so that the critical ratios equal $\frac{r-c}{r}=50\%$ and  $\frac{r-c}{r}=10\%$, respectively. We observe that the other alternative policies also start to perform well for a wide range of parameter settings as the critical ratio goes down. In particular, we see that the $U$ policy no longer clearly outperforms the newsvendor prescription for low CV values. We also see a similar yet less pronounced effect in the comparison between the $U$ policy and the policy that disregards the randomness in $\Lambda$. Finally, we note the $U$ policy continues to perform extremely well with respect to the optimal policy and it is only its relative performance with respect to the two other policies that deteriorates as the critical ratio goes down.

\subsection{Effect of distribution asymmetry}\label{Subsection: varying_skewness}

Our numerical result thus far have assumed symmetric Uniform arrival rate distributions. Next, we generalize our results by considering arrival rate distributions that are asymmetric and follow a Beta distribution to study the effect of skewness on the performance of $U$ and the other policies.  {Specifically, we assume  that $\Lambda\sim Beta(\alpha_1,\alpha_2,\lambda-\underline{b}\sqrt{\lambda},\lambda+\overline{b}\sqrt{\lambda})$, where the first two arguments are the scale parameters of the distribution and the last two arguments are the lower and upper bounds of the support.  We let $\underline{b}$ and $\overline{b}$ be arbitrarily large so that the arrival rate may not realize in the QED regime (i.e., the assumption~\ref{eq:QED} is not necessarily satisfied).  Our proposed policy $U$ is defined for $X\sim Beta(\alpha_1,\alpha_2,\underline{b},\overline{b})$ from (\ref{Eq: UniversalX}).}

We keep the mean arrival rate  fixed at $E[\Lambda]=\lambda=100$ (i.e., $E[X]=0$) throughout this section and we consider three cases where we keep the variance of the arrival rate fixed at three levels: The low CV case keeps the variance of $\Lambda$ fixed and equal to that of a $\mathcal{U}[90,110]$ random variable (i.e., $Var\left(X\right)=Var\left(\mathcal{U}\left[-1,1\right]\right)$),  the moderate CV case keeps the variance of $\Lambda$ fixed and equal to that of a $\mathcal{U}[50,150]$ random variable (i.e., $Var\left(X\right)=Var\left(\mathcal{U}\left[-5,5\right]\right)$), and the high CV case keeps the variance of $\Lambda$ fixed and equal to that of a $\mathcal{U}[10,190]$ random variable (i.e., $Var\left(X\right)=Var\left(\mathcal{U}\left[-9,9\right]\right)$).

\begin{figure}[htb*]
\subfloat[low CV ($CV_\Lambda=0.0144$)]{
\includegraphics[scale=0.425]{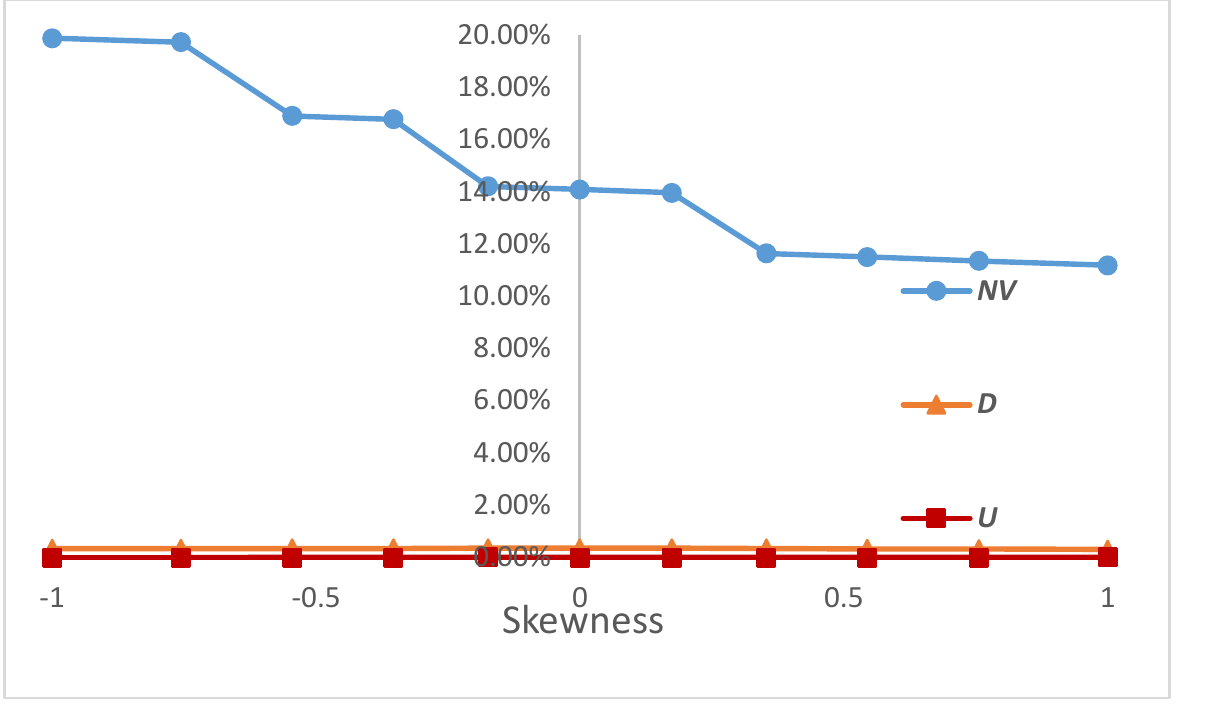}
\label{fig:skewlowcv_cost}
}
\subfloat[moderate CV ($CV_\Lambda=0.0722$)]{
\includegraphics[scale=0.425]{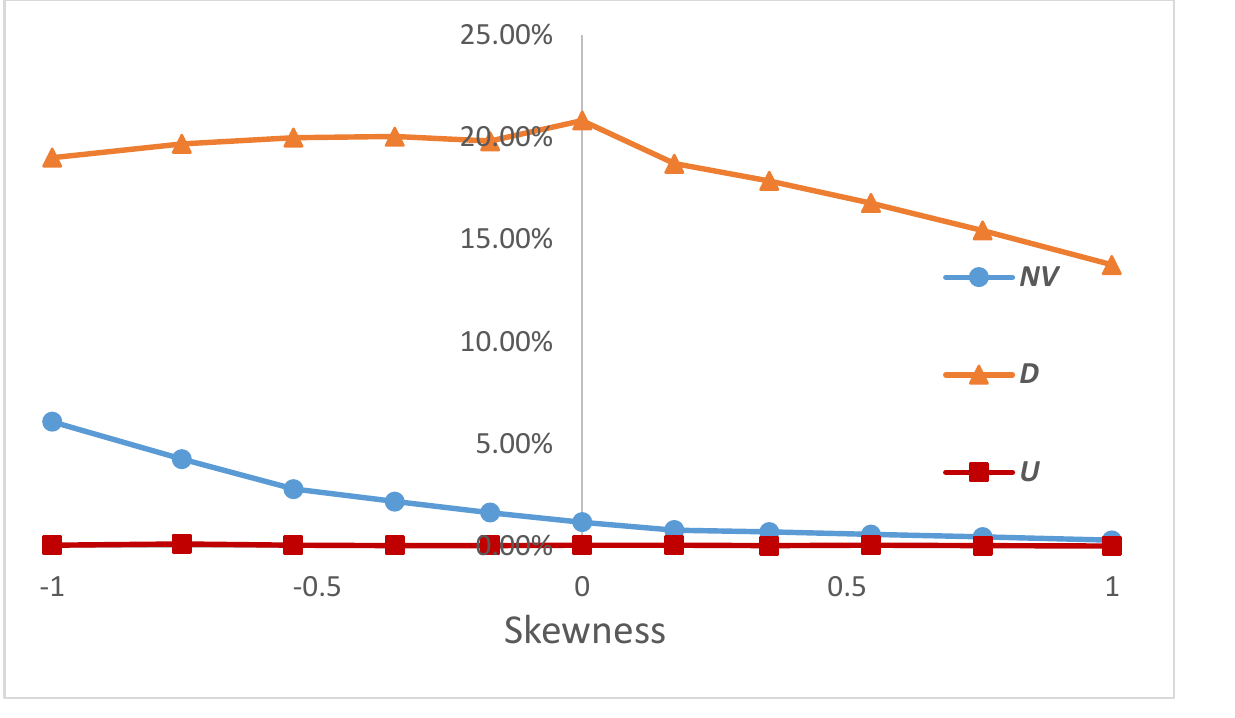}
\label{fig:skewmedcv_cost}
}
\subfloat[high CV ($CV_\Lambda=0.1299$)]{
\includegraphics[scale=0.425]{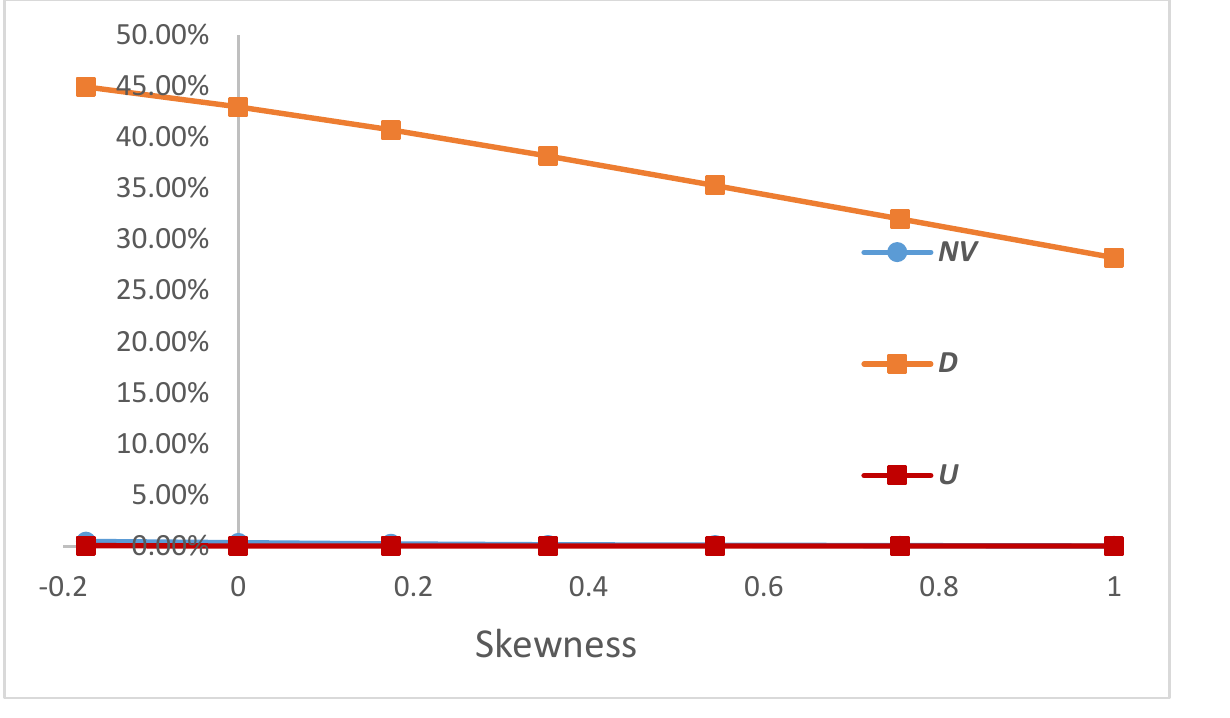}
\label{fig:skewhighcv_cost}
}
\caption{\% cost error for varying skewness levels at three levels of arrival rate uncertainty}
\label{fig: skew_cost}
\end{figure}

\begin{figure}[htb*]
\subfloat[low CV ($CV_\Lambda=0.0144$)]{
\includegraphics[scale=0.425]{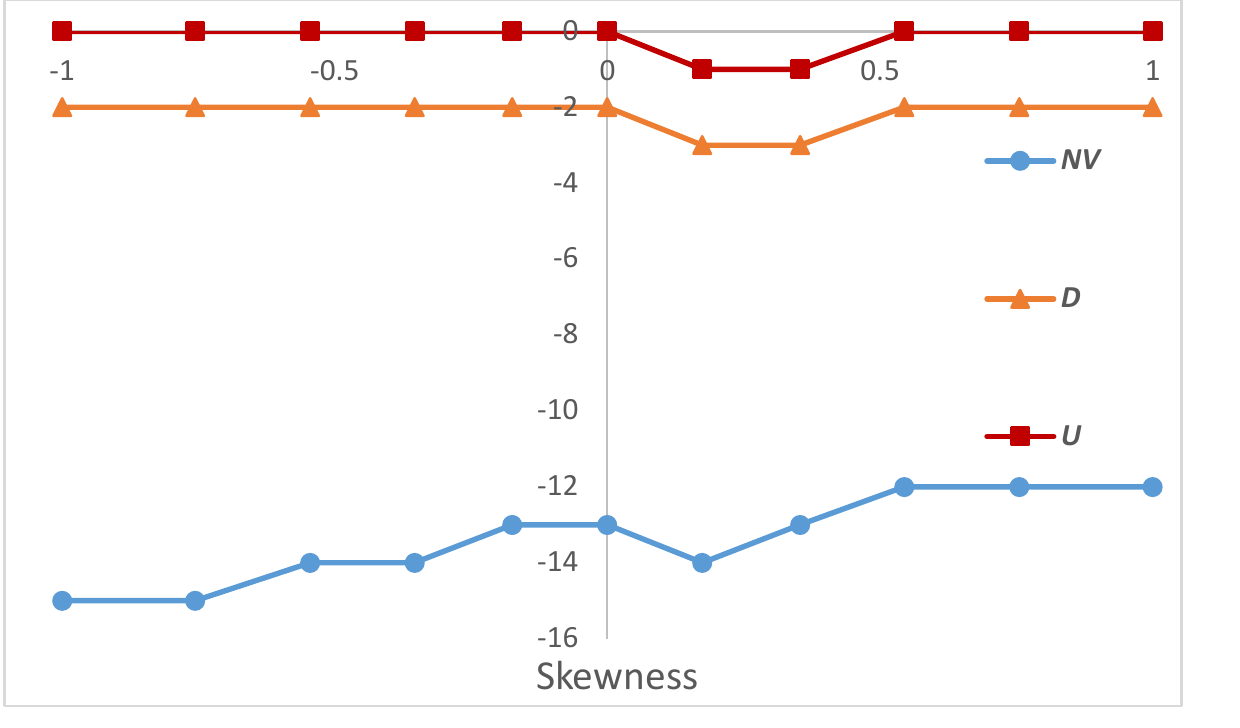}
\label{fig:skewlowcv_staff}
}
\subfloat[moderate CV ($CV_\Lambda=0.0722$)]{
\includegraphics[scale=0.425]{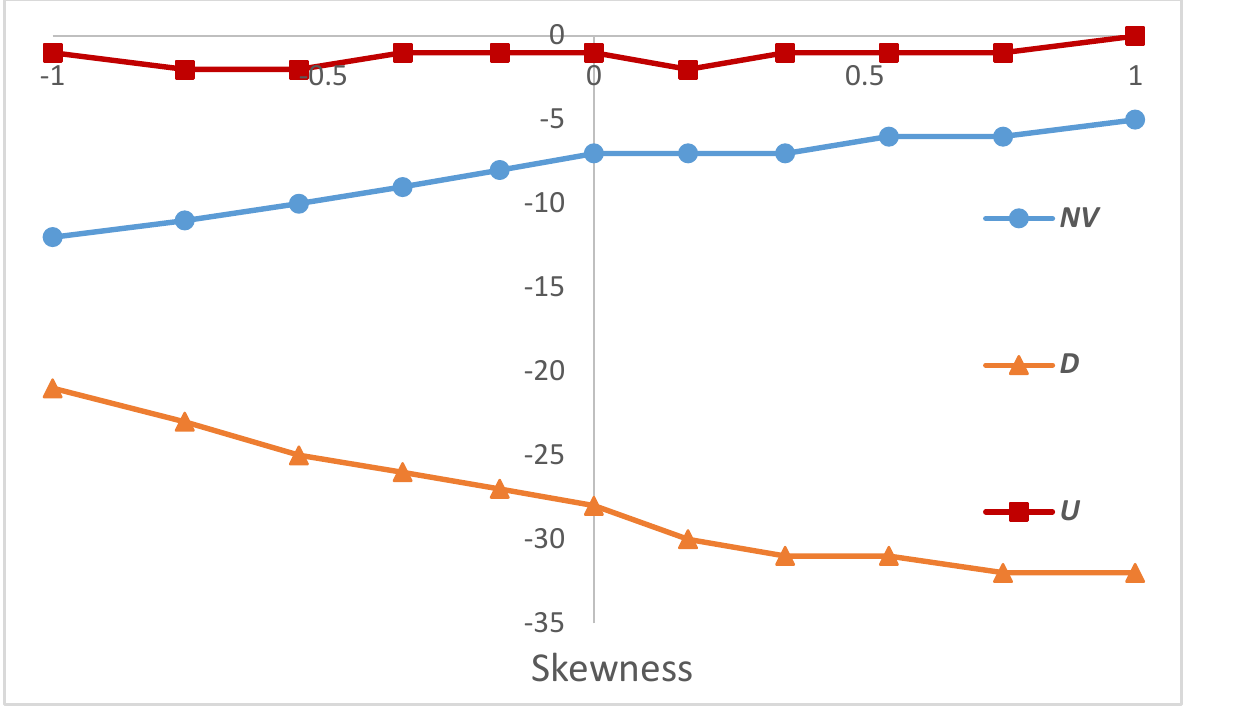}
\label{fig:skewmedcv_staff}
}
\subfloat[high CV ($CV_\Lambda=0.1299$)]{
\includegraphics[scale=0.425]{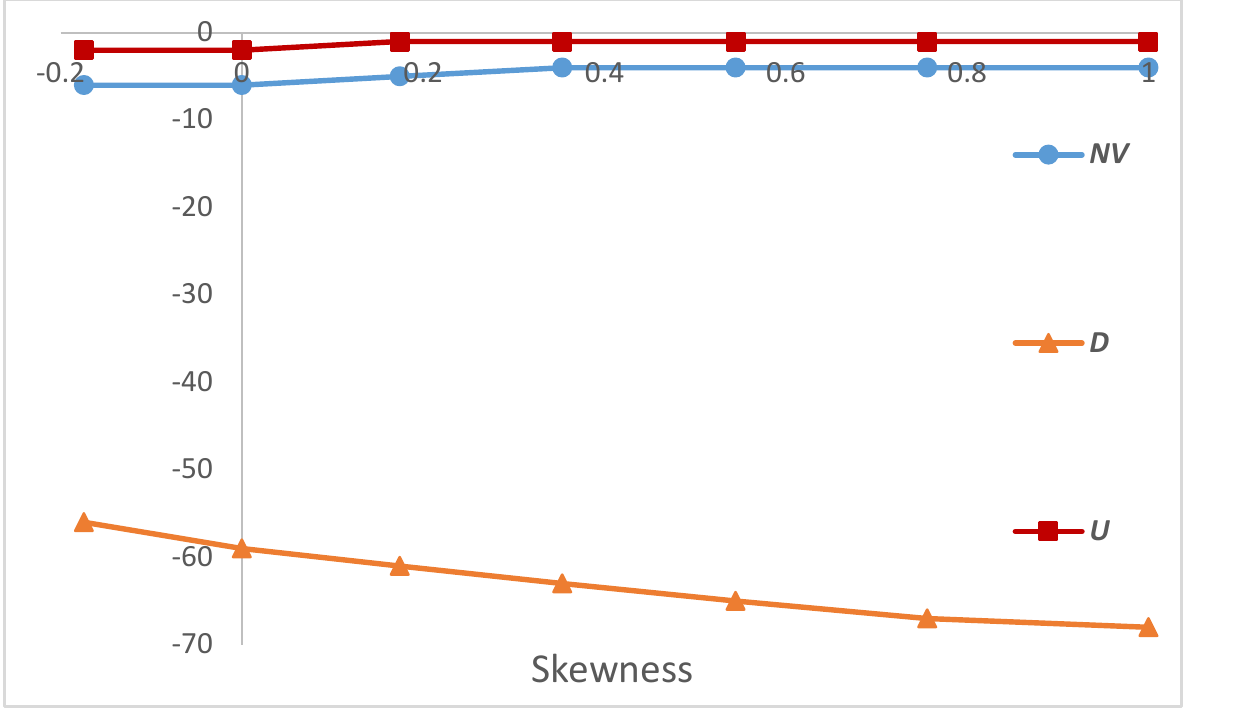}
\label{fig:skewhighcv_staff}
}
\caption{staffing error (i.e., $N^{\mbox{\tiny{opt}}} - \triangle$ where $\triangle \in \{ N_U, N_{D}, N_{NV} \}$) for varying skewness at three levels of arrival rate uncertainty}
\label{fig: skew_staff}
\end{figure}

We study the effect of asymmetry by changing the skewness of the Beta distribution  through its scale parameters $\alpha_1$ and $\alpha_2$. In particular, we set $E[X]=\frac{\alpha_1 \overline{b}+\alpha_2 \underline{b}}{\alpha_1+\alpha_2}=0$ and $Var\left(X\right)=\frac{\alpha_1\alpha_2\left(\overline{b}-\underline{b}\right)^2}{\left(\alpha_1+\alpha_2\right)^2\left(\alpha_1+\alpha_2+1\right)}=\sigma^2$, where $\sigma^2$ denotes the variance of the associated CV level, and we choose\footnote{Note that setting $E[X]=0$ yields $\frac{\overline{b}}{\underline{b}}=-\frac{\alpha_2}{\alpha_1}$, which together with $Var\left(X\right))=\sigma^2$ yields $\underline{b}^2=\frac{\sigma^2\alpha_1(\alpha_1+\alpha_2+1)}{\alpha_2}$. Hence the values of $\alpha_1$, $\alpha_2$, $\underline{b}$, and $\overline{b}$ are not fully determined and we arbitrarily set $\alpha_1+\alpha_2=2$ and change $\alpha_1$ and $\alpha_2$ accordingly, which also changes $\underline{b}$ and $\overline{b}$.} $\alpha_1$ and $\alpha_2$ such that $\alpha_1+\alpha_2=2$. We start with a negative-skewed (left-skewed) Beta distribution with scale parameters $\alpha_1=1.5$ and $\alpha_2=0.5$ with a corresponding skewness of $-1$. Then we decrease $\alpha_1$ and increase $\alpha_2$ so that the skewness of the Beta distribution increases.\footnote{Recall that the skewness of Beta distribution is given by $\frac{2(\alpha_2-\alpha_1)(\sqrt{\alpha_1+\alpha_2+1})}{(\alpha_1+\alpha_2+2)(\sqrt{\alpha_1\alpha_2})}$.} The mid-point where $\alpha_1=\alpha_2=1$ and so the skewness equals 0, corresponds to the symmetric Uniform distribution. After $\alpha_1=\alpha_2=1$, the distribution becomes positive-skewed (right-skewed) as we decrease $\alpha_1$ and increase $\alpha_2$ and we continue until $\alpha_1=0.5$ and $\alpha_2=1.5$ which corresponds to a skewness of $+1$. We plot the percentage cost error and staffing errors of $U$ and the other policies in Figure \ref{fig: skew_cost} and Figure \ref{fig: skew_staff}, respectively. Tables \ref{tab: T2R1Skew}-\ref{tab: T2R9Skew} in EC provide further details (the exact costs and staffing levels as well as the shape parameters of the Beta distribution).

From Figure \ref{fig: skew_cost} and Figure \ref{fig: skew_staff}, we first see that $U$ continues to perform well under asymmetric arrival rate distributions and across varying levels of skewness. In line with our observations for the uniform distribution, $D$ does not perform well except for low CV values. On the other hand, the newsvendor policy performs well for high levels of variability while its performance deteriorates for lower levels of variability and in particular left-skewed distribution. This is because the newsvendor staffing is given by $N_{NV}:=\lambda+\beta_2^\star\sqrt{\lambda}$ for $\beta_2^\star:=F_{X}^{-1}(\frac{p-c}{p})$, which decreases   as the skewness decreases. Hence, the newsvendor staffs less as skewness decreases and the distribution gets more left-skewed, as seen in Figure \ref{fig: skew_staff}. Therefore the newsvendor performs worse when the distribution is left-skewed because its understaffing is more severe, resulting in higher abandonment and routing control costs.

\subsection{Discussion} \label{subsec:discussion}

Our numerical results in Sections \ref{Subsection: finite_system_size} - \ref{Subsection: varying_skewness} show that $U$ is extremely robust, and achieves close to minimum cost over a large range of parameters and assumptions on the amount of the arrival rate uncertainty, and its distribution form. In fact, in virtually all of our experiments $U$ outperformed $D$ and $NV$ and achieved a cost that was very close to the true optimal cost. To better understand the reason for the extremely robust performance of $U$, we compare the actual expected cost to the diffusion approximation of the expected cost, for a wide range of staffing levels. We keep the staffing level $N$ fixed, and approximate the cost using the expression that appears in the limit in Theorem~\ref{Thm: Cost Convergence}.  However, we do not assume that the form of the arrival rate uncertainty is consistent with (\ref{eq:QED}) (as assumed by Theorem~\ref{Thm: Cost Convergence}).  Specifically, for
\[
\beta=\frac{N-\lambda}{\sqrt{\lambda}} \mbox{ and } X=\frac{\Lambda-\lambda}{\sqrt{\lambda}},
\]
it follows that the limiting diffusion cost $c\beta+E \left[ \hat{z}(\beta - X, \hat{T}^\star(\beta - X)) \right]$ when re-scaled gives the following approximation for the actual cost:
\begin{eqnarray}
cN + E \left[ z^{\mbox{\tiny{opt}}}\left(N,\Lambda\right) \right] & \approx & c\lambda ~+~\sqrt{\lambda} \left(c\beta+E \left[ \hat{z}(\beta - X, \hat{T}^\star(\beta - X)) \right]\right) \nonumber \\
& = & cN + \sqrt{\lambda} E \left[ \hat{z}(\beta-X, \hat{T}^\star(\beta - X)) \right]. \label{eq:approx_cost}
\end{eqnarray}

\begin{figure}[htb*]
\centering
\subfloat{
\includegraphics[scale=0.4]{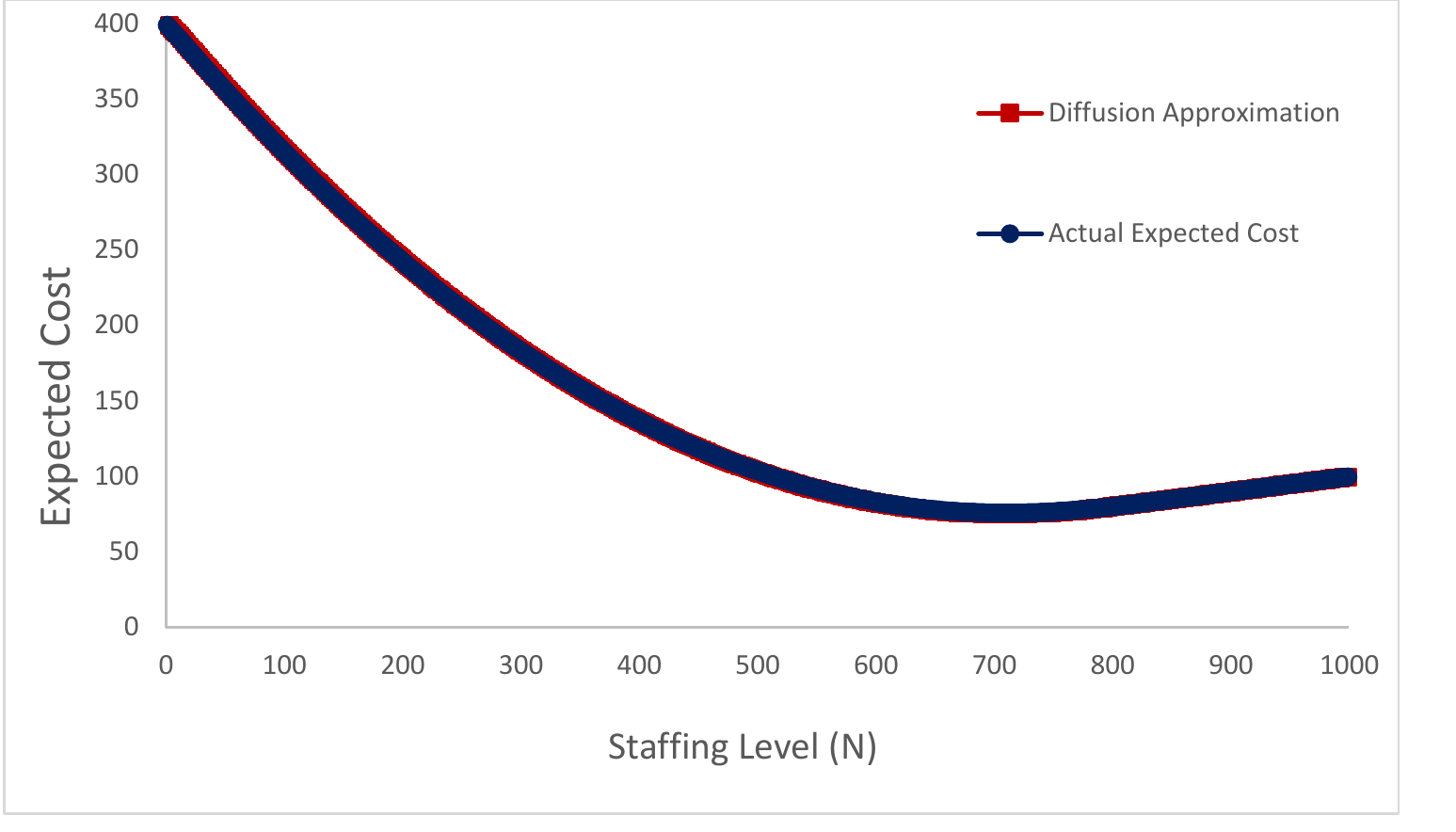}
}
\subfloat{
\includegraphics[scale=0.43]{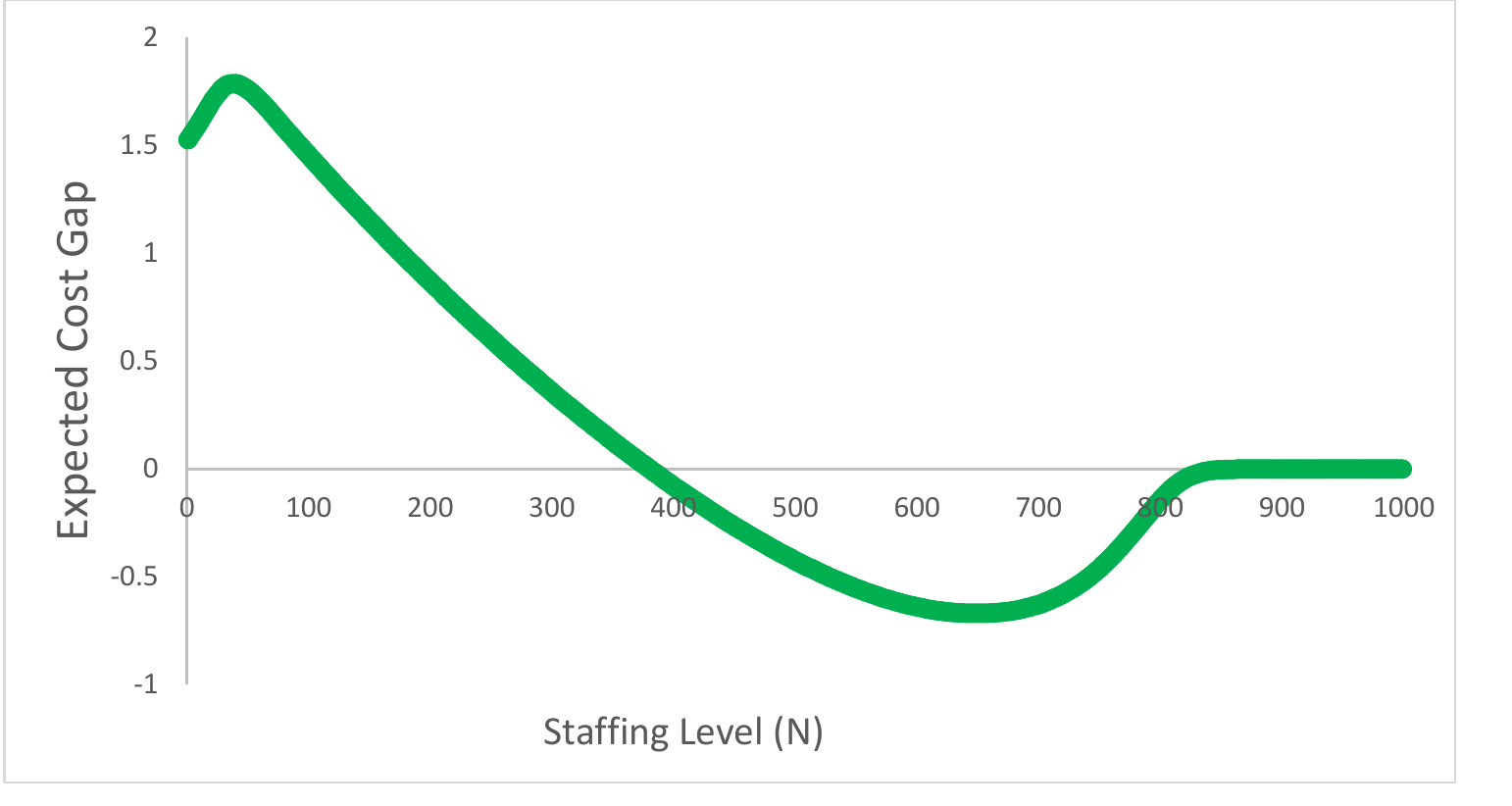}
}
\caption{Expected Cost Difference: Actual vs Diffusion. [$\Lambda\sim\mathcal{U}[20,780]$, $\lambda=400$, $c=0.1$, $p=1$ and $a=5$.]}
\label{fig:expcostdif}
\end{figure}

Figure~\ref{fig:expcostdif} demonstrates numerically that the approximation (\ref{eq:approx_cost}) is very accurate, far beyond what is proven in Theorem~\ref{Thm: Cost Convergence}. Specifically, Figure~\ref{fig:expcostdif} plots the actual expected cost (the left-hand side of (\ref{eq:approx_cost})) and the re-scaled diffusion cost (the right-hand side of (\ref{eq:approx_cost})), and the difference between the two.    It is clear that the approximation in (\ref{eq:approx_cost}) is extremely accurate, over a wide range of staffing levels.  This helps explain the robustness in the performance of $U$.

Figure~\ref{fig:expcostdif} suggests that the performance of square root safety staffing policies can be approximated well without making special assumptions on the limit regime; that is, there is a ``universal'' approximation. This is because we do not restrict $\beta$ values to a particular range which can also be evidenced from the $\beta^{\star}$ values that we observed in our numerical studies in Sections~\ref{Subsection: finite_system_size}-\ref{Subsection: varying_skewness}. In particular,  the $\beta^{\star}$ values in Figures 2-4 come from a wide range from -12 to 10 (see Table \ref{tab: T2R159CR beta} for details). Recalling that $\sqrt{\lambda}=10$ in Figures 2-4, we observe that such extremely low or high values of $\beta^{\star}$ that are essentially of the same order of magnitude as $\sqrt{\lambda}$ (so that the resulting safety staffing is in fact of order $\lambda$) allow us to capture heavily overloaded or underloaded systems, and thus allows us to approximate parameter regimes beyond what is assumed in our asymptotic analysis in Section \ref{Section:aoHW}.

Our observation is also consistent with the universal approximation result of \cite{GHM2012} for a M/M/N+M model with {\em deterministic} arrival rate and {\em no} routing control.  Although it is tempting to think that  \cite{GHM2012} can be used to explain Figure~\ref{fig:expcostdif}, the modeling generalization from  a deterministic to a random arrival rate is not immediate, even if we do not allow for outsourcing.   Our goal is to show that
\begin{equation}
E \left[z_{\tau(\infty)}(N,\Lambda) \right] -\sqrt{\lambda}  E\left[ \hat{z} \left( \frac{N-\Lambda}{\sqrt{\lambda}},\infty \right) \right] = \mathcal{O}(1),
\end{equation}
where $\hat{z}(m,\infty)$ is as defined in Lemma~\ref{lemma:KW2010},  when the limit as $\hat{T} \rightarrow \infty$ is taken.
Under the policy $\tau(\infty)$ that outsources no customers (and so besides staffing only incurs costs through customer abandonment), then $z_{\tau(\infty)}(N,l) = a \gamma \overline{Q}_{\tau(\infty)}(l)$, and so
it follows from Corollary 1 in  \cite{GHM2012} and algebraic manipulation that
\begin{equation} \label{eq:zUniversal}
z_{\tau(\infty)}(N,l) - \sqrt{l}\hat{z} \left( \frac{N-l}{\sqrt{l}},\infty\right) = \mathcal{O}(1).
\end{equation}
Although it seems reasonable that a technical argument would enable us to show that (\ref{eq:zUniversal}) implies
\[
E \left[z_{\tau(\infty)}(N,\Lambda) \right] - E\left[ \sqrt{\Lambda} \hat{z} \left( \frac{N-\Lambda}{\sqrt{\Lambda}},\infty \right) \right] = \mathcal{O}(1),
\]
it is not clear how to establish that
\begin{equation} \label{eq:hard_to_show}
\sqrt{\lambda}E\left[  \hat{z} \left( \frac{N-\Lambda}{\sqrt{\lambda}},\infty \right) \right] - E\left[ \sqrt{\Lambda} \hat{z} \left( \frac{N-\Lambda}{\sqrt{\Lambda}},\infty \right) \right]= \mathcal{O}(1),
\end{equation}
Intuitively, we expect something like (\ref{eq:hard_to_show}) to be true because as the system becomes more understaffed $\hat{z}(\cdot,\infty)$ starts to look linear and as the system becomes more overstaffed $\hat{z}(\cdot,\infty)$ should be negligible (and when the staffing is neither severely understaffed or overstaffed we can appeal to Theorem~\ref{Thm: Cost Convergence}).  However, the analytic argument is not straightforward, and remains as an open question.

%We have tested the performance of $U$} for a wide range of parameters settings. All of our results show that $U$} performs well across all parameter settings establishing that $U$} is robust and performs well in parameter settings beyond which it is asymptotically optimal.On the other, we see that the alternate policies can perform arbitrarily bad. The non-random policy starts to perform poorly even for moderate levels of arrival uncertainty. The newsvendor policy performs poorly for lower levels of arrival uncertainty and especially lower staffing levels (higher critical ratios). Given the extent of outsourcing in call center industry, there should be a non-negligible fraction of call centers that face outsourcing costs that are higher than staffing costs in-house. Although outsourcing costs can be arbitrarily close to staffing costs in-house at first site, hidden costs of outsourcing will actually increase the cost of outsourcing putting many centers in operating setting with higher critical rations. In this case, the simpler newsvendor policy will perform poorly while $U$} will continue to perform very well.

\section{Conclusions}
\label{Section: Conclusion}
Customer call centers face the difficult challenge of having to determine staffing levels in the face of a high level of variability in customer demand. Making these upfront decisions may result in either over- or under-staffing, which implies large unnecessary operating costs. One strategy that call centers use to mitigate these risks is to co-source. Then there is a dynamic decision of whether to handle an arriving call in-house or to outsource it. Our paper deals with the joint staffing and dynamic outsourcing decision for a co-sourced customer call center. The key feature in our model is the time scales differentiation, where staffing decisions are made upfront when the arrival rate is uncertain, while the outsourcing decisions are made dynamically in real-time with full knowledge of the arrival rate realization as well as the system state.

We propose the policy $U$, which staffs the center with the {mean} offered load plus a safety staffing that is of the order of square root of this {mean} offered load. The policy also outsources calls whenever the number of customers in line exceeds a threshold. When the order of magnitude of arrival rate randomness is the same order as the inherent system fluctuations in the queue length {(which is on the order of the square root of the realized arrival rate)}, we show that our proposed policy is asymptotically optimal as the mean arrival rate grows large. Then, we perform an extensive numerical experiment to study the performance of $U$ beyond the regime for which it is proved to be asymptotically optimal. In all of our numerical experiments {in which the system has more than a few servers} we did not encounter even one case in which the performance of $U$ was not superb. In contrast, our two benchmark policies each have parameter regimes in which their performance can be arbitrarily bad. It is, therefore, the {\em robustness} of the performance of $U$ that we would like to stress as the main takeaway from this paper. One does not need to identify the ``right" operating regime in order to determine which policy to use. The $U$ policy appears to be a ``one-size-fits-all".

We further observe from our numerical results that the newsvendor policy performs well except when the in-house staffing cost $c$ is low relative to the costs of outsourcing $p$ and abandonment $a$. Since the newsvendor policy is a simpler policy than $U$, the question arises: can we use the newsvendor instead of $U$? Although the answer is ``yes'' for a non-trivial portion of the parameter space, the answer is ``no'' when considering the entire parameter space. This is important because outsourcing costs can be very high when the hidden costs of outsourcing are included \citep{Kharif2003} and \citep{Rubin2013}. A recent study \citep{InfoTech2011} shows that ``hidden costs'' of outsourcing add, on average, 25\% to the overall cost showing that the hidden costs of outsourcing can indeed be significant.

Several important extensions are worth pursuing. In this paper, the recourse action we consider, in case the center is under-staffed is outsourcing some calls to an outside vendor (other actions have been considered in the literature such as utilizing ``on-call" agents, \cite{MehrotraOzlukSaltzman2010}, \cite{Gans12}, or by sharing resources with other sites~\cite{PerryWhitt09}). On the flip side, one might consider recourse actions for scenarios in which the call center ends up being over-staffed. Those may include sending agents home (as in \cite{MehrotraOzlukSaltzman2010}, \cite{Gans12}) or, if this is not desirable or feasible, diverting agents to other activities such as making outbound calls (\cite{GansZhou2007,BhulaiKoole03}) or cross-selling (\cite{GurvichArmonyMaglaras09,ArmonyGurvich10}). It would be interesting to examine a framework in which recourse actions in both directions may be used.

Another interesting extension is to model the staffing decisions of the outside vendor more explicitly, as in \cite{GurvichPerry11} in the case of known arrival rate. In our model we assume that a call outsourced incurs a fixed cost no matter how many calls are sent to the outside vendor. This is consistent with the assumption made in \cite{Aksin2008}, and is equivalent to assuming either that the outside service provider has ample service capacity or that it pools demand from a large enough client base so that the calls the company sends to the outside vendor do not have much impact. For the case where the outsourcing is not preferred or does not exist, we see that the performance of $U$ remains superb. (We do not report these results due to space limitations.) However, if the outsourcing capacity is positive and limited we do not know how $U$ will perform, or how it would need to be modified to maintain this superior performance.

The aforementioned performance deterioration may be mitigated by designing contracts that ensure the outside vendor has ample capacity.  But is such a design optimal?  It would be interesting to incorporate arrival rate uncertainty when studying contract design questions such as those studied in \cite{Aksin2008}, and \cite{ZhouRen2011}.

Beyond the co-sourcing application considered in this paper, one might use a similar framework to examine joint staffing and control decisions with other system topologies  and other types of control such as the problems considered in~\cite{GurvichArmonyMandelbaum08},~\cite{DaiTezcan08},~\cite{TezcanDai10},~\cite{GurvichWhitt10}, and~\cite{ArmonyMandelbaum11}. In particular, it will sometimes be possible to incorporate arrival rate uncertainty without changing the control that has been proven to be asymptotically optimal when the arrival rate is known. Ideally, the robustness of the policy performance with respect to the assumptions on the arrival rate uncertainty in this model will be true in much more generality. Skill-based routing is one particularly interesting direction, where not only does one have to determine the size of multiple pools of agents under arrival rate uncertainty, but also, once the arrival rate is known, dynamic control is used to determine the assignment of calls to agents. This setting raises another interesting modeling question that has to do with how to model the arrival rate uncertainty in the presence of multiple streams of arrivals (\cite{ShenHuang08}, \cite{Ibrahim13}).

Finally, our numerical results with respect to the remarkable robustness of $U$  suggest that there might be an underlying theoretical justification to this robustness, in the spirit of the universal approximation of \cite{GHM2012}. In Section \ref{subsec:discussion} we have discussed why their results in a framework that assumes a known arrival rate and no dynamic control may not be readily applied to our framework. But one wonders whether similar universal approximation principles apply when a random arrival rate and dynamic control
are incorporated into the model.

% Acknowledgments here
%\ACKNOWLEDGMENT{AMY:  Do we want to add any acknowledgements?  Maybe thank the Anonymous referees?}

\bibliographystyle{ormsv080}
\bibliography{AKW_04_07_2014}

\begin{thebibliography}{49}
\expandafter\ifx\csname natexlab\endcsname\relax\def\natexlab#1{#1}\fi
\expandafter\ifx\csname url\endcsname\relax
  \def\url#1{{\tt #1}}\fi
\expandafter\ifx\csname urlprefix\endcsname\relax\def\urlprefix{URL }\fi
\expandafter\ifx\csname urlstyle\endcsname\relax
  \expandafter\ifx\csname doi\endcsname\relax
  \def\doi#1{doi:\discretionary{}{}{}#1}\fi \else
  \expandafter\ifx\csname doi\endcsname\relax
  \def\doi{doi:\discretionary{}{}{}\begingroup \urlstyle{rm}\Url}\fi \fi

\bibitem[{Aksin et~al.(2007)Aksin, Armony, and Mehrotra}]{Aksin07}
Aksin, Z., M.~Armony, V.~Mehrotra. 2007.
\newblock The modern call center: A multi-disciplinary perspective on
  operations management research.
\newblock {\it Production and Operations Management\/} {\bf 16}(6) 665--688.

\bibitem[{Ak{\c{s}}in et~al.(2008)Ak{\c{s}}in, de~V{\'e}ricourt, and
  Karaesmen}]{Aksin2008}
Ak{\c{s}}in, Z., F.~de~V{\'e}ricourt, F.~Karaesmen. 2008.
\newblock Call center outsourcing contract analysis and choice.
\newblock {\it Management Science\/} {\bf 54}(2) 354--368.

\bibitem[{Armony and Gurvich(2010)}]{ArmonyGurvich10}
Armony, M., I.~Gurvich. 2010.
\newblock When promotions meet operations: Cross selling and its effect on
  call-center performance.
\newblock {\it Manufacturing \& Service Operations Management\/} {\bf 12}(3)
  470--488.

\bibitem[{Armony and Mandelbaum(2011)}]{ArmonyMandelbaum11}
Armony, M., A.~Mandelbaum. 2011.
\newblock Routing and staffing in large-scale service systems: The case of
  homogeneous impatient customers and heterogeneous servers.
\newblock {\it Operations research\/} {\bf 59}(1) 50--65.

\bibitem[{Armony et~al.(2009)Armony, Plambeck, and Seshadri}]{Armony09}
Armony, M., E.~Plambeck, S.~Seshadri. 2009.
\newblock Sensitivity of optimal capacity to customer impatience in an
  unobservable m/m/s queue (why you shouldn't shout at the dmv).
\newblock {\it Manufacturing \& Service Operations Management\/} {\bf 11}(1)
  19--32.

\bibitem[{Bassamboo et~al.(2005)Bassamboo, Harrison, and
  Zeevi}]{BassambooHarrisonZeevi05}
Bassamboo, A., J.M. Harrison, A.~Zeevi. 2005.
\newblock Dynamic routing and admission control in high-volume service systems:
  Asymptotic analysis via multi-scale fluid limits.
\newblock {\it Queueing Systems\/} {\bf 51}(3) 249--285.

\bibitem[{Bassamboo et~al.(2006)Bassamboo, Harrison, and
  Zeevi}]{BassambooHarrisonZeevi06}
Bassamboo, A., J.M. Harrison, A.~Zeevi. 2006.
\newblock {Design and Control of a Large Call Center: Asymptotic Analysis of an
  LP-Based Method}.
\newblock {\it Operations Research\/} {\bf 54}(3) 419--435.

\bibitem[{Bassamboo et~al.(2010)Bassamboo, Randhawa, and
  Zeevi}]{BassambooRandhawaZeevi10}
Bassamboo, A., R.S. Randhawa, A.~Zeevi. 2010.
\newblock Capacity sizing under parameter uncertainty: Safety staffing
  principles revisited.
\newblock {\it Management Science\/} {\bf 56}(10) 1668--1686.

\bibitem[{Bassamboo and Zeevi(2009)}]{BassambooZeevi09}
Bassamboo, A., A.~Zeevi. 2009.
\newblock {On a Data-Driven Method for Staffing Large Call Centers}.
\newblock {\it Operations Research\/} {\bf 57}(3) 714--726.

\bibitem[{Bhulai and Koole(2003)}]{BhulaiKoole03}
Bhulai, S., G.~Koole. 2003.
\newblock A queueing model for call blending in call centers.
\newblock {\it IEEE transactions on automatic control\/} {\bf 48} 1434--1438.

\bibitem[{Borst et~al.(2004)Borst, Mandelbaum, and Reiman}]{Borst04}
Borst, S., A.~Mandelbaum, M.I. Reiman. 2004.
\newblock {Dimensioning large call centers}.
\newblock {\it Operations research\/} {\bf 52}(1) 17--34.

\bibitem[{Brown et~al.(2005)Brown, Gans, Mandelbaum, Sakov, Shen, Zeltyn, and
  Zhao}]{Brown05}
Brown, L., N.~Gans, A.~Mandelbaum, A.~Sakov, H.~Shen, S.~Zeltyn, L.~Zhao. 2005.
\newblock {Statistical Analysis of a Telephone Call Center}.
\newblock {\it Journal of the American Statistical Association\/} {\bf
  100}(469) 36--50.

\bibitem[{Chen and Henderson(2001)}]{Chen01}
Chen, B.P.K., S.G. Henderson. 2001.
\newblock {Two issues in setting call centre staffing levels}.
\newblock {\it Annals of Operations Research\/} {\bf 108}(1) 175--192.

\bibitem[{Dai and Tezcan(2008)}]{DaiTezcan08}
Dai, J.G., T.~Tezcan. 2008.
\newblock Optimal control of parallel server systems with many servers in heavy
  traffic.
\newblock {\it Queueing Systems\/} {\bf 59} 95--134.

\bibitem[{Gans et~al.(2003)Gans, Koole, and Mandelbaum}]{Gans03}
Gans, N., G.~Koole, A.~Mandelbaum. 2003.
\newblock {Telephone call centers: Tutorial, review, and research prospects}.
\newblock {\it Manufacturing \& Service Operations Management\/} {\bf 5}(2)
  79--141.

\bibitem[{Gans et~al.(2012)Gans, Shen, Zhou, Korolev, McCord, and
  Ristock}]{Gans12}
Gans, N., H.~Shen, Y-P Zhou, K.~Korolev, A.~McCord, H.~Ristock. 2012.
\newblock Parametric stochastic programming models for call-center workforce
  scheduling.
\newblock Working Paper.

\bibitem[{Gans and Zhou(2007)}]{GansZhou2007}
Gans, N., Y.P. Zhou. 2007.
\newblock {Call-routing schemes for call-center outsourcing}.
\newblock {\it Manufacturing \& Service Operations Management\/} {\bf 9}(1)
  33--50.

\bibitem[{Garnett et~al.(2002)Garnett, Mandelbaum, and
  Reiman}]{GarnettMandelbaumReiman02}
Garnett, O., A.~Mandelbaum, M.~Reiman. 2002.
\newblock {Designing a call center with impatient customers}.
\newblock {\it Manufacturing \& Service Operations Management\/} {\bf 4}(3)
  208--227.

\bibitem[{Gilson and Khandelwal(2005)}]{GiKh2005}
Gilson, K.~A., D.~K. Khandelwal. 2005.
\newblock {Getting more from call centers. The McKinsey Quaterly}.

\bibitem[{Green et~al.(2001)Green, Kolesar, and Soares}]{Green_etal_2001}
Green, L.V., P.J. Kolesar, J.~Soares. 2001.
\newblock Improving the sipp approach for staffing service systems that have
  cyclic demands.
\newblock {\it Operations Research\/} {\bf 49}(4) 549--564.

\bibitem[{Gurvich et~al.(2000)Gurvich, Armony, and
  Maglaras}]{GurvichArmonyMaglaras09}
Gurvich, I., M.~Armony, C.~Maglaras. 2000.
\newblock Cross-selling in a call center with a heterogeneous customer
  population.
\newblock {\it Operations Research\/} {\bf 57}(2) 299--313.

\bibitem[{Gurvich et~al.(2008)Gurvich, Armony, and
  Mandelbaum}]{GurvichArmonyMandelbaum08}
Gurvich, I., M.~Armony, A.~Mandelbaum. 2008.
\newblock Service-level differentiation in call centers with fully flexible
  servers.
\newblock {\it Management Science\/} {\bf 54}(2) 279--294.

\bibitem[{Gurvich et~al.(2013)Gurvich, Huang, and Mandelbaum}]{GHM2012}
Gurvich, I., J.~Huang, A.~Mandelbaum. 2013.
\newblock Excursion-based universal approximations for the erlang-a queue in
  steady-state.
\newblock {\it Mathematics of Operations Research.\/} Forthcoming.

\bibitem[{Gurvich et~al.(2010)Gurvich, Luedtke, and
  Tezcan}]{GurvichLuedtkeTezcan09}
Gurvich, I., J.~Luedtke, T.~Tezcan. 2010.
\newblock Staffing call centers with uncertain demand forecasts: A
  chance-constrained optimization approach.
\newblock {\it Management Science\/} {\bf 56}(7) 1093--1115.

\bibitem[{Gurvich and Perry(2012)}]{GurvichPerry11}
Gurvich, I., O.~Perry. 2012.
\newblock Overflow networks: Approximations and implications to call-center
  outsourcing.
\newblock {\it Operations Research\/} {\bf 60}(4) 996--1009.

\bibitem[{Gurvich and Whitt(2010)}]{GurvichWhitt10}
Gurvich, I., W.~Whitt. 2010.
\newblock Service-level differentiation in many-server service systems via
  queue-ratio routing.
\newblock {\it Operations Research\/} {\bf 58}(2) 316--328.

\bibitem[{Halfin and Whitt(1981)}]{HalfinWhitt81}
Halfin, S., W.~Whitt. 1981.
\newblock {Heavy-traffic limits for queues with many exponential servers}.
\newblock {\it Operations research\/} {\bf 29}(3) 567--588.

\bibitem[{Ibrahim and L'Ecuyer(2013)}]{Ibrahim13}
Ibrahim, R., P.~L'Ecuyer. 2013.
\newblock Forecasting call center arrivals: Fixed-effects, mixed-effects, and
  bivariate models.
\newblock {\it Manufacturing \& Service Operations Management\/} {\bf 15}(1)
  72--85.

\bibitem[{ICMI(2006)}]{ICMI2006}
ICMI. 2006.
\newblock {ICMI's Contact Center Outsourcing Report}.
\newblock
  \url{http://www.icmi.com/files/ICMI/members/ccmr/ccmr2006/ccmr06/June2006_issue.pdf}.

\bibitem[{InfoTech(2011)}]{InfoTech2011}
InfoTech. 2011.
\newblock Discover the hidden costs of outsourcing.
\newblock
  \url{http://www.infotech.com/research/ss/it-discover-the-hidden-costs-of-outsourcing}.

\bibitem[{Jongbloed and Koole(2001)}]{Jongbloed01}
Jongbloed, G., G.~Koole. 2001.
\newblock {Managing uncertainty in call centres using Poisson mixtures}.
\newblock {\it Applied Stochastic Models in Business and Industry\/} {\bf
  17}(4) 307--318.

\bibitem[{Kharif(2003)}]{Kharif2003}
Kharif, O. 2003.
\newblock The hidden costs of it outsourcing.
\newblock {\it BusinessWeek\/}  (October 26 2003).

\bibitem[{Ko{\c{c}}a\u{g}a and Ward(2010)}]{KocagaWard10}
Ko{\c{c}}a\u{g}a, Y.L., A.R. Ward. 2010.
\newblock {Admission Control for a Multi-Server Queue with Abandonment}.
\newblock {\it Queueing Systems\/} {\bf 65}(3) 275--323.

\bibitem[{Koole and Pot(2011)}]{KoolePot11}
Koole, G., A.~Pot. 2011.
\newblock A note on profit maximization and monotonicity for inbound call
  centers.
\newblock {\it Operations Research\/} {\bf 59}(5) 1304--1308.

\bibitem[{Liu and Whitt(2012)}]{LiuWhitt2012}
Liu, Y., W.~Whitt. 2012.
\newblock Stabilizing customer abandonment in many-server queues with
  time-varying arrivals.
\newblock {\it Operations Research\/} {\bf 60}(6) 1551--1564.

\bibitem[{Maman(2009)}]{Maman09}
Maman, S. 2009.
\newblock {Uncertainty in the demand for service: The case of call centers and
  emergency departments}.
\newblock Master's thesis, Technion-Israel Institute of Technology.

\bibitem[{Mandelbaum and Zeltyn(2009)}]{MandelbaumZeltyn2009}
Mandelbaum, A., S.~Zeltyn. 2009.
\newblock Staffing many-server queues with impatient customers: Constraint
  satisfaction in call centers.
\newblock {\it Operations Research\/} {\bf 57}(5) 1189--1205.

\bibitem[{Mehrotra et~al.(2010)Mehrotra, Ozl{\"u}k, and
  Saltzman}]{MehrotraOzlukSaltzman2010}
Mehrotra, V., O.~Ozl{\"u}k, R.~Saltzman. 2010.
\newblock Intelligent procedures for intra-day updating of call center agent
  schedules.
\newblock {\it Production and Operations Management\/} {\bf 19}(3) 353--367.

\bibitem[{Perry and Whitt(2009)}]{PerryWhitt09}
Perry, O., W.~Whitt. 2009.
\newblock Responding to unexpected overloads in large-scale service systems.
\newblock {\it Management Science\/} {\bf 55}(8) 1353--1367.

\bibitem[{Robbins and Harrison(2010)}]{RobbinsHarrison2010}
Robbins, T.R., T.P. Harrison. 2010.
\newblock A stochastic programming model for scheduling call centers with
  global service level agreements.
\newblock {\it European Journal of Operational Research\/} {\bf 207}(3)
  1608--1619.

\bibitem[{Ross(2001)}]{Ross01}
Ross, A.M. 2001.
\newblock {Queueing systems with daily cycles and stochastic demand with
  uncertain parameters}.
\newblock Ph.D. thesis, University of California, Berkeley.

\bibitem[{Rubin(2013)}]{Rubin2013}
Rubin, J. 2013.
\newblock {Hidden Costs of Outsourcing. {\it Forbes} (March 29)}.
\newblock
  \url{http://www.forbes.com/sites/forbesinsights/2013/03/29/the-hidden-costs-of-outsourcing/}.

\bibitem[{Shen and Huang(2008)}]{ShenHuang08}
Shen, H., J.Z. Huang. 2008.
\newblock Interday forecasting and intraday updating of call center arrivals.
\newblock {\it Manufacturing \& Service Operations Management\/} {\bf 10}(3)
  391--410.

\bibitem[{Steckley et~al.(2009)Steckley, Henderson, and
  Mehrotra}]{Steckley2009}
Steckley, S.G., S.G. Henderson, V.~Mehrotra. 2009.
\newblock {Forecast errors in service systems}.
\newblock {\it Probability in the Engineering and Informational Sciences\/}
  {\bf 23}(2) 305--332.

\bibitem[{Tezcan and Dai(2010)}]{TezcanDai10}
Tezcan, T., J.G. Dai. 2010.
\newblock Dynamic control of n-systems with many servers: asymptotic optimality
  of a static priority policy in heavy traffic.
\newblock {\it Operations Research\/} {\bf 58} 94--110.

\bibitem[{Whitt(2005)}]{Whitt05}
Whitt, W. 2005.
\newblock {Engineering Solution of a Basic Call-Center Model}.
\newblock {\it Management Science\/} {\bf 51}(2) 221--235.

\bibitem[{Whitt(2006)}]{Whitt06}
Whitt, W. 2006.
\newblock {Stafffng a Call Center with Uncertain Arrival Rate and Absenteeism}.
\newblock {\it Production and Operations Management\/} {\bf 15}(1) 88--102.

\bibitem[{Zan et~al.(2013)Zan, Hasenbein, and Morton}]{ZanHasenbeinMorton2013}
Zan, J., J.J. Hasenbein, D.P. Morton. 2013.
\newblock Staffing large service systems under arrival-rate uncertainty.
\newblock {\it arXiv preprint arXiv:1304.6701\/} .

\bibitem[{Zhou and Ren(2011)}]{ZhouRen2011}
Zhou, Y-P., Z.J. Ren. 2011.
\newblock Service outsourcing.
\newblock {\it Wiley Encyclopedia of Operations Research and Management
  Science\/} .

\end{thebibliography}


\begin{thebibliography}{5}
\expandafter\ifx\csname natexlab\endcsname\relax\def\natexlab#1{#1}\fi
\expandafter\ifx\csname url\endcsname\relax
  \def\url#1{{\tt #1}}\fi
\expandafter\ifx\csname urlprefix\endcsname\relax\def\urlprefix{URL }\fi
\expandafter\ifx\csname urlstyle\endcsname\relax
  \expandafter\ifx\csname doi\endcsname\relax
  \def\doi#1{doi:\discretionary{}{}{}#1}\fi \else
  \expandafter\ifx\csname doi\endcsname\relax
  \def\doi{doi:\discretionary{}{}{}\begingroup \urlstyle{rm}\Url}\fi \fi

\bibitem[{Abramowitz and Stegun(1964)}]{AbramowitzStegun64}
Abramowitz, M., I.A. Stegun. 1964.
\newblock {\it Handbook of mathematical functions with formulas, graphs, and
  mathematical tables\/}.
\newblock Dover publications.

\bibitem[{Garnett et~al.(2002)Garnett, Mandelbaum, and
  Reiman}]{GarnettMandelbaumReiman02EC}
Garnett, O., A.~Mandelbaum, M.~Reiman. 2002.
\newblock {Designing a call center with impatient customers}.
\newblock {\it Manufacturing \& Service Operations Management\/} {\bf 4}(3)
  208--227.

\bibitem[{Ko{\c{c}}a\u{g}a and Ward(2010)}]{KocagaWard10EC}
Ko{\c{c}}a\u{g}a, Y.L., A.R. Ward. 2010.
\newblock {Admission Control for a Multi-Server Queue with Abandonment}.
\newblock {\it Queueing Systems\/} {\bf 65}(3) 275--323.

\bibitem[{Smith and Whitt(1981)}]{SmithWhitt1981}
Smith, D.~R., W.~Whitt. 1981.
\newblock {Resource sharing for efficiency in traffic systems}.
\newblock {\it Bell Systems Tech. J.\/} {\bf 60} 39--55.

\bibitem[{Whitt(1984)}]{Whitt84}
Whitt, W. 1984.
\newblock Heavy-traffic approximations for service systems with blocking.
\newblock {\it AT\&T Bell Laboratories Technical Journal\/} {\bf 63}(5)
  689--708.

\end{thebibliography}
\newpage
\appendix{\pagenumbering{arabic}}

\begin{center}
{\LARGE \textbf{Electronic Companion to ``Staffing Call Centers with Uncertain Arrival Rates and Co-sourcing'' }}\\[12pt]
\end{center}
\bigskip
In this electronic companion, we provide supporting material for our manuscript titled ``Staffing Call Centers with Uncertain Arrival Rates and Co-sourcing''. We first provide the proofs of the Theorems, Propositions, and Lemmas in the main paper body (Section~\ref{Section: Proofs}) and then provide additional numerical results to support the figures in the main paper body (Section~\ref{Section:AdditionalNumerics}).

\section{Proofs}\label{Section: Proofs}
\noindent{\newline \bf Proof of Proposition~\ref{prop: prelimit optimal}:}
We first prove part (i) and then prove part (ii). \\
\noindent{\newline \bf Proof of part (i):} It is helpful to first obtain a lower bound on the cost of any staffing policy, regardless of the routing policy $\pi \in \Pi$ that is followed.  For this, let $\overline{S}_\pi(l)$ be the expected number of busy servers when the arrival rate realization is $l$, and let $\overline{S}_\pi(\Lambda)$ be the associated random variable.  Then, since $N \geq \overline{S}_\pi(l)$,
\[
C(N,\pi)  \geq  c E \left[ \overline{S}_\pi(\Lambda) \right] + p E\left[ \Lambda P_\pi(\mbox{out};\Lambda) \right] + a E\left[ \Lambda P_\pi(\mbox{ab};\Lambda) \right],
\]
so that the assumption $c \geq \min(a,p)$ then yields
\begin{equation} \label{Eq: lower_bound_useful}
C(N,\pi)  \geq  \min(a,p) E \left[ \overline{S}_\pi(\Lambda) + \Lambda \left( P_\pi(\mbox{out};\Lambda) + P_\pi(\mbox{ab};\Lambda)  \right) \right].
\end{equation}
For any realization $l$ of $\Lambda$, the arrival rate into the system must equal the departure rate from the system (due to both abandonments and service completion), so that
\[
l(1-P_\pi(\mbox{out};l)) = \overline{S}_\pi(l) + lP_\pi(\mbox{ab};l),
\]
or, equivalently,
\[
l = \overline{S}_\pi(l) + l \left( P_\pi(\mbox{out};l) + P_\pi(\mbox{ab};l) \right).
\]
Taking expectations and recalling that $E[\Lambda]=\lambda$ shows that
\begin{equation} \label{Eq: Balance}
\lambda = E \left[ \overline{S}_\pi(\Lambda) + \Lambda \left( P_\pi(\mbox{out};\Lambda) + P_\pi(\mbox{ab};\Lambda)  \right)   \right].
\end{equation}
Substituting the equality (\ref{Eq: Balance}) into (\ref{Eq: lower_bound_useful}) shows
\[
C(N,\Lambda) \geq \min(a,p) \lambda.
\]

Suppose $a\leq p$ ($a>p$).  Then, the minimum in (\ref{Eq: lower_bound_useful}) is attained by setting a zero staffing level and letting everyone abandon (outsourcing everyone), because the cost associated with that policy is $C(0,\Lambda) = a\lambda$ ($C(0,\Lambda) = p\lambda$).  Hence $N^{\mbox{\tiny{opt}}} = 0$ and
\begin{eqnarray*}
\pi^{\mbox{\tiny{opt}}}(0,\Lambda) & = & (0,0,\ldots,0) \mbox{ if } a >p \\
\pi^{\mbox{\tiny{opt}}}(0,\Lambda) & = & (1,1,\ldots,1) \mbox{ if } a \leq p.
\end{eqnarray*}

\noindent{\newline \bf Proof of part (ii):}  We must show that
\[
E[z_\pi(N,\Lambda)] = p E\left[ \Lambda P_\pi( \mbox{out}; \Lambda) \right] + a E \left[ \Lambda P_\pi(\mbox{ab}; \Lambda) \right].
\]
 is lower bounded by the cost associated with the policy $\tau(\infty)$ that does not outsource any calls (has $P_{\tau(\infty)}(\mbox{out};\Lambda)=0$); i.e., we must show that
\[
p E\left[ \Lambda P_\pi(\mbox{out};\Lambda) \right] + a E \left[ \Lambda P_\pi(\mbox{ab}; \Lambda) \right] \geq a E \left[ \Lambda P_\tau(\infty) (\mbox{ab}; \Lambda) \right]
\]
Since $p \geq a$ by assumption, it is sufficient to show
\begin{equation} \label{Eq: sufficient_to_show}
E \left[ \Lambda \left(  P_\pi(\mbox{out};\Lambda) + P_\pi(\mbox{ab}; \Lambda) \right) \right] \geq E \left[ \Lambda P_\tau(\infty) (\mbox{ab}; \Lambda) \right].
\end{equation}
It can be seen through a coupling argument that the number of busy servers is stochastically larger in the system that does not outsource any calls, and so
\begin{equation} \label{Eq: needed_bound}
E\left[ \overline{S}_\pi(\Lambda) \right] \leq E\left[ \overline{S}_{\tau(\infty)}(\Lambda) \right].
\end{equation}
The inequality (\ref{Eq: needed_bound}) combined with the equality (\ref{Eq: Balance}) in the proof of part (i) implies (\ref{Eq: sufficient_to_show}), and so the proof is complete.
\eProof

\noindent{\bf Proof of Proposition~\ref{prop:fluid_bound}:}  We first prove part (i) and then prove part (ii). \\
\noindent{\newline \bf Proof of part (i):}
Since $E\left[ z_{\boldsymbol{\pi}}^\lambda(\mathbf{N}, \Lambda^\lambda(X)) \right] \geq 0$ for all $\lambda$,
\[
\liminf_{\lambda \rightarrow \infty} \frac{cN^\lambda + E \left[ z_{\boldsymbol{\pi}}(\mathbf{N}, \Lambda^\lambda(X)) \right]}{\lambda} \geq c \liminf_{\lambda \rightarrow \infty} \frac{N^\lambda}{\lambda}.
\]
Hence if $\liminf_{\lambda \rightarrow \infty} N^\lambda/ \lambda \geq 1$, the proof is complete.  Assume otherwise, that $\liminf_{\lambda \rightarrow \infty} N^\lambda/ \lambda = \delta < 1$.  As in the proof of Proposition 1, define $\overline{S}_{\boldsymbol{\pi}}^\lambda(l^\lambda)$ as the expected number of busy servers in the system with mean arrival rate $\lambda$ and arrival rate realization $l^\lambda = \lambda + x \sqrt{\lambda}$, and let $\overline{S}_{\boldsymbol{\pi}}^\lambda(\Lambda^\lambda)$ be the associated random variable.  Recall the equality (\ref{Eq: Balance}) in the proof of Proposition 1 and note that it holds for each $\lambda$, and so,
\begin{equation} \label{Eq: lower_bound_Prop2_proof}
E\left[ \Lambda^\lambda\left( P_{\boldsymbol{\pi}}^\lambda(\mbox{out}; \Lambda^\lambda) + P_{\boldsymbol{\pi}}^\lambda(\mbox{ab};\Lambda^\lambda) \right) \right] = \lambda - E\left[ \overline{S}_{\boldsymbol{\pi}}^\lambda \left( \Lambda^\lambda \right) \right] \geq \lambda - N^\lambda.
\end{equation}
(In words, (\ref{Eq: lower_bound_Prop2_proof}) states that the expected steady-state rate at which customers abandon or are outsourced must equal or exceed the difference between the mean arrival rate and the service capacity.)  From (\ref{Eq: lower_bound_Prop2_proof}) and our assumption that $\min(a,p) >c$,
\begin{eqnarray*}
E\left[ z_{\boldsymbol{\pi}}^\lambda \left( N^\lambda, \Lambda^\lambda(X) \right) \right] & \geq & \min(a,p) E \left[ \Lambda^\lambda\left(  P_{\boldsymbol{\pi}}^\lambda(\mbox{out}; \Lambda^\lambda) + P_{\boldsymbol{\pi}}^\lambda(\mbox{ab};\Lambda^\lambda)  \right) \right] \\
& \geq & c\left( \lambda - N^\lambda \right).
\end{eqnarray*}
The above inequality then implies
\[
\liminf_{\lambda \rightarrow \infty} \frac{cN^\lambda + E \left[ z_{\boldsymbol{\pi}}^\lambda \left( \mathbf{N}, \Lambda^\lambda(X) \right) \right]}{\lambda} \geq c \delta + c(1-\delta) = c.
\]

\noindent{\newline \bf Proof of part (ii):}  When no customers are outsourced, for every realization $l^\lambda = l^\lambda(x) = \lambda + \sqrt{\lambda}x$ of $\Lambda^\lambda$, the system operates as a $M/M/N^\lambda+M$ queue, which was analyzed in~\citeECtex{GarnettMandelbaumReiman02EC}.  Since the staffing assumption implies
\[
\sqrt{N^\lambda} \left( 1-\frac{l^\lambda(x)}{N^\lambda} \right) \rightarrow \beta-x \mbox{ as } \lambda \rightarrow \infty,
\]
the condition of Theorem 4 in~\citeECtex{GarnettMandelbaumReiman02EC}  is satisfied, and so\footnote{Although $\triangle$ is explicitly specified in~\citeECtex{GarnettMandelbaumReiman02EC}, we do not provide the expression here, because knowledge that the limit is finite is enough for our purposes.}
\begin{equation} \label{Eq:GMR}
\sqrt{N^\lambda} P_{\boldsymbol{\tau(\infty)}}^\lambda (\mbox{ab}; l^\lambda) \rightarrow \triangle \in (0,\infty) \mbox{ as } \lambda \rightarrow \infty.
\end{equation}
Since $P_{\tau(\infty)}^\lambda(\mbox{out};l^\lambda)=0$,
\begin{equation} \label{Eq:zbound}
z_{\boldsymbol{\tau(\infty)}}^\lambda(\mathbf{N},l^\lambda) = a l^\lambda P_{\boldsymbol{\tau(\infty)}}^\lambda(\mbox{ab};l^\lambda).
\end{equation}
It follows from (\ref{Eq:GMR}) and (\ref{Eq:zbound}) that
\[
\frac{z_{\boldsymbol{\tau(\infty)}}^\lambda(N^\lambda,l^\lambda)}{\lambda} \rightarrow 0 \mbox{ as } \lambda \rightarrow \infty,
\]
and so
\[
\frac{z_{\boldsymbol{\tau(\infty)}}^\lambda(N^\lambda,\Lambda^\lambda)}{\lambda} \rightarrow 0 \mbox{ almost surely, as } \lambda \rightarrow \infty.
\]
Since $P_{\boldsymbol{\tau(\infty)}}^\lambda(\mbox{ab};l^\lambda) \leq 1$ for any realization $l^\lambda$ of $\Lambda^\lambda$, (\ref{Eq:zbound}) implies
\[
\frac{z_{\boldsymbol{\tau(\infty)}}^\lambda(N^\lambda,\Lambda^\lambda)}{\lambda}  \leq a \left( 1 + \frac{1}{\sqrt{\lambda}}X \right) \leq a (1+|X|)
\]
for all $\lambda \geq 1$.  Since $E|X|<\infty$ by assumption, the dominated convergence theorem implies that
\[
\frac{E\left[ z_{\boldsymbol{\tau(\infty)}}^\lambda\left( \mathbf{N}, \Lambda^\lambda \right) \right]}{\lambda} \rightarrow 0 \mbox{ as } \lambda \rightarrow \infty,
\]
which is sufficient to complete the proof.
\eProof

\noindent{\bf Proof of Lemma~\ref{lemma:KW2010}:}
This is Theorem 5.2 from~\citeECtex{KocagaWard10EC} adapted to the setting of this paper, and, to read this proof, the reader is advised to have a copy of that paper on hand. We abbreviate ~\citeECtex{KocagaWard10EC} to KW for the remainder of this proof. First note that Theorem 5.2 in KW holds for any threshold admission policy $\theta^N$ defined in Theorem 5.1 in that paper, and not only for the policy $\theta^{\star,N}$ that appears in Theorem 5.2.  Next, when $X$ in this paper realizes as $x$, the arrival rate in KW is $l^\lambda = \lambda = x \sqrt{\lambda}$, so that
\[
l^\lambda(x) - N^\lambda = -m\sqrt{\lambda} + o\left( \sqrt{\lambda} \right) \mbox{ for } m = \beta-x.
\]
The above display implies that the conditions in Theorem 5.2 in KW are satisfied, and so, noting that $\lim_{t\rightarrow \infty} E[\xi(t,\theta^{\star,N})]/t$ in their notation is $z_{\boldsymbol{\tau}}(\mathbf{N},l^\lambda(x))$ in ours, for $\boldsymbol{\tau} = \{\tau(T^\lambda): \lambda \geq 0\}$ and
\[
T^\lambda = N^\lambda + \hat{T} \sqrt{l^\lambda(x)},
\]
it follows that
\begin{equation} \label{eq:kappaFromKW}
\frac{1}{\sqrt{\lambda}} z_{\boldsymbol{\tau}}(\mathbf{N},l^\lambda(x)) \rightarrow \kappa, \mbox{ as } \lambda \rightarrow \infty,
\end{equation}
where $\kappa$ is as defined in (4.9) in KW.

Finally, we match the notation and show the algebra that establishes $\kappa$ in (4.9) in KW is exactly $\hat{z}(m,\hat{T})$.  Substituting the notation in Table~\ref{tab:notationMatch} into $\kappa$ in (4.9) in KW shows
\[
\kappa = \frac{\left( \begin{array}{l} p\exp\left( \frac{-\gamma}{2} \left( \hat{T}^2+2\frac{m}{\gamma}\hat{T} \right) \right) \\
+ a \left[ 1-\exp\left( \frac{-\gamma}{2} \left( \hat{T}^2+2\frac{m}{\gamma}\hat{T} \right)  \right)
+ \sqrt{\frac{2\pi}{\gamma}} m \exp\left( \frac{m^2}{2\gamma} \right)\left( \Phi\left( \frac{m}{\sqrt{\gamma}} \right) - \Phi\left( \sqrt{\gamma}\left( \hat{T} + \frac{m}{\gamma} \right) \right) \right) \right] \end{array} \right) }{\sqrt{2\pi}\left[ \exp\left( \frac{m^2}{2} \right) \Phi(m) + \frac{1}{\sqrt{\gamma}}\exp\left( \frac{m^2}{2\gamma} \right) \left( \Phi\left( \sqrt{\gamma}\left( \hat{T} + \frac{m}{\gamma} \right) \right) - \Phi\left( \frac{m}{\sqrt{\gamma}} \right) \right) \right]}
\]
Then, multiplying by $\phi(m/\sqrt{\gamma}) / \phi(m/\sqrt{\gamma})$ and recalling that $\phi(x) = (1/\sqrt{2\pi}) \exp(-x^2/2)$ for $x \in (-\infty,\infty)$ yields
\[
\kappa = \frac{p \phi\left( \sqrt{\gamma}\left( \hat{T}+\frac{m}{\gamma} \right) \right)
+ a\left[ \phi \left( \frac{m}{\sqrt{\gamma}}  \right) - \phi\left( \sqrt{\gamma}\left( \hat{T}+\frac{m}{\gamma} \right) \right)
 + \frac{m}{\sqrt{\gamma}} \left[ \Phi \left( \frac{m}{\sqrt{\gamma}}  \right) - \Phi\left( \sqrt{\gamma}\left( \hat{T}+\frac{m}{\gamma} \right) \right) \right]\right] }{\frac{\phi\left( \frac{m}{\sqrt{\gamma}} \right)}{\phi(m)}\Phi(m) + \frac{1}{\sqrt{\gamma}}\left(\Phi\left( \sqrt{\gamma}\left( \hat{T}+\frac{m}{\gamma} \right) \right) - \Phi \left( \frac{m}{\sqrt{\gamma}}  \right) \right)},
\]
which is exactly $\hat{z}(m,\hat{T})$.

\begin{table}[h*]
  \centering
  \begin{tabular}{cc}
    \hline \hline
    {\bf KW notation} & {\bf This paper's notation} \\ \hline
  $\sigma^2$ & 2 \\
  $\mu$ & 1 \\
  $m$ & $\beta - x$ \\
  $h_I$ & 0 \\
  $a$ & $a$ \\
  $c$ & $p$ \\
  $l$ & $\hat{T}$ \\
\hline \hline
\end{tabular}
\caption[The notation match]{The notation match between~\citeECtex{KocagaWard10EC} and this paper.}
%\caption{The notation match between~Kocaga and Ward (2010) and this paper.}
\label{tab:notationMatch}
\end{table}

\eProof

\noindent{\bf Proof of Corollary~\ref{corollary:KW2010}:}
This follows because a careful reading of the proof of Theorem 5.2 in KW shows that that result holds whenever the term multiplying $\hat{T}$ in (\ref{eq:KW2010lemmaT}) is of the same order as the square root of the arrival rate.
\eProof

\noindent{\bf Proof of Theorem \ref{Thm: Cost Convergence}:  }
It follows from Lemma~\ref{lemma:KW2010}  that for a given realization $x$ of the random variable $X$
\[
\sqrt{\lambda} \left(\frac{cN^\lambda + z_{\boldsymbol{\tau}}^\lambda(\mathbf{N},x)}{\lambda} - c \right) \rightarrow c\beta+\hat{z}(\beta-x,\hat{T}) \mbox{ as } \lambda \rightarrow \infty.
\]
For this proof, we must argue that the interchange of limit and expectation is valid; in particular, it is enough to show that
\begin{equation} \label{eq:toShow}
\frac{1}{\sqrt{\lambda}} E\left[ z_{\boldsymbol{\tau}}^\lambda(\mathbf{N}, \Lambda^\lambda(X)) \right] \rightarrow E \left[ \hat{z}(\beta-X,\hat{T}) \right], \mbox{ as } \lambda \rightarrow \infty.
\end{equation}
Recall that
\[
z_{\boldsymbol{\tau}}^\lambda(\mathbf{N},\Lambda^\lambda(X)) = p \Lambda^\lambda P_{\boldsymbol{\tau}}^\lambda(\mbox{out}; \Lambda^\lambda(X)) + a \gamma \overline{Q}_{\boldsymbol{\tau}}^\lambda(\Lambda^\lambda(X)).
\]
Suppose we can show
\begin{equation} \label{eq:toShow1}
\frac{1}{\sqrt{\lambda}} E \left[ \Lambda^\lambda(X) P_{\boldsymbol{\tau}}^\lambda(\mbox{out},\Lambda^\lambda(X)) \right] \rightarrow EL_1(X) \mbox{ as } \lambda \rightarrow \infty
\end{equation}
and
\begin{equation} \label{eq:toShow2}
 \frac{1}{\sqrt{\lambda}}\gamma E \left[ \overline{Q}_{\boldsymbol{\tau}}^\lambda(\Lambda^\lambda(X)) \right] \rightarrow EL_2(X) \mbox{ as } \lambda \rightarrow \infty,
\end{equation}
where, for any $x \in (-\infty,\infty)$,
\begin{eqnarray*}
L_1(x) & := & \frac{ \phi\left( \sqrt{\gamma}\left( \hat{T}+\frac{m}{\gamma} \right) \right) }{\frac{\phi\left( \frac{m}{\sqrt{\gamma}} \right)}{\phi(m)}\Phi(m) + \frac{1}{\sqrt{\gamma}}\left(\Phi\left( \sqrt{\gamma}\left( \hat{T}+\frac{m}{\gamma} \right) \right) - \Phi \left( \frac{m}{\sqrt{\gamma}}  \right) \right)} \\[.1in]
L_2(x) & := & \frac{  \left[ \phi \left( \frac{m}{\sqrt{\gamma}}  \right) - \phi\left( \sqrt{\gamma}\left( \hat{T}+\frac{m}{\gamma} \right) \right)
 + \frac{m}{\sqrt{\gamma}} \left[ \Phi \left( \frac{m}{\sqrt{\gamma}}  \right) - \Phi\left( \sqrt{\gamma}\left( \hat{T}+\frac{m}{\gamma} \right) \right) \right]\right] }{\frac{\phi\left( \frac{m}{\sqrt{\gamma}} \right)}{\phi(m)}\Phi(m) + \frac{1}{\sqrt{\gamma}}\left(\Phi\left( \sqrt{\gamma}\left( \hat{T}+\frac{m}{\gamma} \right) \right) - \Phi \left( \frac{m}{\sqrt{\gamma}}  \right) \right)}.
\end{eqnarray*}
Since it is straightforward to see that $L_1(x)+L_2(x) = z(\beta-x,\hat{T})$, which implies $E[L_1(X)]+E[L_2(X)] = E\left[ \hat{z}(\beta-X,\hat{T}) \right]$, to complete the proof, it is enough to show (\ref{eq:toShow1}) and (\ref{eq:toShow2}).

\noindent {\bf The argument to show (\ref{eq:toShow1}):}

We begin by observing that for any realization $x$ of $X$, an argument similar to the proof of part (c) of Theorem 5.2 (on page 311 of \citeECtex{KocagaWard10EC}) shows that
\[
\frac{l^\lambda(x)}{\sqrt{\lambda}} P_{\boldsymbol{\tau}}^\lambda(\mbox{out},l^\lambda(x)) \rightarrow L_1(x),
\]
as $\lambda \rightarrow \infty$, almost surely.
Then, to show (\ref{eq:toShow1}), it is sufficient to show that $\Lambda^\lambda(X) P_{\boldsymbol{\tau}}^\lambda(\mbox{out},\Lambda^\lambda(X))/\sqrt{\lambda}$ is bounded by an integrable random variable. For this, first note that
\[
P_{\boldsymbol{\tau}}^\lambda(\mbox{out},\Lambda^\lambda(X)) \leq B(N^\lambda,\Lambda^\lambda(X)),
\]
where $B(m,\lambda)$ is the Erlang blocking probability in a $M$/$M$/$m$/$m$ model that has offered load $\lambda$ (which is exactly equal to the arrival rate in our model since $\mu=1$).  Also,
\[
B(N^\lambda,\Lambda^\lambda(X)) \leq B(N^\lambda,\Lambda^\lambda(|X|)),
\]
because the Erlang blocking probability is increasing in the offered load.
Furthermore, the Erlang loss function $L(m,\lambda):=\lambda B(m,\lambda)$ (again, when the offered load equals the arrival rate) is subadditive by Theorem 1 in \citeECtex{SmithWhitt1981}, and so
\[
L(N^\lambda, \Lambda^\lambda(X) \leq L(N^\lambda,\lambda) + L(N^\lambda, \sqrt{\lambda}|X|),
\]
which implies
\begin{equation} \label{eq:ErlangLossBound}
\Lambda^\lambda(X) B(N^\lambda, \Lambda^\lambda(X)) \leq \lambda B(N^\lambda,\lambda) + \sqrt{\lambda} |X| B(N^\lambda, \sqrt{\lambda}|X|).
\end{equation}
Then, also noting that $B(m,\lambda) \leq 1$ for any positive integer $m$ and finite $\lambda \geq0$, we have that
\begin{equation} \label{eq:ErlangLossBound2}
\frac{\Lambda^\lambda(X)}{\sqrt{\lambda}} P_{\boldsymbol{\tau}}^\lambda(\mbox{out};\Lambda^\lambda(X)) \leq \sqrt{\lambda}B(N^\lambda, \lambda) + |X|.
\end{equation}
Since (see for example, \citeECtex{Whitt84})
\[
\sqrt{\lambda} B(N^\lambda, \lambda) \rightarrow \frac{\phi(\beta)}{\Phi(\beta)} \mbox{ as } \lambda \rightarrow \infty,
\]
and recalling that $E|X|<\infty$ by assumption, it follows that the right-hand-side of (\ref{eq:ErlangLossBound2}) is bounded by an integrable random variable, so that (\ref{eq:toShow1}) is justified.

\noindent {\bf The argument so show (\ref{eq:toShow2}):}

We begin by observing that for any realization $x$ of $X$, an argument similar to the proof of part (b) of Theorem 5.2 in \citeECtex{KocagaWard10EC} shows that
\[
\frac{\gamma \overline{Q}_{\tau(T^\lambda)}(l^\lambda(x))}{\sqrt{\lambda}} \rightarrow L_2(x),
\]
as $\lambda \rightarrow \infty$, almost surely.
Then, to show (\ref{eq:toShow2}), it is sufficient to show that $\gamma \overline{Q}^{\lambda}_{\boldsymbol{\tau}}(\Lambda^\lambda(X))/\sqrt{\lambda}$ is bounded by an integrable random
variable.  Since $\gamma \overline{Q}_{\boldsymbol{\tau}}^\lambda(l^\lambda(x)) = l^\lambda(x)
 P_{\boldsymbol{\tau}}^\lambda(\mbox{ab};l^\lambda(x))$ for any realization $x$ of $X$ (as observed in Remark~\ref{remark:waitingCost}), it is sufficient to show that $\Lambda^\lambda(X) P_{\boldsymbol{\tau}}^\lambda(\mbox{ab};l^\lambda(X))/\sqrt{\lambda}$ is bounded by an integrable random variable.  It follows from a coupling argument that for any realization $x$ of $X$,
\[
P_{\boldsymbol{\tau}}^\lambda(\mbox{ab};l^\lambda(X)) \leq B(N^\lambda, l^\lambda(X)),
\]
where $B$ is the Erlang blocking probability defined in the argument to show (\ref{eq:toShow1}).  Therefore, the bound in (\ref{eq:ErlangLossBound}) implies
\[
\frac{\Lambda^\lambda(X)}{\sqrt{\lambda}}P_{\boldsymbol{\tau}}^\lambda(\mbox{ab};l^\lambda(X)) \leq \sqrt{\lambda}B(N^\lambda, \lambda) + |X|.
\]
The right-hand side of the above expression is bounded by an integrable random variable for the exact same reason that (\ref{eq:ErlangLossBound2}) was.

\eProof

\noindent{\bf Proof of Proposition~\ref{prop: beta_cost_finite}:  }
We first observe that
\begin{eqnarray} \label{eq:expBound}
\lefteqn{E \left[ \hat{z}(\beta-X,\hat{T}^\star(\beta-X)) \right] \leq} \\ & & E\left[ \hat{z}(\beta-X,\hat{T}^\star(\beta-X)) | X \leq \beta+1 \right] + E \left[ \hat{z}(\beta-X,\hat{T}^\star(\beta-X)) | X > \beta +1\right]. \nonumber
\end{eqnarray}
Next, we use the following two claims. The proofs of the claims can be found immediately following this proof.

\begin{claim} \label{claim:xincreasing}
For any fixed $\beta \in (-\infty,\infty)$, $\hat{z}(\beta-x,\hat{T}^\star(\beta-x))$ is an increasing function of $x$ on $(-\infty,\infty)$.
\end{claim}

\begin{claim} \label{claim:betasupfinite}
For any fixed $\beta \in (-\infty,\infty)$ and $x > \beta+1$,
\[
 \hat{z}(\beta-x,0) < p(x-\beta+1).
\]
\end{claim}
\noindent It follows from Claim~\ref{claim:xincreasing} that
\begin{equation} \label{eq:boundbeta1}
E \left[ \hat{z}(\beta-X,\hat{T}^\star(\beta-X)) | X \leq \beta+1 \right] \leq \hat{z}(-1,\hat{T}^\star(-1)),
\end{equation}
and it follows from Claim~\ref{claim:betasupfinite} that
\begin{equation} \label{eq:boundbeta2}
 E \left[ \hat{z}(\beta-X,\hat{T}^\star(\beta-X)) | X > \beta +1\right] \leq  p \left( E|X| + \beta +1 \right).
\end{equation}
Since $\hat{z}(-1,\hat{T}^\star(-1)) < \infty$ from the definition of $\hat{z}$ in (\ref{eq:zdef}) and $\hat{T}^\star$ in (\ref{Eq: T def}) (recalling that there exists a unique $\hat{T}^\star$ that solves (\ref{Eq: T def}) by Proposition 4.1 in~\citeECtex{KocagaWard10EC}) and $E|X| < \infty$ by assumption, we conclude from (\ref{eq:expBound}), (\ref{eq:boundbeta1}), and (\ref{eq:boundbeta2}) that
\[
E\left[ \hat{z}\left( \beta - X, \hat{T}^\star(\beta-X) \right) \right] \leq \hat{z}(-1,\hat{T}^\star(-1)) + p\left( E|X| + \beta +1 \right) < \infty.
\]

\eProof

\noindent{\bf Proof of Claim~\ref{claim:xincreasing}:  }
It follows from a coupling argument that the expected steady-state number of customers waiting in a $M$/$M$/$N$/$B$+$M$ queue, as well as the expected steady-state loss proportion, is increasing as the arrival rate increases, but $N$ and $B$ remain fixed.  Next, recall from Section~\ref{Section:Exact} that for any fixed realization $x$ of $X$, the long-run average operating cost associated with an admissible routing policy $\boldsymbol{\pi}$ in the system with mean arrival rate $\lambda$ is
\[
z_{\boldsymbol{\pi}}^\lambda(\mathbf{N}, l^\lambda(x)) = p l^\lambda(x) P_{\boldsymbol{\pi}}^\lambda(\mbox{out}; l^\lambda(x)) + a\gamma \overline{Q}_{\boldsymbol{\pi}}^\lambda(l^\lambda(x)).
\]
Since $l^\lambda(x) = \lambda + x \sqrt{\lambda}$ is increasing in $x$, it follows that if $x_1 < x_2$, then
\begin{equation} \label{eq:zboundarrival}
z_{\boldsymbol{\pi}}^\lambda(\mathbf{N}, l^\lambda(x_1)) \leq z_{\boldsymbol{\pi}}^\lambda(\mathbf{N}, l^\lambda(x_2)).
\end{equation}
Suppose $(N^\lambda - \lambda)/\sqrt{\lambda} \rightarrow \beta$ as $\lambda \rightarrow \infty$, define $\hat{T}_2^\star$ to satisfy (\ref{Eq: T def}) with $x$ replaced by $x_2$, and let $\boldsymbol{\tau} = \{ \tau(T^\lambda): \lambda \geq 0\}$ where
\[
T^\lambda = N^\lambda + \hat{T}_2^\star\sqrt{\lambda}.
\]
Then, Corollary~\ref{corollary:KW2010} implies that
\begin{equation} \label{eq:limitComparison}
\lim_{\lambda \rightarrow \infty} \frac{1}{\sqrt{\lambda}} z_{\boldsymbol{\tau}}^\lambda(\mathbf{N}, l^\lambda(x_1)) = \hat{z}(\beta-x_1, \hat{T}_2^\star)
\mbox{ and }
\lim_{\lambda \rightarrow \infty} \frac{1}{\sqrt{\lambda}} z_{\boldsymbol{\tau}}^\lambda(\mathbf{N}, l^\lambda(x_2)) = \hat{z}(\beta-x_2,\hat{T}_2^\star).
\end{equation}
It follows from (\ref{eq:zboundarrival}) and (\ref{eq:limitComparison}) that
\[
\hat{z}(\beta-x,\hat{T}_2^\star) \leq \hat{z}(\beta-x_2,\hat{T}_2^\star).
\]
Next, from (\ref{Eq: Dif Cost LB}), since $\hat{T}_1^\star$ is the minimizer of $\hat{z}(\beta-x_1,\hat{T})$ over $\hat{T} \in [0,\infty)$,
\[
\hat{z}(\beta-x_1, \hat{T}_1^\star) \leq \hat{z}(\beta-x_1, \hat{T}_2^\star).
\]
We conclude from the previous two displays that
\[
\hat{z}(\beta-x_1, \hat{T}_1^\star) \leq \hat{z}(\beta-x_2, \hat{T}_2^\star),
\]
which completes the proof, since $x_1$ and $x_2$ are arbitrary.
\eProof

\noindent{\bf Proof of Claim~\ref{claim:betasupfinite}: }
We use the following inequality, that is given in 7.1.13 of \citeECtex{AbramowitzStegun64},
\begin{equation} \label{eq:AbSt}
\frac{1}{x+\sqrt{x^2+2}} < e^{x^2} \int_x^\infty e^{-t^2} dt \mbox{ for } x \geq 0.
\end{equation}
From the definition of $\hat{z}$ in (\ref{eq:zdef}), the symmetry of the normal distribution, and algebra,
\begin{eqnarray*}
\hat{z}(\beta-x,0) & = & p \frac{\phi(\beta-x)}{\Phi(\beta-x)} \\
& = & p \frac{\phi(x-\beta)}{1-\Phi(x-\beta)} \\
& = & \frac{p}{\sqrt{2}} \times \left[ \exp\left( \left(\frac{x-\beta}{\sqrt{2}} \right)^2 \right)  \int_{\frac{x-\beta}{\sqrt{2}}}^\infty \exp(-t^2) dt \right]^{-1}.
\end{eqnarray*}
Then, the inequality (\ref{eq:AbSt}) implies that
\[
\hat{z}(\beta-x,0) < \frac{p}{\sqrt{2}} \left( \frac{x-\beta}{\sqrt{2}} + \sqrt{\left( \frac{x-\beta}{\sqrt{2}} \right)^2 +2} \right) \mbox{ for all } x \geq \beta +1.
\]
Finally, since
\[
\frac{x-\beta}{\sqrt{2}} + \sqrt{\left( \frac{x-\beta}{\sqrt{2}} \right)^2 +2}  \leq \frac{x-\beta}{\sqrt{2}} + \sqrt{\left( \frac{x-\beta}{\sqrt{2}} + \sqrt{2} \right)^2} = \sqrt{2}(x-\beta+1),
\]
the stated claim is established.
\eProof

\noindent{\bf Proof of Proposition~\ref{prop: beta_star_finite}:  }

In light of Proposition~\ref{prop: beta_cost_finite}, it is sufficient to show that
\begin{enumerate}
\item[(a)] $c\beta + E \left[ \hat{z}\left( \beta - X, \hat{T}^\star(\beta-X) \right) \right] \rightarrow \infty$ as $|\beta| \rightarrow \infty$; and,
\item[(b)] $c\beta + E \left[ \hat{z}\left( \beta - X, \hat{T}^\star(\beta-X) \right) \right]$ is continuous in $\beta$ on $(\beta_0,\infty)$ for $\beta_0<0$ and $|\beta_0|$ arbitrarily large.
\end{enumerate}
We first show (a) and then show (b).  To do this, the following two claims are useful.  Their proofs can be found immediately after the proof of this proposition.

\begin{claim} \label{claim:continuity}
Let $p<a$. Then, for any realization $x$ of $X$, $\hat{T}^\star$ defined in (\ref{Eq: T def}) is continuous in $\beta$.
\end{claim}

\begin{claim} \label{claim:nonincreasing}
For any realization $x$ of $X$, $\hat{z}(\beta-x, \hat{T}^\star(\beta-x))$ is a decreasing function of $\beta$ on $(-\infty,\infty)$.
\end{claim}

\noindent{\newline \bf Proof of (a):}  For any realization $x$ of $X$, any $\beta \in \Re$, and $N^\lambda$ that satisfies the conditions of Lemma 1, by Lemma 1, $\hat{z}(\beta-x,\hat{T}^\star(\beta-x))$ is obtained as the limit of the non-negative cost function $z_{\boldsymbol{\tau^\star}}^\lambda(\mathbf{N}, l^\lambda(x))$, scaled by $1/\sqrt{\lambda}$, for $\boldsymbol{\tau^\star} = \{ T^{\lambda,\star}: \lambda \geq 0\}$ and $T^{\lambda,\star}$ defined in (\ref{eq:l_value}).  Hence $\hat{z}(\beta-x, \hat{T}^\star(\beta-x))$ is a non-negative function of $\beta$, and so
\[
c\beta + E \left[ \hat{z}(\beta-X, \hat{T}^\star(\beta-X)) \right] \rightarrow \infty \mbox{ as } \beta \rightarrow \infty.
\]
Next, we handle the case that $\beta \rightarrow -\infty$. When $p<a$, recall from (\ref{Eq: T def}) that  $\hat{T}^\star(\beta-x)$ solves
\[
\hat{z}(\beta-x,\hat{T}^\star(\beta-x)) = -p(\beta-x) + (a-p) \gamma \hat{T}^\star(\beta-x),
\]
for any realization $x$ of $X$. Then, adding $c\beta$ to both sides, taking expectations, and recalling that $EX=0$,
\[
c\beta + E \left[ \hat{z}\left( \beta-X,\hat{T}^\star(\beta-X) \right) \right] = (c-p)\beta + (a-p) \gamma E\left[\hat{T}^\star(\beta-X)\right].
\]
Since $c-p <0$ by assumption and $\hat{T}^\star(\beta-x) \geq 0$ for any realization $x$, it follows that
\[
c\beta + E\left[ \hat{z}(\beta-X, \hat{T}^\star(\beta-X)) \right] \rightarrow \infty \mbox{ as } \beta \rightarrow -\infty,
\]
when $p<a$.

When $a\leq p$, we first observe that $T^\star(\beta -x)=\infty$ for any realization $x$ of $X$ and so
\[\hat{z }(\beta -x,T^\star(\beta -x))=\hat{z }(\beta -x,\infty).
\]
Next, we let $\hat{\underline{z}}(\beta -x,\hat{\underline{T}}^\star(\beta -x))$ be the minimum cost and $\hat{\underline{T}}^\star(\beta -x)$ the optimal threshold of the diffusion control problem with costs $\underline{a}$ and $\underline{p}$ such that $\underline{p}<\underline{a}=a\leq p$. It follows from the optimality of $\hat{\underline{T}}^\star$ that, for any realization $X=x$,

\[\hat{\underline{z}}(\beta -x,\hat{\underline{T}}^\star(\beta -x))\leq \hat{\underline{z}}(\beta -x,\infty)=\hat{z}(\beta -x,\infty),
\]
where the equality follows since $\hat{\underline{z}}(\beta -x,\infty)$ and $\hat{z}(\beta -x,\infty)$ do not depend on $p$ or $\underline{p}$. Hence
\[c\beta+E\left[\hat{\underline{z}}(\beta -X,\hat{\underline{T}}^\star(\beta -X))\right] \leq c\beta+E\left[\hat{z }(\beta -X,\infty)\right].
\]
Since $\underline{p}<\underline{a}=a$, we can also repeat the same arguments as in the case $p<a$ to get
\[
c\beta + E\left[ \hat{\underline{z}}(\beta-X, \hat{\underline{T}}^\star(\beta-X)) \right] \rightarrow \infty \mbox{ as } \beta \rightarrow -\infty.
\]
Combining the previous two displays we get
\[
c\beta + E\left[ \hat{z}(\beta-X, \hat{T}^\star(\beta-X)) \right] =c\beta + E\left[ \hat{z}(\beta-X, \infty) \right]\rightarrow \infty \mbox{ as } \beta \rightarrow -\infty,
\]
when $a\leq p$.

\noindent{\newline \bf Proof of (b):}  Set $\beta_0<0$ and $|\beta_0|$ arbitrarily large.  From the definition of continuity, it is enough to show that if $\{\beta_n\}$ is a sequence in $(\beta_0,\infty)$ that converges to $\beta \in (\beta_0,\infty)$ as $n \rightarrow \infty$, then
\begin{equation}\label{Eq:contiuityCond}
c \beta_n + E \left[ \hat{z}\left( \beta_n-X, \hat{T}^\star(\beta_n-X) \right) \right] \rightarrow c \beta + E \left[ \hat{z}\left( \beta-X, \hat{T}^\star(\beta-X) \right) \right] \mbox{ as } n \rightarrow \infty.
\end{equation}
Since $\hat{z}$ is a continuous function of its arguments (as can be seen immediately from its definition in (\ref{eq:zdef})) and $\hat{T}^\star$ is a continuous function of $\beta$ by Claim~\ref{claim:continuity}, it follows that for any realization $x$ of $X$,
\[
\hat{z}(\beta_n-x,\hat{T}^\star(\beta_n - x)) \rightarrow \hat{z}(\beta-x,\hat{T}^\star(\beta-x)) \mbox{ as } n \rightarrow \infty.
\]
For every $n$, by Claim~\ref{claim:nonincreasing},
\[
\hat{z}(\beta_n - x, \hat{T}^\star(\beta_n-x)) \leq \hat{z}(\beta_0-x, \hat{T}^\star(\beta_0-x)), \mbox{ for any } x \in (-\infty,\infty).
\]
Since by Proposition~\ref{prop: beta_cost_finite},
\[
E\left[ \hat{z}\left( \beta_0-X,\hat{T}^\star(\beta_0-X) \right) \right] < \infty,
\]
the dominated convergence theorem implies that
\[
E\left[ \hat{z}\left( \beta_n - X, \hat{T}^\star(\beta_n-X)\right) \right] \rightarrow E \left[ \hat{z}\left( \beta-X,\hat{T}^\star(\beta-X) \right) \right] \mbox{ as } n \rightarrow \infty,
\]
from which (\ref{Eq:contiuityCond}) follows.

\eProof

\noindent{\bf Proof of Claim~\ref{claim:continuity}:  }
Let $p<a$, and for a fixed realization of $X$ as $x$, define
\[
g(\beta,\hat{T}):= (a-p) \gamma \hat{T} - \hat{z}(\beta-x,\hat{T}) - p(\beta-x).
\]
For any fixed point $(\beta_1,\hat{T}^\star) \in \Re \times[0,\infty)$, (\ref{eq:zdef}) shows that $\hat{z}(\beta-x,\hat{T})$ is differentiable. Hence
\[
\left. \frac{\partial g}{\partial \hat{T}} \right\vert_{\hat{T}=\hat{T}^\star} = (a-p) \gamma - \left. \frac{d\hat{z}(\beta-x,\hat{T})}{d\hat{T}} \right\vert_{\hat{T}=\hat{T}^\star}.
\]
Recall from (\ref{Eq: T def}) and (\ref{Eq: Dif Cost LB}) that when $g(\beta_1,\hat{T}^\star)=0$, then
\[
\hat{z}(\beta_1 -x,\hat{T}^\star) \leq \hat{z}(\beta_1 - x,\hat{T}) \mbox{ for any other } \hat{T} \geq 0.
\]
Hence the first-order conditions imply that
\[
\left. \frac{d\hat{z}(\beta-x,\hat{T})}{d\hat{T}}\right\vert_{\hat{T}=\hat{T}^\star} = 0,
\]
and so
\[
\left. \frac{\partial g}{\partial \hat{T}} \right\vert_{\hat{T}=\hat{T}^\star} = (a-p) \gamma.
\]
Because $a<p$ by assumption, the conditions of the implicit function theorem are satisfied.  We use the implicit function theorem to conclude that there exists an open set $\mathcal{S}_{\beta_1}$ containing $\beta_1$ and an open set $\mathcal{S}_{\hat{T}^\star}$ containing $\hat{T}^\star$, and a unique continuously invertible function $h: \mathcal{S}_{\beta_1} \rightarrow \mathcal{S}_{\hat{T}^\star}$ such that
\[
\left\{ (\beta,h(\beta)) : \beta \in \mathcal{S}_{\beta_1} \right\} = \{ (\beta,\hat{T}^\star) \in \mathcal{S}_{\beta_1} \times \mathcal{S}_{\hat{T}^\star} \}.
\]
The existence of a unique continuously invertible function $h$ holds for any $\beta_1 \in \Re$, and so the proof is complete.
\eProof

\noindent{\bf Proof of Claim~\ref{claim:nonincreasing}:  }
As in the proof of Claim~\ref{claim:xincreasing}, it follows from a coupling argument that
 the expected steady-state number of customers waiting in a $M/M/N/B+M$ queue, as well as the expected steady-state loss proportion, is decreasing as $N$ increases but $B$ remains fixed.    Next, recall from Section 3 that for any fixed realization $x$ of $X$, the long-run average operating cost associated with an admissible routing policy $\boldsymbol{\pi}$ in the system with mean arrival rate $\lambda$ is
\[
z_{\boldsymbol{\pi}}^\lambda(\mathbf{N},l^\lambda(x)) = pl^\lambda(x) P_{\boldsymbol{\pi}}^\lambda(\mbox{out}; l^\lambda(x)) + a\gamma \overline{Q}_{\boldsymbol{\pi}}^\lambda(l^\lambda(x)).
\]
Therefore, for two staffing policies
\[
\mathbf{N_1} = \{  \lambda + \beta_1\sqrt{\lambda}: \lambda \geq 0\} \mbox{ and } \mathbf{N_2} = \{ \lambda + \beta_2 \sqrt{\lambda}: \lambda \geq 0 \} \mbox{ with } \beta_1 \leq \beta_2
\]
and the identical threshold routing policy (equivalently, the identical $B$ because the threshold is on the total number of customers in the system) $\boldsymbol{\tau} = \{T^\lambda: \lambda \geq 0 \}$ with
\[
T^\lambda = N_1^\lambda + \hat{T}_1^\star\sqrt{l^\lambda(x)},
\]
for $\hat{T}_1^\star$ that satisfies (\ref{Eq: T def}) with $\beta_1$ replacing $\beta$, it follows that
\[
z_{\boldsymbol{\tau}}^\lambda (\mathbf{N_2},l^\lambda(x)) < z_{\boldsymbol{\tau}}(\mathbf{N_1},l^\lambda(x)).
\]
Dividing both sides of the above inequality by $\sqrt{\lambda}$, taking the limit as $\lambda \rightarrow \infty$, and applying Claim~\ref{claim:continuity}, shows that
\begin{equation}\label{Eq:claim2LB1}
\hat{z}\left( \beta_2-x,\hat{T}_1^\star \right) < \hat{z}\left( \beta_1 - x, \hat{T}_1^\star \right).
\end{equation}
Since from (\ref{Eq: Dif Cost LB}), $\hat{T}_2^\star$ is a minimizer when $\beta_2$ replaces $\beta$ (not $\beta_1$),
\begin{equation}\label{Eq:claim2LB2}
\hat{z}\left( \beta_2-x,\hat{T}_2^\star \right) < \hat{z}\left( \beta_2-x,\hat{T}_1^\star \right).
\end{equation}
The proof is complete from (\ref{Eq:claim2LB1}), (\ref{Eq:claim2LB2}), and the fact that $\beta_1$ and $\beta_2$ are arbitrary.

\eProof

\noindent{\bf Proof of Theorem \ref{Thm: AO}:  }
It follows from Theorem~\ref{Thm: Cost Convergence} that
\[
\mathcal{\hat{C}}^\lambda(U) \rightarrow \hat{\mathcal{C}}^\star < \infty, \mbox{ as } \lambda \rightarrow \infty.
\]
Therefore, to show asymptotic optimality (see Definition~\ref{asymptotic optimality}), it is enough to show
\begin{equation} \label{eq:ao_toshow}
\liminf_{\lambda \rightarrow \infty} \hat{\mathcal{C}}^\lambda(\mathbf{u}) \geq \hat{\mathcal{C}}^\star
\end{equation}
under any arbitrary admissible policy
\[
\mathbf{u} = (\mathbf{N},\boldsymbol{\pi}) = \{ (N^\lambda, \pi^\lambda): \lambda \geq 0 \}.
\]
We first establish (\ref{eq:ao_toshow}) and then prove
\begin{equation} \label{eq:littleo}
\frac{cN^{\lambda,\star} + E \left[ z_{\boldsymbol{\tau}^\star}^\lambda\left( \mathbf{N^\star}, \Lambda^\lambda(X) \right) \right]-\mathcal{C}^{\lambda,\mbox{\tiny{opt}}}(\lambda) }{\sqrt{\lambda}} \rightarrow 0, \mbox{ as } \lambda \rightarrow \infty.
\end{equation}

{\bf The argument that our proposed policy is asymptotically optimal (that (\ref{eq:ao_toshow}) holds).}
We first argue that we need only consider admissible policies under which
\begin{equation} \label{eq:liminfinfinite}
\liminf_{\lambda \rightarrow \infty} \frac{N^\lambda - \lambda}{\sqrt{\lambda}} > - \infty
\end{equation}
holds.  To see this, assume the bound (established at the end of this proof)
\begin{equation} \label{eq:assumedbound}
E \left[ z_{\boldsymbol{\pi}}^\lambda(\mathbf{N}, \Lambda^\lambda(X)) \right] \geq \min(a,p)(\lambda - N^\lambda)
\end{equation}
is valid.  Then, from the definition of $\hat{C}^\lambda(\mathbf{u})$ and (\ref{eq:assumedbound}),
\begin{eqnarray*}
\hat{\mathcal{C}}^\lambda( \mathbf{u} ) & = & c \frac{N^\lambda - \lambda}{\sqrt{\lambda}} + \frac{E \left[ z_{\boldsymbol{\pi}}^\lambda(\mathbf{N}, \Lambda^\lambda(X)) \right]}{\sqrt{\lambda}} \\
& \geq & \left( \min(a,p) - c \right) \frac{\lambda - N^\lambda}{\sqrt{\lambda}}.
\end{eqnarray*}
If (\ref{eq:liminfinfinite}) does not hold, then it follows from the assumption $\min(a,p) >c$ and the above display that
\[
\liminf_{\lambda \rightarrow \infty} \hat{\mathcal{C}}^\lambda(\mathbf{u}) = \infty,
\]
which trivially satisfies (\ref{eq:ao_toshow}).  In summary, it is enough to show (\ref{eq:ao_toshow}) holds for the subset of admissible policies that satisfy (\ref{eq:liminfinfinite}).

Consider any subsequence $\lambda_i$ on which (\ref{eq:liminfinfinite}) holds and also on which the $\liminf$ in (\ref{eq:ao_toshow}) is attained, so that
\[
\lim_{\lambda_i \rightarrow \infty} \hat{C}^{\lambda_i}(\mathbf{u}) = \liminf_{\lambda \rightarrow \infty} \hat{C}^\lambda(\mathbf{u}) < \infty.
\]
We may assume the limit is finite because otherwise (\ref{eq:ao_toshow}) holds trivially.  On this subsequence, since the limit is finite and
\[
\hat{C}^{\lambda_i}(\mathbf{u}) = c \frac{N^{\lambda_i} - \lambda_i}{\sqrt{\lambda_i}} + \frac{E \left[ z_{\boldsymbol{\pi}}^\lambda(\mathbf{N}, \Lambda^\lambda(X)) \right]}{\sqrt{\lambda}}
\]
has its second term positive, it must be the case that
\[
\limsup_{\lambda_i \rightarrow \infty} \frac{N^{\lambda_i} - \lambda_i}{\sqrt{\lambda_i}} < \infty.
\]
Since
\[
-\infty < \liminf_{\lambda_i \rightarrow \infty} \frac{N^{\lambda_i} - \lambda_i}{\sqrt{\lambda_i}} < \limsup_{\lambda_i \rightarrow \infty} \frac{N^{\lambda_i} - \lambda_i}{\sqrt{\lambda_i}} < \infty,
\]
the Bolzano-Weierstrass theorem guarantees that any subsequence has a further convergent subsequence $\lambda_{i_j}$ on which
\begin{equation} \label{eq:BW}
\frac{N^{\lambda_{i_j}} - \lambda_{i_j}}{\sqrt{\lambda_{i_j}}} \rightarrow \beta \in \Re \mbox{ as } \lambda_{i_j} \rightarrow \infty.
\end{equation}
Since from the properties of the limit
\[
\lim_{\lambda_{i_j}\rightarrow \infty} \hat{C}(\mathbf{u}) \geq c \lim_{\lambda_{i_j} \rightarrow \infty} \frac{N^{\lambda_{i_j}} - \lambda_{i_j}}{\sqrt{\lambda_{i_j}}} + \liminf_{\lambda_{i_j} \rightarrow \infty} \frac{E \left[ z_{\boldsymbol{\pi}}^{\lambda_{i_j}}(\mathbf{N}, \Lambda^{\lambda_{i_j}}(X)) \right]}{\sqrt{\lambda_{i_j}}},
\]
and Fatou's lemma guarantees
\[
\liminf_{\lambda_{i_j} \rightarrow \infty} \frac{E \left[ z_{\boldsymbol{\pi}}^{\lambda_{i_j}}(\mathbf{N}, \Lambda^{\lambda_{i_j}}(X)) \right]}{\sqrt{\lambda_{i_j}}} \geq E \left[ \liminf_{\lambda_{i_j} \rightarrow \infty} \frac{z_{\boldsymbol{\pi}}^{\lambda_{i_j}}(\mathbf{N}, \Lambda^{\lambda_{i_j}}(X)) }{\sqrt{\lambda_{i_j}}} \right]
\]
it follows that
\begin{equation} \label{eq:boundFatou}
\lim_{\lambda_{i_j}\rightarrow \infty} \hat{C}(\mathbf{u}) \geq c \beta +  E \left[ \liminf_{\lambda_{i_j}  \rightarrow \infty} \frac{z_{\boldsymbol{\pi}}^{\lambda_{i_j}}(\mathbf{N}, \Lambda^{\lambda_{i_j}}(X)) }{\sqrt{\lambda_{i_j}}} \right].
\end{equation}
Next, for any realization $x$ of $X$, it is straightforward to see Theorem 5.2 part (ii) in~\citeECtex{KocagaWard10EC} can be used to conclude
\begin{equation} \label{eq:boundKW}
\liminf_{\lambda_{i_j}  \rightarrow \infty} \frac{z_{\boldsymbol{\pi}}^{\lambda_{i_j}}(\mathbf{N}, \l^{\lambda_{i_j}}(x)) }{\sqrt{\lambda_{i_j}}} \geq \hat{z} \left( \beta - x, \hat{T}^\star(\beta-x) \right),
\end{equation}
because (\ref{eq:BW}) implies the conditions of that theorem are satisfied since\footnote{Please see Table~\ref{tab:notationMatch} to see how to match the notation between~\citeECtex{KocagaWard10EC} and this paper.}
\[
\frac{N^{\lambda_{i_j}}-l^{\lambda_{i_j}}(x)}{\sqrt{\lambda_{i_j}}} = \frac{N^{\lambda_{i_j}}-\lambda_{i_j}}{\sqrt{\lambda_{i_j}}} - x \rightarrow \beta - x, \mbox{ as } \lambda_{i_j} \rightarrow \infty.
\]
It follows from (\ref{eq:boundFatou}) and (\ref{eq:boundKW}) that
\[
\lim_{\lambda_{i_j} \rightarrow \infty} \hat{C}^{\lambda_{i_j}}(\mathbf{u}) \geq c \beta + E \left[\hat{z}\left( \beta-X, \hat{T}^\star(\beta-X)\right) \right].
\]
Since the definitions of $\beta^\star$ and $\hat{\mathcal{C}}^\star$ imply that
\[
c \beta + E \left[\hat{z}\left( \beta-X, \hat{T}^\star(\beta-X)\right) \right] \geq \hat{\mathcal{C}}^\star = c \beta^\star + E \left[\hat{z}\left( \beta^\star-X, \hat{T}^\star(\beta^\star-X)\right) \right],
\]
we conclude that (\ref{eq:ao_toshow}) is satisfied.  Since the subsequence $\lambda_{i_j}$ was arbitrary, the proof is complete once we establish the earlier assumed bound (\ref{eq:assumedbound}).

To see the bound (\ref{eq:assumedbound}) holds, let $\overline{S}^\lambda_{\boldsymbol{\pi}}(l^\lambda(x))$ denote the expected number of busy servers when the realized arrival rate is $l^\lambda(x) = \lambda + x \sqrt{\lambda}$.  Since the arrival rate into the system must equal the departure rate from the system (due to both abandonments and service completions),
\[
l^\lambda(x) \left( 1-P^\lambda_{\mathbf{\pi}}(\mbox{out}; l^\lambda(x)) \right) = \overline{S}_{\boldsymbol{\pi}}^\lambda(l^\lambda(x)) + l^\lambda(x) P_{\boldsymbol{\pi}}^\lambda(\mbox{ab};l^\lambda(x)),
\]
or, equivalently,
\[
l^\lambda(x) \left( P^\lambda_{\mathbf{\pi}}(\mbox{out}; l^\lambda(x)) + P_{\boldsymbol{\pi}}^\lambda(\mbox{ab};l^\lambda(x))  \right) = l^\lambda(x) - \overline{S}_{\boldsymbol{\pi}}^\lambda(l^\lambda(x)).
\]
Since $\overline{S}_{\boldsymbol{\pi}}^\lambda(l^\lambda(x)) \leq N^\lambda$, it follows that
\[
l^\lambda(x) \left( P^\lambda_{\mathbf{\pi}}(\mbox{out}; l^\lambda(x)) + P_{\boldsymbol{\pi}}^\lambda(\mbox{ab};l^\lambda(x))  \right) \geq l^\lambda(x) - N^\lambda.
\]
From the definition of $z_{\boldsymbol{\pi}}^\lambda(\mathbf{N},l^\lambda(x))$, it follows that
\[
z_{\boldsymbol{\pi}}^\lambda(\mathbf{N},l^\lambda(x)) \geq \min(a,p) l^\lambda(x)  \left( P^\lambda_{\mathbf{\pi}}(\mbox{out}; l^\lambda(x)) + P_{\boldsymbol{\pi}}^\lambda(\mbox{ab};l^\lambda(x))  \right).
\]
Hence
\[
z_{\boldsymbol{\pi}}^\lambda(\mathbf{N},l^\lambda(x)) \geq \min(a,p) \left( l^\lambda(x) - N^\lambda \right),
\]
and (\ref{eq:assumedbound}) follows by taking expectations and recalling $EX=0$.

{\bf The argument that our proposed policy achieves cost $o(\sqrt{\lambda})$ higher than the minimum cost (that (\ref{eq:littleo}) holds).}

Recall that $\mathcal{C}^{\lambda,\mbox{\tiny{opt}}}$ is the minimum achievable cost defined in (\ref{eq:objective}) for the system with mean arrival rate $\lambda$.  Since
\[
\frac{cN^{\lambda,\star} + E \left[ z_{\boldsymbol{\tau}^\star}^\lambda\left( \mathbf{N^\star}, \Lambda^\lambda(X) \right) \right]-\mathcal{C}^{\lambda,\mbox{\tiny{opt}}} }{\sqrt{\lambda}}  = \hat{\mathcal{C}}^\lambda(\mathbf{u}^\star) - \frac{\mathcal{C}^{\lambda,\mbox{\tiny{opt}}} - c}{\sqrt{\lambda}},
\]
and from Theorem~\ref{Thm: Cost Convergence}
\[
\hat{\mathcal{C}}^\lambda(\mathbf{u}^\star) \rightarrow \hat{\mathcal{C}}^\star \mbox{ as } \lambda \rightarrow \infty,
\]
it is enough to establish
\begin{equation} \label{eq:cor1}
\frac{\mathcal{C}^{\lambda,\mbox{\tiny{opt}}} - c}{\sqrt{\lambda}} \rightarrow \hat{\mathcal{C}}^\star \mbox{ as } \lambda \rightarrow \infty.
\end{equation}
It follows from Theorem~\ref{Thm: AO} that
\begin{equation} \label{eq:cor2}
\liminf_{\lambda \rightarrow \infty} \frac{\mathcal{C}^{\lambda,\mbox{\tiny{opt}}} - c}{\sqrt{\lambda}}  \geq \hat{\mathcal{C}}^\star.
\end{equation}
Also, a policy $\mathbf{u^{\mbox{\tiny{opt}}}}$ that consists of a sequence $(N^{\lambda,\mbox{\tiny{opt}}}, \pi^{\lambda,\mbox{\tiny{opt}}})$ in which each element of the sequence is exactly optimal for each $\lambda$, and so achieves the minimum cost $\mathcal{C}^{\mbox{\tiny{opt}}}$, must be asymptotically optimal.  Hence, because $\mathbf{u}^{\mbox{\tiny{opt}}}$ is an admissible policy, from the definition of asymptotic optimality,
\begin{equation} \label{eq:cor3}
\limsup_{\lambda \rightarrow \infty} \frac{\mathcal{C}^{\lambda,\mbox{\tiny{opt}}} - c}{\sqrt{\lambda}} = \limsup_{\lambda \rightarrow \infty} \hat{\mathcal{C}}^\lambda(\mathbf{u^{\mbox{\tiny{opt}}}}) \leq \lim_{\lambda \rightarrow \infty} \hat{\mathcal{C}}^\lambda(\mathbf{u^\star}) = \hat{\mathcal{C}}^\star.
\end{equation}
The limit (\ref{eq:cor1}) follows from (\ref{eq:cor2}) and (\ref{eq:cor3}).
\eProof

\section{Supporting Numerical Tables}
\label{Section:AdditionalNumerics}

In this section we provide detailed numerical results in tabular format which support our findings in Section \ref{Section: Model}. We set the mean service time and the mean patience time equal to 1 and fix the cost parameters at $c=0.1$, $p=1$ and $a=5$ unless specified otherwise.

Table \ref{tab: CV1} provides the details of the numerical study shown in Figure \ref{fig: T2}. The first column shows the Uniform arrival rate distribution which has its mean $\lambda$ fixed at 100 and has an increasing CV as we go down in the table. The second column shows the optimal staffing level. Columns three for and five show the the staffing level, associated cost and cost percentage error of $U$, respectively. Columns six, seven, and eight report the same numbers for the first alternative policy, $D$, while columns nine, ten and eleven report the same numbers for the second alternative policy $NV$.

Tables \ref{tab: T2R1CR}, \ref{tab: T2R5CR}  and \ref{tab: T2R9CR} provide details for the numerical studies associated with Figures \ref{fig:crlowcv_cost} and \ref{fig:crlowcv_staff}, Figures \ref{fig:crmedcv_cost} and \ref{fig:crmedcv_staff}, and Figures \ref{fig:crhighcv_cost} and \ref{fig:crhighcv_staff}, respectively. The first column shows the varying staffing cost while the other columns are the same as in Table \ref{tab: CV1}. Table \ref{tab: T2R159CR beta} provides the $\beta^{\star}$ values used in Tables \ref{tab: T2R1CR}-\ref{tab: T2R9CR}.

Tables \ref{tab: T2R1Skew}, \ref{tab: T2R5Skew} and \ref{tab: T2R9Skew} provide details for the numerical studies associated with Figures \ref{fig:skewlowcv_cost} and \ref{fig:skewlowcv_staff}, Figures \ref{fig:skewmedcv_cost} and \ref{fig:skewmedcv_staff}, and Figures \ref{fig:skewhighcv_cost} and \ref{fig:skewhighcv_staff}, respectively. The first column shows the parameters of the Beta distribution that $\Lambda$ follows while the other columns are the same as before.

\begin{table}[htb*]
    \begin{tabular}{c||c|ccc|ccc|ccc}
    \hline
    \multicolumn{1}{c|}{Arrival rate}  &\multicolumn{1}{c|}{\multirow{2}{*}{$N^{\mbox{\tiny{opt}}}$}}& \multicolumn{3}{c|}{$U$}  & \multicolumn{3}{c|}{$D$}  & \multicolumn{3}{c}{NV}\\ \cline{3-11}
    \multicolumn{1}{c|}{distribution}  &\multicolumn{1}{c|}{} &$N_{U}$  & $\mathcal{C}(N_{U})$ &\% error&  $N_{D}$ & $\mathcal{C}(N_{D})$ &\% error&  $N_{NV}$ & $\mathcal{C}(N_{NV})$ &\% error \\ \hline\hline
	$\Lambda=100$ & 119   & 119   & 12.41 & 0.01\% & 119   & 12.41 & 0.00\% & 101   & 17.32 & 39.65\% \\
	$\mathcal{U}[99,101]$ & 119   & 119   & 12.41 & 0.01\% & 119   & 12.41 & 0.00\% & 101   & 16.82 & 35.59\% \\
    $\mathcal{U}[90,110]$ & 121   & 121   & 12.71 & 0.01\% & 119   & 12.76 & 0.36\% & 108   & 14.73 & 15.90\% \\
    $\mathcal{U}[80,120]$ & 127   & 126   & 13.38 & 0.03\% & 119   & 13.75 & 2.76\% & 116   & 14.12 & 5.57\% \\
    $\mathcal{U}[70,130]$ & 133   & 132   & 14.16 & 0.07\% & 119   & 15.19 & 7.33\% & 124   & 14.56 & 2.94\% \\
    $\mathcal{U}[60,140]$ & 140   & 139   & 14.97 & 0.05\% & 119   & 16.93 & 13.14\% & 132   & 15.24 & 1.79\% \\
    $\mathcal{U}[50,150]$ & 147   & 146   & 15.82 & 0.06\% & 119   & 18.88 & 19.38\% & 140   & 16.00 & 1.19\% \\
    $\mathcal{U}[40,160]$ & 155   & 154   & 16.67 & 0.03\% & 119   & 20.95 & 25.65\% & 148   & 16.81 & 0.84\% \\
    $\mathcal{U}[30,170]$ & 162   & 161   & 17.54 & 0.05\% & 119   & 23.10 & 31.73\% & 156   & 17.65 & 0.62\% \\
    $\mathcal{U}[20,180]$ & 170   & 169   & 18.42 & 0.04\% & 119   & 25.32 & 37.53\% & 164   & 18.50 & 0.47\% \\
    $\mathcal{U}[10,190]$ & 178   & 176   & 19.30 & 0.06\% & 119   & 27.59 & 43.02\% & 172   & 19.36 & 0.37\% \\
    \hline\hline
\end{tabular}
    \caption{Performance of $U$ vs other policies: Varying CV (Figure~\ref{fig: T2} in the main body)}
\label{tab: CV1}
\end{table}

\begin{table}[htb*]
  \centering
  \begin{tabular}{c||c|ccc|ccc|ccc}
    \hline
    \multicolumn{1}{c|}{Staffing cost}  &\multicolumn{1}{c|}{\multirow{2}{*}{$N^{\mbox{\tiny{opt}}}$}}& \multicolumn{3}{c|}{$U$}  & \multicolumn{3}{c|}{$D$}  & \multicolumn{3}{c}{NV}\\ \cline{3-11}
    \multicolumn{1}{c|}{($c$)}  &\multicolumn{1}{c|}{} &$N_{U}$  & $\mathcal{C}(N_{U})$ &\% error&  $N_{D}$ & $\mathcal{C}(N_{D})$ &\% error&  $N_{NV}$ & $\mathcal{C}(N_{NV})$ &\% error \\ \hline\hline

    0.01  & 134   & 132   & 1.38  & 0.31\% & 129   & 1.41  & 2.68\% & 110   & 4.07  & 195.70\% \\
    0.05  & 126   & 125   & 6.54  & 0.06\% & 122   & 6.61  & 1.18\% & 109   & 8.77  & 34.24\% \\
    0.1   & 121   & 121   & 12.71 & 0.01\% & 119   & 12.76 & 0.36\% & 108   & 14.73 & 15.90\% \\
    0.2   & 116   & 116   & 24.56 & 0.02\% & 115   & 24.57 & 0.03\% & 106   & 25.75 & 4.83\% \\
    0.3   & 112   & 112   & 35.93 & 0.02\% & 112   & 35.93 & 0.00\% & 104   & 36.70 & 2.13\% \\
    0.4   & 108   & 108   & 46.91 & 0.03\% & 108   & 46.91 & 0.00\% & 102   & 47.35 & 0.95\% \\
    0.5   & 104   & 105   & 57.51 & 0.08\% & 105   & 57.51 & 0.02\% & 100   & 57.70 & 0.36\% \\
    0.6   & 100   & 101   & 67.72 & 0.08\% & 102   & 67.75 & 0.07\% & 98    & 67.75 & 0.07\% \\
    0.7   & 95    & 96    & 77.48 & 0.05\% & 97    & 77.51 & 0.03\% & 96    & 77.48 & 0.00\% \\
    0.8   & 89    & 90    & 86.71 & 0.07\% & 91    & 86.74 & 0.04\% & 94    & 86.90 & 0.23\% \\
    0.9   & 75    & 78    & 94.99 & 0.35\% & 79    & 95.01 & 0.05\% & 92    & 95.98 & 1.07\% \\
    0.95  & 59    & 63    & 98.38 & 0.75\% & 64    & 98.39 & 0.03\% & 91    & 100.39 & 2.05\% \\
    0.99  & 1     & 8     & 101.01 & 0.01\% & 11    & 100.01 & 0.01\% & 90    & 103.81 & 3.81\% \\

\hline \hline
\end{tabular}
\caption{Performance of $U$ vs other policies: Varying staffing cost ($c$) under low arrival rate variability (Figure~\ref{fig:crlowcv_cost} (a) and~\ref{fig:crlowcv_staff} (a) in the main body)}
\label{tab: T2R1CR}
\end{table}

\begin{table}[htb*]
  \centering
  \begin{tabular}{c||c|ccc|ccc|ccc}
    \hline
    \multicolumn{1}{c|}{Staffing cost}  &\multicolumn{1}{c|}{\multirow{2}{*}{$N^{\mbox{\tiny{opt}}}$}}& \multicolumn{3}{c|}{$U$}  & \multicolumn{3}{c|}{$D$}  & \multicolumn{3}{c}{NV}\\ \cline{3-11}
    \multicolumn{1}{c|}{($c$)}  &\multicolumn{1}{c|}{} &$N_{U}$  & $\mathcal{C}(N_{U})$ &\% error&  $N_{D}$ & $\mathcal{C}(N_{D})$ &\% error&  $N_{NV}$ & $\mathcal{C}(N_{NV})$ &\% error \\ \hline\hline
    0.01  & 170   & 165   & 1.77  & 1.15\% & 129   & 5.44  & 210.53\% & 149   & 2.41  & 37.62\% \\
    0.05  & 156   & 154   & 8.24  & 0.18\% & 122   & 12.14 & 47.53\% & 145   & 8.58  & 4.31\% \\
    0.1   & 147   & 146   & 15.82 & 0.06\% & 119   & 18.88 & 19.38\% & 140   & 16.00 & 1.19\% \\
    0.2   & 134   & 134   & 29.85 & 0.04\% & 115   & 31.35 & 5.04\% & 130   & 29.91 & 0.20\% \\
    0.3   & 122   & 123   & 42.64 & 0.04\% & 112   & 43.08 & 1.04\% & 120   & 42.66 & 0.04\% \\
    0.4   & 111   & 112   & 54.28 & 0.10\% & 108   & 54.29 & 0.03\% & 110   & 54.28 & 0.00\% \\
    0.5   & 100   & 102   & 64.83 & 0.21\% & 105   & 64.93 & 0.19\% & 100   & 64.81 & 0.00\% \\
    0.6   & 89    & 91    & 74.28 & 0.35\% & 102   & 75.03 & 1.03\% & 90    & 74.27 & 0.00\% \\
    0.7   & 79    & 81    & 82.69 & 0.37\% & 97    & 84.25 & 1.92\% & 80    & 82.67 & 0.01\% \\
    0.8   & 68    & 71    & 90.05 & 0.31\% & 91    & 92.48 & 2.74\% & 70    & 90.03 & 0.02\% \\
    0.9   & 56    & 58    & 96.27 & 0.21\% & 79    & 98.46 & 2.29\% & 60    & 96.32 & 0.06\% \\
    0.95  & 44    & 45    & 98.81 & 0.14\% & 64    & 99.69 & 0.90\% & 55    & 99.01 & 0.21\% \\
    0.99  & 0     & 0     & 100.00 & 0.00\% & 11    & 100.03 & 0.03\% & 51    & 100.92 & 0.92\% \\
\hline \hline
\end{tabular}
\caption{Performance of $U$ vs other policies: Varying staffing cost ($c$) under moderate arrival rate variability (Figure~\ref{fig:crlowcv_cost} (b) and~\ref{fig:crlowcv_staff} (b) in the main body)}
\label{tab: T2R5CR}
\end{table}

\begin{table}[htb*]
  \centering
  \begin{tabular}{c||c|ccc|ccc|ccc}
    \hline
    \multicolumn{1}{c|}{Staffing cost}  &\multicolumn{1}{c|}{\multirow{2}{*}{$N^{\mbox{\tiny{opt}}}$}}& \multicolumn{3}{c|}{$U$}  & \multicolumn{3}{c|}{$D$}  & \multicolumn{3}{c}{NV}\\ \cline{3-11}
    \multicolumn{1}{c|}{($c$)}  &\multicolumn{1}{c|}{} &$N_{U}$  & $\mathcal{C}(N_{U})$ &\% error&  $N_{D}$ & $\mathcal{C}(N_{D})$ &\% error&  $N_{NV}$ & $\mathcal{C}(N_{NV})$ &\% error \\ \hline\hline
    0.01  & 209   & 202   & 2.18  & 1.49\% & 129   & 13.33 & 519.46\% & 188   & 2.58  & 20.02\% \\
    0.05  & 191   & 188   & 10.10 & 0.16\% & 122   & 26.74 & 165.06\% & 181   & 10.26 & 1.74\% \\
    0.1   & 178   & 176   & 19.30 & 0.06\% & 119   & 27.59 & 43.02\% & 172   & 19.36 & 0.37\% \\
    0.2   & 156   & 156   & 35.97 & 0.03\% & 115   & 40.30 & 12.05\% & 154   & 35.98 & 0.04\% \\
    0.3   & 136   & 138   & 50.60 & 0.09\% & 112   & 52.16 & 3.09\% & 136   & 50.60 & 0.00\% \\
    0.4   & 117   & 119   & 63.29 & 0.23\% & 108   & 63.52 & 0.37\% & 118   & 63.29 & 0.00\% \\
    0.5   & 99    & 101   & 74.10 & 0.29\% & 105   & 74.19 & 0.15\% & 100   & 74.09 & 0.01\% \\
    0.6   & 80    & 83    & 83.03 & 0.22\% & 102   & 84.31 & 1.57\% & 82    & 83.02 & 0.01\% \\
    0.7   & 61    & 65    & 90.12 & 0.17\% & 97    & 93.49 & 3.78\% & 64    & 90.10 & 0.02\% \\
    0.8   & 43    & 47    & 95.35 & 0.13\% & 91    & 101.54 & 6.54\% & 46    & 95.33 & 0.03\% \\
    0.9   & 24    & 28    & 98.71 & 0.08\% & 79    & 106.71 & 8.14\% & 28    & 98.71 & 0.03\% \\
    0.95  & 15    & 15    & 99.67 & 0.02\% & 64    & 106.10 & 6.45\% & 19    & 99.71 & 0.04\% \\
    0.99  & 0     & 0     & 100.00 & 0.00\% & 11    & 100.14 & 0.14\% & 12    & 100.17 & 0.17\% \\
\hline \hline
\end{tabular}
\caption{Performance of $U$ vs other policies: Varying staffing cost ($c$) under high arrival rate variability (Figure~\ref{fig:crlowcv_cost} (c) and~\ref{fig:crlowcv_staff} (c) in the main body)}
\label{tab: T2R9CR}
\end{table}

\begin{table}[htb*]
  \centering
  \begin{tabular}{c||c|c|c|}
    \hline
    \multicolumn{1}{c|}{Staffing cost}  & \multicolumn{3}{c|}{Arrival Rate Variability}  \\ \cline{2-4}
    \multicolumn{1}{c|}{($c$)}  & Low (Table \ref{tab: T2R1CR})& Moderate (Table \ref{tab: T2R5CR}) &High (Table \ref{tab: T2R9CR})\\ \hline\hline
    0.01  & 3.2164 & 6.5123 & 10.1808 \\
    0.05  & 2.5108 & 5.4114 & 8.7650 \\
    0.10  & 2.1109 & 4.6235 & 7.6149 \\
    0.20  & 1.5948 & 3.3824 & 5.6329 \\
    0.30  & 1.1972 & 2.2735 & 3.7603 \\
    0.40  & 0.8368 & 1.2118 & 1.9217 \\
    0.50  & 0.4777 & 0.1723 & 0.0980 \\
    0.60  & 0.0881 & -0.8550 & -1.7180 \\
    0.70  & -0.3778 & -1.8761 & -3.5296 \\
    0.80  & -1.0220 & -2.9188 & -5.3385 \\
    0.90  & -2.2158 & -4.2349 & -7.2004 \\
    0.95  & -3.6768 & -5.5266 & -8.5063 \\
    0.99  & -9.2008 & -10.3327 & -12.5916 \\
\hline \hline
\end{tabular}
\caption{$\beta^{\star}$ values for changing staffing costs and arrival rate variabilities given in Tables \ref{tab: T2R1CR}-\ref{tab: T2R9CR}}
\label{tab: T2R159CR beta}
\end{table}

\begin{table}[htb*]
  \centering
  \begin{tabular}{c||c|ccc|ccc|ccc}
    \hline
    \multicolumn{1}{c|}{Arrival rate}  &\multicolumn{1}{c|}{\multirow{2}{*}{$N^{\mbox{\tiny{opt}}}$}}& \multicolumn{3}{c|}{$U$}  & \multicolumn{3}{c|}{$D$}  & \multicolumn{3}{c}{NV}\\ \cline{3-11}
    \multicolumn{1}{c|}{distribution}  &\multicolumn{1}{c|}{} &$N_{U}$  & $\mathcal{C}(N_{U})$ &\% error&  $N_{D}$ & $\mathcal{C}(N_{D})$ &\% error&  $N_{NV}$ & $\mathcal{C}(N_{NV})$ &\% error \\ \hline\hline
    $Beta(1.5,0.5)$ & 121   & 121   & 12.65 & 0.00\% & 119   & 12.70 & 0.34\% & 106   & 15.17 & 19.90\% \\
    $Beta(1.4,0.6)$ & 121   & 121   & 12.67 & 0.01\% & 119   & 12.71 & 0.35\% & 106   & 15.17 & 19.74\% \\
    $Beta(1.3,0.7)$ & 121   & 121   & 12.68 & 0.01\% & 119   & 12.72 & 0.35\% & 107   & 14.83 & 16.92\% \\
    $Beta(1.2,0.8)$ & 121   & 121   & 12.69 & 0.01\% & 119   & 12.74 & 0.36\% & 107   & 14.82 & 16.79\% \\
    $Beta(1.1,0.9)$ & 121   & 121   & 12.70 & 0.02\% & 119   & 12.75 & 0.36\% & 108   & 14.51 & 14.22\% \\
    $Beta(1.0,1.0)$ & 121   & 121   & 12.71 & 0.01\% & 119   & 12.76 & 0.36\% & 108   & 14.51 & 14.10\% \\
    $Beta(0.9,1.1)$ & 122   & 121   & 12.72 & 0.01\% & 119   & 12.77 & 0.36\% & 108   & 14.50 & 13.98\% \\
    $Beta(0.8,1.2)$ & 122   & 121   & 12.74 & 0.01\% & 119   & 12.78 & 0.35\% & 109   & 14.22 & 11.65\% \\
    $Beta(0.7,1.3)$ & 121   & 121   & 12.75 & 0.01\% & 119   & 12.79 & 0.34\% & 109   & 14.21 & 11.52\% \\
    $Beta(0.6,1.4)$ & 121   & 121   & 12.76 & 0.01\% & 119   & 12.80 & 0.33\% & 109   & 14.21 & 11.37\% \\
    $Beta(0.5,1.5)$ & 121   & 121   & 12.77 & 0.02\% & 119   & 12.81 & 0.31\% & 109   & 14.20 & 11.20\% \\
\hline \hline
\end{tabular}
\caption{Performance of $U$ vs other policies: Varying skewness under low arrival rate variability (Figure~\ref{fig: skew_cost} (a) and~\ref{fig: skew_staff} (a) in the main body)}
\label{tab: T2R1Skew}
\end{table}

\begin{table}[htb*]
  \centering
  \begin{tabular}{c||c|ccc|ccc|ccc}
    \hline
    \multicolumn{1}{c|}{Arrival rate}  &\multicolumn{1}{c|}{\multirow{2}{*}{$N^{\mbox{\tiny{opt}}}$}}& \multicolumn{3}{c|}{$U$}  & \multicolumn{3}{c|}{$D$}  & \multicolumn{3}{c}{NV}\\ \cline{3-11}
    \multicolumn{1}{c|}{distribution}  &\multicolumn{1}{c|}{} &$N_{U}$  & $\mathcal{C}(N_{U})$ &\% error&  $N_{D}$ & $\mathcal{C}(N_{D})$ &\% error&  $N_{NV}$ & $\mathcal{C}(N_{NV})$ &\% error \\ \hline\hline
    $Beta(1.5,0.5)$ & 140   & 139   & 14.73 & 0.06\% & 119   & 17.53 & 19.03\% & 128   & 15.62 & 6.10\% \\
    $Beta(1.4,0.6)$ & 142   & 140   & 14.97 & 0.14\% & 119   & 17.90 & 19.70\% & 131   & 15.59 & 4.28\% \\
    $Beta(1.3,0.7)$ & 144   & 142   & 15.18 & 0.07\% & 119   & 18.20 & 20.01\% & 134   & 15.60 & 2.81\% \\
    $Beta(1.2,0.8)$ & 145   & 144   & 15.38 & 0.05\% & 119   & 18.47 & 20.06\% & 136   & 15.72 & 2.20\% \\
    $Beta(1.1,0.9)$ & 146   & 145   & 15.60 & 0.05\% & 119   & 18.69 & 19.83\% & 138   & 15.85 & 1.66\% \\
    $Beta(1.0,1.0)$ & 147   & 146   & 15.82 & 0.06\% & 119   & 19.11 & 20.84\% & 140   & 16.00 & 1.19\% \\
    $Beta(0.9,1.1)$ & 149   & 147   & 16.04 & 0.06\% & 119   & 19.04 & 18.74\% & 142   & 16.16 & 0.81\% \\
    $Beta(0.8,1.2)$ & 150   & 149   & 16.27 & 0.04\% & 119   & 19.18 & 17.88\% & 143   & 16.38 & 0.71\% \\
    $Beta(0.7,1.3)$ & 150   & 149   & 16.52 & 0.06\% & 119   & 19.29 & 16.80\% & 144   & 16.61 & 0.60\% \\
    $Beta(0.6,1.4)$ & 151   & 150   & 16.78 & 0.04\% & 119   & 19.37 & 15.46\% & 145   & 16.85 & 0.46\% \\
    $Beta(0.5,1.5)$ & 151   & 151   & 17.06 & 0.03\% & 119   & 19.41 & 13.77\% & 146   & 17.11 & 0.31\% \\
\hline \hline

\end{tabular}
\caption{Performance of $U$ vs other policies: Varying skewness under moderate arrival rate variability (Figure~\ref{fig: skew_cost} (b) and~\ref{fig: skew_staff} (b) in the main body)}
\label{tab: T2R5Skew}
\end{table}

\begin{table}[htb*]
  \centering
  \begin{tabular}{c||c|ccc|ccc|ccc}
    \hline
    \multicolumn{1}{c|}{Arrival rate}  &\multicolumn{1}{c|}{\multirow{2}{*}{$N^{\mbox{\tiny{opt}}}$}}& \multicolumn{3}{c|}{$U$}  & \multicolumn{3}{c|}{$D$}  & \multicolumn{3}{c}{NV}\\ \cline{3-11}
    \multicolumn{1}{c|}{distribution}  &\multicolumn{1}{c|}{} &$N_{U}$  & $\mathcal{C}(N_{U})$ &\% error&  $N_{D}$ & $\mathcal{C}(N_{D})$ &\% error&  $N_{NV}$ & $\mathcal{C}(N_{NV})$ &\% error \\ \hline\hline
    Beta(1.1,0.9) & 175   & 173   & 18.85 & 0.07\% & 119   & 27.31 & 44.98\% & 169   & 18.93 & 0.48\% \\
    $Beta(1.0,1.0)$ & 178   & 176   & 19.30 & 0.06\% & 119   & 27.59 & 43.02\% & 172   & 19.36 & 0.37\% \\
    $Beta(0.9,1.1)$ & 180   & 179   & 19.76 & 0.05\% & 119   & 27.81 & 40.77\% & 175   & 19.81 & 0.26\% \\
    $Beta(0.8,1.2)$ & 182   & 181   & 20.24 & 0.04\% & 119   & 27.98 & 38.22\% & 178   & 20.27 & 0.16\% \\
    $Beta(0.7,1.3)$ & 184   & 183   & 20.75 & 0.04\% & 119   & 28.08 & 35.33\% & 180   & 20.77 & 0.13\% \\
    $Beta(0.6,1.4)$ & 186   & 185   & 21.29 & 0.03\% & 119   & 28.11 & 32.05\% & 182   & 21.30 & 0.09\% \\
    $Beta(0.5,1.5)$ & 187   & 186   & 21.87 & 0.04\% & 119   & 28.04 & 28.26\% & 183   & 21.88 & 0.06\% \\
\hline \hline
\end{tabular}
\caption{Performance of $U$ vs other policies: Varying skewness under high arrival rate variability (Figure~\ref{fig: skew_cost} (c) and~\ref{fig: skew_staff} (c) in the main body)}
\label{tab: T2R9Skew}
\end{table}

\clearpage

\renewcommand{\refname}{Electronic Compantion References}
\bibliographystyleECtex{ormsv080}
\bibliographyECtex{AKW_04_07_2014_EC}

\end{document}